\documentclass[11pt]{amsart}
\pdfoutput=1

\usepackage[usenames,dvipsnames]{xcolor}   %

\usepackage[utf8]{inputenc}
\usepackage[T1]{fontenc}

\usepackage[letterpaper,margin=1in]{geometry}

\usepackage{lmodern}
\usepackage{microtype}
\usepackage{etoolbox}                 %
\usepackage{lettrine}

\usepackage{mathtools}
\usepackage{amsmath,amssymb,amsfonts,amsxtra,bm}
\usepackage{amsthm}
\usepackage{mathrsfs}
\usepackage{latexsym}
\usepackage{nicefrac}       %

\usepackage{booktabs}
\usepackage{makecell,multirow}       %
\usepackage{colortbl}
\usepackage{graphicx}
\usepackage{epstopdf}
\usepackage{float}
\floatstyle{plaintop}
\restylefloat{table}
\usepackage{footnote, tablefootnote}

\usepackage{enumitem}
\usepackage{algorithm}
\usepackage{algorithmic}
\usepackage{comment}
\usepackage{xspace}
\usepackage{url}
\usepackage[most]{tcolorbox}

\usepackage[round]{natbib}

\usepackage{hyperref}

\numberwithin{equation}{section}

\makeatletter

\definecolor{darkblue}{rgb}{0.0,0.0,0.65}
\definecolor{darkred}{rgb}{0.68,0.05,0.0}
\definecolor{darkgreen}{rgb}{0.0,0.29,0.29}
\definecolor{darkpurple}{rgb}{0.47,0.09,0.29}

\patchcmd{\section}{\scshape}{\fontfamily{ppl}\bfseries\scshape\selectfont}{}{%
  \PackageWarning{aiproofs}{Could not bold the section font; amsart may have changed}}

\newcommand{\aiauthor}{%
  \author[Suvrit Sra]{Suvrit Sra\textsuperscript{\normalfont*}}%
  \thanks{*\,TU Munich, Dept.~of Mathematics, School of CIT; Garching, Germany.
    Email: \texttt{s.sra@tum.de}.}}

\newcommand{\arxivlink}[1]{%
  \thanks{Also on arXiv: \url{https://arxiv.org/abs/#1}.}}
\newcommand{\preprintsfooter}[2]{%
  \def\ps@firstpage{\ps@plain
    \def\@oddfoot{\hfil\parbox[t]{\textwidth}{\centering\normalfont\scriptsize
      First public release: #2, as v1 in the author's preprint repository,
      \url{https://github.com/suvrit/preprints/tree/main/preprints/#1};
      every version is kept there with its timestamp.\\[3pt]\thepage}\hfil
      \global\topskip\normaltopskip}%
    \let\@evenfoot\@oddfoot
    \def\@oddhead{\@serieslogo\hss}%
    \let\@evenhead\@oddhead}}

\makeatother

\hypersetup{colorlinks=true,linkcolor=red!60!black,citecolor=green!48!black,urlcolor=blue!60!black}
\setlist{nosep,leftmargin=*}

\newenvironment{aiusage}{\section*{AI usage}}{}

\newcommand{\aiacknowledgements}[1][]{%
  \section*{Acknowledgements}
  \if\relax\detokenize{#1}\relax\else #1\ \fi
  The author (SS) thanks the Alexander von Humboldt (AvH) Foundation for
  their generous support via a Humboldt Professorship in AI.%
}

\newcommand{\reals}{\mathbb{R}}
\newcommand{\N}{\mathbb{N}}

\newcommand{\E}{\mathbb{E}}

\newcommand{\C}{\mathbb{C}}

\newcommand{\nlsum}{\sum\nolimits}

\DeclareMathOperator{\real}{Re}
\DeclareMathOperator{\imag}{Im}

\DeclareMathOperator{\supp}{supp}

\newtheorem{theorem}{Theorem}[section]

\newtheorem{prop}[theorem]{Proposition}
\newtheorem{lemma}[theorem]{Lemma}
\newtheorem{corr}[theorem]{Corollary}

\theoremstyle{definition}

\definecolor{cdarkred}{RGB}{180,70,30} %
\definecolor{cdarkgreen}{RGB}{0,130,0} %
\definecolor{azure}{RGB}{240,255,255}
\definecolor{cmidgreen}{RGB}{201,219,216}
\definecolor{clightgreen}{RGB}{200,255,200}
\definecolor{clightgray}{RGB}{240,240,240}
\definecolor{cmygreen}{RGB}{10,120,0}
\definecolor{cmyblue}{RGB}{50,50,255}
\definecolor{sblue}{RGB}{100,100,255}
\definecolor{cmypurple}{RGB}{150,50,180}
\definecolor{cmyred}{RGB}{238, 64,  0}
\definecolor{cdarkblue}{RGB}{30,75,170}
\definecolor{csharpred}{RGB}{255,0,0}

\tcbset{dashedrounded/.style={no borderline,borderline={0.5pt}{0pt}{dashed},arc=4pt}}
\tcbset{supporting/.style={colback=gray!8!white}}
\newtcolorbox{dashedbox}[1][]{
  before={\par\medskip\noindent},
  after={\par\medskip},
  colback=yellow!6!white,
  boxrule=0pt,
  enhanced jigsaw,
  top=2pt,
  left=2pt,
  right=2pt,
  bottom=0pt,
  dashedrounded,
  #1
}

\newtcbox{\tcbmathbox}[1][]{
  nobeforeafter,
  colback=yellow!30!white,
  boxrule=0pt,
  enhanced jigsaw,
  top=0pt,
  left=1pt,
  right=1pt,
  bottom=0pt,
  borderline horizontal={0.5pt}{0pt}{dashed},
  borderline vertical={0.5pt}{0pt}{dashed},
  #1
}

\usepackage{needspace}
\usepackage{tabularx}
\newcolumntype{Y}{>{\raggedright\arraybackslash}X}
\allowdisplaybreaks[2]
\setlist[enumerate]{itemsep=3pt,topsep=6pt}

\newcommand{\Zp}{\mathcal Z_+}
\newcommand{\dd}{\,\mathrm d}
\newcommand{\PF}{\mathrm{PF}}
\newcommand{\CBF}{\mathrm{CBF}}
\newcommand{\Tstar}{T_*}
\DeclareMathOperator{\Exp}{Exp}
\DeclareMathOperator{\GammaLaw}{Gamma}\DeclareMathOperator{\Lap}{Lap}

\hypersetup{
  pdftitle={On some positivity properties of the Riemann Xi function},
  pdfauthor={Suvrit Sra},
  bookmarksdepth=2
}

\begin{document}
\title[On some positivity properties of the Riemann $\Xi$ function]{On some positivity properties of the Riemann $\Xi$ function}
\aiauthor
\begin{abstract}
  We prove strict log-concavity of the positive functions $\Theta$ and $\Lambda$ with one-sided and
  bilateral Laplace transforms $1/\xi(\tfrac12+\sqrt s)$ and $1/\Xi$, respectively, where
  $\Xi(z)=\xi(\tfrac12+iz)$. For either function, log-concavity is the P\'olya frequency property of
  order two (PF$_2$). %
  Without assuming the Riemann hypothesis (RH), but using its published verification up to height
  about $3\cdot10^{12}$, we further prove that both densities are $\PF_{k_*}$ for $k_*=1\,810\,540\,931\,075$.
  The proof views $\Theta$ as a small signed perturbation of a density that is PF of every order,
  with an error bound independent of the order, so the numerical check reduces to scalar
  inequalities. Since $\Lambda$ is a Gaussian mixture of $\Theta$, it inherits the same order.
  Conversely, a pole obstruction shows that if either density is PF$_k$, then every zeta zero
  lower than the $(k+2)$-nd lies on the critical line. Finally,
  we prove that if RH holds up to height $T>\tfrac12$, then $\Xi$ satisfies the generalized Laguerre
  inequalities $L_n[\Xi](x)\ge0$ of every order $n$ for $|x|\le T-\tfrac12$, extending the case $n=1$
  of Csordas, Ruttan and Varga (1991).
\end{abstract}
\keywords{Riemann Xi; P\'olya frequency; log-concavity; generalized Laguerre; complete monotonicity.}
\maketitle

\addtocontents{toc}{\protect\setcounter{tocdepth}{-1}}

\vspace*{-8pt}
\section{Introduction}\label{sec:headline}
For a connoisseur of positivity, the Riemann $\Xi$ function offers a tantalizing array of questions. While the hardest of these are equivalent to the Riemann hypothesis (RH) itself, we find that even weak necessary conditions for RH already repay study. The ideas permeating our study are classical: Fourier--Laplace transforms, the Laguerre--P\'olya class, P\'olya frequency functions, among others. Our goal is a richer understanding of these positivity properties without assuming RH, and ideally with as little information about the zeta zeros as possible.

\vskip4pt
\noindent We use the usual normalization
\citep[Section~2.1]{titchmarsh1986theory} for $\Xi$:
\begin{equation}\label{eq:normalization}
 \xi(s)=\tfrac12s(s-1)\pi^{-s/2}\Gamma(s/2)\zeta(s),\qquad
 c_0=\xi(\tfrac12)>0,\qquad \Xi(z)=\xi(\tfrac12+iz),
\end{equation}
and define the following functions
\begin{equation}\label{eq:FhB}
 F(z) :=\Xi(z)/c_0,\qquad
 h(z) :=\xi(\tfrac12+\sqrt z)/c_0,\qquad
 B(s) :=\log h(s)\quad(s\ge0).
\end{equation}
Let $\Zp$ denote the multiset of
nontrivial zeta zeros $\rho=\beta+i\gamma$ with $\gamma>0$, and define
\begin{equation*}
 H(t):=\sum_{\rho\in\Zp}e^{-t(\gamma-i(\beta-\frac12))^2}\qquad(t>0).
\end{equation*}
Under RH three positivity properties hold: (a) $H$ is \emph{completely monotone} (i.e.,
$(-1)^nH^{(n)}\ge0$ for all $n\ge 0$); (b) by Schoenberg's theorem
\citep{schoenberg1951polya,grochenig2020schoenberg}, $1/(c_0h)$ and $1/\Xi$ are the one-sided and
bilateral Laplace transforms of P\'olya frequency (PF) functions of every order; and (c) $\Xi$ satisfies
every generalized Laguerre inequality on $\reals$. Conversely, each of these properties implies RH
(Proposition~\ref{thm:RH-boundary}; \citealp[Theorem~2.9]{csordas1990necessary}). Without RH, an
off-line zero $\rho$ and its reflection $1-\bar\rho$ both lie in $\Zp$ and contribute complex-conjugate
terms to $H$. So $H$ stays real, but these terms oscillate in $t$ and can in principle make it
negative. We ask what parts of (a)--(c) follow from the available bounds on those
oscillations. %

\vskip4pt
\noindent\textbf{\emph{Our main results}} are the following. Items marked~$^*$ overlap in part with earlier
or parallel work; the precise relation is described in
Section~\ref{sec:related} and in the discussion of the corresponding result.
\begin{enumerate}[leftmargin=2.8em,labelsep=0.5em]
\item[(1)$^*$] Let $\Theta$ on $(0,\infty)$ and $\Lambda$ on $\reals$ be the
positive densities with one-sided and bilateral Laplace transforms
$1/(c_0h)$ and $1/\Xi$, respectively. They exist without RH. More generally, for every $\alpha>0$, $h^{-\alpha}$ is the Laplace
transform of a positive self-decomposable law, and
$[c_0/\xi(\tfrac12+\|u\|)]^\alpha$ is strictly positive definite on
$\reals^d$ for every $d$ (Theorem~\ref{thm:prob-main}). The existence
question of \citet{roynette2005wald} was answered by
\citet[v24, Theorem~A]{polson2026hilbert}; the proof here uses zeros only
through height $100$ and adds the regularity of the densities. This result answers a positive-definiteness question raised by
\citet{grochenig2020schoenberg}.

\item[(2)] The densities $\Theta$ and $\Lambda$ are strictly log-concave;
the proof uses zeros only through height $100$. Using higher zero verification through
height $3\,000\,175\,332\,800$, we prove
$\Theta,\Lambda\in\PF_{k_*}$ for $k_*=1\,810\,540\,931\,075$
(Theorem~\ref{cor:Lambda-PF6}). For each $k$, the $k\times k$ minors formed by
$\Theta(x-y)$ are also positive when all differences are sufficiently
small and positive (Proposition~\ref{prop:Theta-local-PF}).

\item[(3)$^*$] If RH holds up to height $T>\tfrac12$, then $\Xi$ satisfies the generalized Laguerre inequalities
$L_n[\Xi](x)\ge0$ of every order $n$ for $|x|\le T-\tfrac12$, where
$L_n[\Xi](x)$ is the coefficient of $y^{2n}$ in $|\Xi(x+iy)|^2$
(Theorem~\ref{thm:Laguerre-headline}). The case $n=1$ goes back to
\citet{csordas1991laguerre}. %
\end{enumerate}

The inputs to our results are the critical strip, explicit zero counting, and published zero verification. Our conclusions do not extend the verified
height, and they are merely weak necessary conditions.
The global Laguerre inequality $L_1[\Xi]\ge0$
(\citealp[Problem~2]{csordas1998convexity}) remains open, and in our view, it is a much stronger necessary condition than those studied in this paper.
Below, in \S\S\ref{sec:thm1}--\ref{sec:lag} we present the formal statements of our main results.

\subsection{Construction of the reciprocal laws $\Theta$ and $\Lambda$}
\label{sec:thm1}
Assuming RH, the densities $\Theta$ and $\Lambda$ exist and are $\PF_\infty$ by Schoenberg's theorem \citep[Theorems~1 and~2]{schoenberg1951polya}, specialized to $\Xi$ in \citet[Theorems~3--5]{grochenig2020schoenberg}. But without RH, at least their existence follows easily. Using the notation from~(\ref{eq:normalization}) and~(\ref{eq:FhB}), from Hadamard's product for $h$ over the zeros (Lemma~\ref{lem:product}) we obtain $B'(s)=\int_0^\infty e^{-st}H(t)\dd t$. So $H>0$ makes $B$ a
\emph{Bernstein function}, nonnegative with completely monotone
derivative, and $h^{-\alpha}=e^{-\alpha B}$ the Laplace transform of a
probability law. %

\begin{theorem}[Reciprocal xi and probability laws]\label{thm:prob-main}
For $t>0$, $H(t)>0$ and $H'(t)<0$. The function $B$ is Bernstein,
with the representation
\begin{equation}\label{eq:main-levy}
 B(s)=\int_0^\infty(1-e^{-st})\tfrac{H(t)}{t}\dd t.
\end{equation}
Moreover, for every $\alpha>0$ there is a positive self-decomposable random variable
$S_\alpha$ with $\E e^{-sS_\alpha}=h(s)^{-\alpha}$. It has a $C^\infty$
density positive on $(0,\infty)$. For every $d\ge1$,
\begin{equation}\label{eq:radial-main}
 h(\|u\|^2)^{-\alpha}
 =c_0^\alpha\xi(\tfrac12+\|u\|)^{-\alpha}
\end{equation}
is strictly positive definite on $\reals^d$ and is the characteristic
function of an infinitely divisible law. Its Gaussian-mixture density is positive on $\reals^d$ and extends to an
entire function on $\C^d$.
\end{theorem}

\subsection{Log-concavity and finite PF orders}
As noted, RH implies $\Theta$ and $\Lambda$ are $\PF_\infty$ functions; without
RH, we can still obtain lower PF orders. A function $f$ is $\PF_m$ if
$\det[f(x_i-y_j)]_{i,j=1}^k\ge0$ for every $1\le k\le m$ and every pair
of increasing real tuples $x_1<\cdots<x_k$ and $y_1<\cdots<y_k$.
For a nonnegative function with interval support, log-concavity is equivalent
to $\PF_2$.

\begin{theorem}[Log-concavity and PF$_k$]
\label{cor:Lambda-PF6}
The densities satisfy $(\log\Theta)''(x)<0$ for $x>0$ and
$(\log\Lambda)''(x)<0$ for $x\in\reals$.
Put $k_*=1\,810\,540\,931\,075$. The extension of $\Theta$ by zero to
$(-\infty,0]$ belongs to $\PF_{k_*}$.
For every $1\le k\le k_*$ and increasing real tuples $x_1<\cdots<x_k$
and $y_1<\cdots<y_k$,
\begin{equation}\label{eq:Theta-PF-main}
 \det[\Theta(x_i-y_j)]_{i,j=1}^k\ge0,
\end{equation}
with strict inequality if and only if $x_i>y_i$ for every $i$.
The density $\Lambda$ belongs to $\PF_{k_*}$ with strict minors:
for every $1\le k\le k_*$ and
increasing real tuples $x_1<\cdots<x_k$ and $y_1<\cdots<y_k$,
\[
 \det[\Lambda(x_i-y_j)]_{i,j=1}^k>0.
\]
\end{theorem}
The large finite order uses the published verified height; its scalar
criterion is Theorem~\ref{thm:signed-PF-criterion}.
Appendix~\ref{app:signed-numerics} establishes that $k_*$ is the maximum
of that explicitly specified sufficient criterion. %

\subsection{Generalized Laguerre inequalities for $\Xi$}
\label{sec:lag}
Under RH, $\Xi$ satisfies the generalized Laguerre inequalities on $\reals$; we prove that it satisfies all of them on an interval, which is once again merely a weak necessary check. For a real entire function $f$, define the (generalized) Laguerre transform
\begin{equation}\label{eq:Laguerre-definition}
 L_n[f](x) :=\frac1{(2n)!}\nlsum_{j=0}^{2n}(-1)^{j+n}
 \tbinom{2n}{j}f^{(j)}(x)f^{(2n-j)}(x),\qquad n\in\N_0.
\end{equation}
These are the coefficients of $y^{2n}$ in $|f(x+iy)|^2$; in particular,
$L_0[f]=f^2$ and $L_1[f]=(f')^2-ff''$. Satisfying the generalized Laguerre inequality $L_n[f]\ge 0$ for each $n\ge0$ implies that $f$ is in the  Laguerre--P\'olya class, so for instance the truth of $L_n[\Xi]\ge 0$ would imply RH. \citet{csordas1991laguerre} proved $L_1[\Xi]\ge0$ on an interval using known zero verification \citep[Remark~4.8]{csordas2015fourier}; we extend that result to show that $L_n[\Xi]\ge 0$ on that interval.

\begin{theorem}[Generalized Laguerre on an interval]
\label{thm:Laguerre-headline}\label{cor:Laguerre-verified}
Let $T>1/2$, and suppose all nontrivial zeta zeros $\beta+i\gamma$ with height
$0<\gamma\le T$ satisfy $\beta=1/2$. Then, 
\begin{equation*}
 L_n[\Xi](x)\ge0
 \qquad\bigl(n=0,1,2,\ldots,\quad |x|\le T-\tfrac12\bigr).
\end{equation*}
Moreover, $L_n[\Xi](x)>0$ unless $x$ is a zero of $\Xi$ of multiplicity exceeding $n$.
\end{theorem}
By \citet[Theorem~1]{platt2021riemann} the critical line hypothesis holds for
$T=\Tstar:=3\cdot10^{12}$, so every $L_n[\Xi]\ge0$ on
$|x|\le\Tstar-\tfrac12$. For $n=1$, this use of verified zeros goes back to \citet[Corollary to Lemma~4]{csordas1991laguerre}, who assume that the zeros are simple; see also \citet[p.~398]{csordas1992level}, \citet[remark after Problem~2]{csordas1998convexity} and \citet[Remark~4.8]{csordas2015fourier}. %

\section{Proof of Theorem~\ref{thm:prob-main}: the reciprocal laws}\label{sec:inputs}\label{sec:ideas}\label{sec:ideas-prob}
This section proves Theorem~\ref{thm:prob-main}. The construction needs two analytic facts: the canonical
product must yield a convergent Laplace integral, and $H$ must stay
positive and decreasing. We establish the convergence first, then
compare the possible negative tail with one retained real zero. Sections~\ref{sec:prob} and~\ref{sec:pf} then construct the laws and the densities $\Theta$ and $\Lambda$. Endpoint asymptotics, used for the regularity of the laws, are proved in
Appendix~\ref{app:heat-details}.

\subsection{Analytic and numerical inputs}
We use the classical analytic continuation, functional equation, zero symmetries, and critical strip of $\zeta$ \citep[Chapter~II]{titchmarsh1986theory}. In particular $\delta_\rho:=\beta-\tfrac12$ satisfies $|\delta_\rho|<1/2$, $\xi$ has order one, and the nontrivial zeros are its only zeros. An elementary theta representation, recalled in Appendix~\ref{app:heat-details}, also gives $c_0>0$.

Let $N(Y)$ count upper-half-plane zeros of ordinate at most $Y$, with multiplicity, and write
\[
 M(Y)=\frac{Y}{2\pi}\left(\log\frac{Y}{2\pi}-1\right),\qquad M(0)=0.
\]
We use the following consequence of \citet[Theorem~1.1]{bellotti2025improved}:
\begin{equation}\label{eq:count}
 |N(Y)-M(Y)|\le\tfrac18\log(Y+2)+4\qquad(Y\ge0).
\end{equation}
For all $U\ge100$, \eqref{eq:count} implies
\begin{equation}\label{eq:tailmass}
 N(Y)\le Y\log Y\quad(Y\ge U),\qquad
 \sum_{\gamma>U}\gamma^{-2}
 \le \frac{2(\log U+1)}U<1.
\end{equation}
The sum estimate follows by Stieltjes integration by parts, dropping the nonpositive boundary term.

Put $U_{\rm PT}=3\,000\,175\,332\,800$, the full published verified
height. We also use the round value $\Tstar=3\cdot10^{12}$.
The proved finite zero inputs are
\begin{equation}\label{eq:finite-input}
 \gamma\le U_{\rm PT}\ \Longrightarrow\ \delta_\rho=0,
 \qquad 14<\gamma_1<15,
\end{equation}
where $\gamma_1$ is the first positive zero ordinate. The first assertion follows from \citet[Theorem~1]{platt2021riemann}, which verifies height $3\,000\,175\,332\,800$. The first-zero enclosure is contained in rigorously isolated zero data \citep{platt2017isolating,lmfdb}. Only the portion through height $100$ is needed for Theorem~\ref{thm:prob-main} and the log-concavity in Theorem~\ref{cor:Lambda-PF6}. Table~\ref{tab:inputs} lists the numerical inputs of each result; the certificates are in Appendix~\ref{app:reciprocal}.

For $\rho\in\Zp$ put $\lambda_\rho:=(\gamma-i\delta_\rho)^2$, so that
$H(t)=\sum_{\rho\in\Zp}e^{-t\lambda_\rho}$.

\begin{lemma}[Product and convergence]\label{lem:product}
The function $h$ is entire of order $1/2$, and
\begin{equation}\label{eq:product}
 h(z)=\prod_{\rho\in\Zp}(1+z/\lambda_\rho).
\end{equation}
The product converges locally uniformly. If $a_\rho=\real\lambda_\rho$, then $a_\rho>0$, $\sum a_\rho^{-1}<\infty$, and $h(s)>0$ for $s\ge0$. The series for $H$ and its derivatives converge locally absolutely on $\real t>0$.
\end{lemma}
\begin{proof}
Evenness turns the order-one Xi function into an order-$1/2$ entire function under squaring. The pair $\rho,1-\rho$ yields one squared zero; exactly one member is in $\Zp$. Hadamard factorization at order less than one therefore yields \eqref{eq:product}, with constant one. For large ordinates $a_\rho=\gamma^2-\delta_\rho^2\ge\gamma^2-1/4$, and zero counting gives reciprocal-square summability. At low ordinates \eqref{eq:finite-input} makes $a_\rho=\gamma^2>0$. Conjugate factors in \eqref{eq:product} combine into positive squared moduli on the positive axis; real factors are positive. On a compact subset of $\real t>0$, the Gaussian decay $e^{-c\gamma^2+C\gamma}$ absorbs every fixed power of $\gamma$, proving the last claim.
\end{proof}

\subsection{Positivity of $H$ and $-H'$}\label{sec:heat}
A high zero can oscillate only after its exponential weight has become
small. The estimate below makes this delay quantitative for both $H$
and $-H'$. Keeping a single low real zero gives a common positive term
against which all negative contributions of the tail can be measured. The
same comparison is the large-time step in
\citet[v24, Proposition~6.5]{polson2026hilbert} and
\citet[v9, proof of Theorem~22]{polson2026duality}; the estimate below
makes it valid for every $t>0$.

Put $q=\gamma_1^2\in(196,225)$ and
\[
 a_\rho=\gamma^2-\delta_\rho^2,\qquad
 b_\rho=2\gamma\delta_\rho,
 \qquad \lambda_\rho=a_\rho-i b_\rho.
\]
Then $|b_\rho|\le\gamma$ and $|\lambda_\rho|=\gamma^2+\delta_\rho^2$.

\begin{lemma}[One retained real zero]\label{lem:anchor}
Suppose $100\le U\le U_{\rm PT}$, and define
\[
 H_{>U}(t)=\sum_{\gamma>U}e^{-t\lambda_\rho},\quad
 K_U(t)=e^{-qt}+H_{>U}(t),\quad \varepsilon_U=U^4e^{-U}.
\]
For every $t>0$,
\begin{equation*}
 K_U(t)\ge(1-\varepsilon_U)e^{-qt}>0,\qquad
 -K_U'(t)\ge(1-\varepsilon_U)q e^{-qt}>0.
\end{equation*}
In particular $H>0$ and $-H'>0$.
\end{lemma}
\begin{proof}
The real summands of $H$ and $-H'$ are respectively
\[
 e^{-a_\rho t}\cos(b_\rho t),\qquad
 e^{-a_\rho t}[a_\rho\cos(b_\rho t)+b_\rho\sin(b_\rho t)].
\]
Both are nonnegative when $|b_\rho|t\le\pi/2$. A negative contribution requires $t>\pi/(2\gamma)$. For $\gamma\ge100$,
$a_\rho-q\ge\gamma^2-1/4-225\ge3\gamma^2/4$.
When the corresponding contribution is negative, we have
\[
 e^{-(a_\rho-q)t}\le e^{-3\pi\gamma/8}\le e^{-\gamma}.
\]
Its negative part is consequently at most $e^{-qt}e^{-\gamma}$ in the first sum, and at most $q e^{-qt}\gamma^2e^{-\gamma}$ in the second, using $|\lambda_\rho|\le q\gamma^2$.
By \eqref{eq:tailmass} and monotonicity of $y^2e^{-y}$ and $y^4e^{-y}$ for $y\ge100$,
\[
 \sum_{\gamma>U}e^{-\gamma}<U^2e^{-U},\qquad
 \sum_{\gamma>U}\gamma^2e^{-\gamma}<U^4e^{-U}.
\]
Subtract from the retained terms $e^{-qt}$ and $qe^{-qt}$. Finally $H-K_U$ is a finite sum of positive real exponentials, proving the assertions about $H$.
\end{proof}

\subsection{The Bernstein representation and the laws}\label{sec:prob}
We first derive the L\'evy representation of $B$ from the canonical
product. Termwise it is Frullani's integral
$\log(1+s/\lambda)=\int_0^\infty(1-e^{-st})e^{-\lambda t}\dd t/t$; for a
positive real rate $\lambda$ this is the L\'evy representation of a gamma law.
When every $\lambda_\rho>0$, summation gives Thorin's representation of a
generalized gamma convolution \citep[Chapter~3]{bondesson1992generalized}. Without RH the exponents
$\lambda_\rho$ may be complex, and the positivity of the resulting L\'evy
density $H(t)/t$ is the content of Lemma~\ref{lem:anchor}. For the xi
function the identity is written in \citet[v24, Theorem~5.3]{polson2026hilbert}.
Lemma~\ref{lem:Sregularity} supplies smooth positive densities and
negative moments. Gaussian subordination gives the radial characteristic
functions, and the growth of $B$ gives their entire inverse Fourier
transforms.

The integral characterization of Bernstein functions and the complete Bernstein subclass are described in \citet[Chapters~3 and~6]{schilling2012bernstein}; here a Poisson random measure constructs the laws directly.

\begin{prop}[Bernstein representation]\label{prop:bernstein}
Equation~\eqref{eq:main-levy} holds, and for $s\ge0$, $n\ge0$,
\begin{equation*}
 (-1)^nB^{(n+1)}(s)=\int_0^\infty t^n e^{-st}H(t)\dd t>0.
\end{equation*}
Moreover $\int_0^\infty|H(t)|\dd t<\infty$.
\end{prop}
\begin{proof}
Logarithmic differentiation of \eqref{eq:product} yields
$B'(s)=\sum_\rho(s+\lambda_\rho)^{-1}$. Interchanging with the elementary Laplace integral is justified by
\[
 \sum_\rho\int_0^\infty|e^{-(s+\lambda_\rho)t}|\dd t
 =\sum_\rho(s+a_\rho)^{-1}<\infty.
\]
Thus $B'(s)=\int e^{-st}H(t)\dd t$ and
$\int|H|\le\sum a_\rho^{-1}<\infty$. Integrate in $s$, using $B(0)=0$, to obtain \eqref{eq:main-levy}. The positivity follows from Lemma~\ref{lem:anchor}; differentiation is justified by integrability at zero and exponential decay at infinity. The L\'evy condition follows from
$\int(1\wedge t)H(t)\dd t/t<\infty$.
\end{proof}

For $\alpha>0$ take a Poisson random measure on $(0,\infty)$ of intensity
$\alpha H(t)\dd t/t$, and let $S_\alpha$ be the sum of its points. It is finite almost surely because its expectation is $\alpha\int H(t)\dd t<\infty$. The Laplace functional proves
\begin{equation}\label{eq:S-laplace}
 \E e^{-sS_\alpha}=e^{-\alpha B(s)}=h(s)^{-\alpha}.
\end{equation}
Independent sums add the parameter $\alpha$, so these laws form a convolution semigroup.
For $0<c<1$,
\begin{equation}\label{eq:SD}
 B(s)-B(cs)=\int_0^\infty(1-e^{-st})\frac{H(t)-H(t/c)}t\dd t.
\end{equation}
Since $H$ is decreasing, the L\'evy density $(H(t)-H(t/c))/t$ in \eqref{eq:SD} is nonnegative. Equation~\eqref{eq:SD} therefore constructs an independent $R_c\ge0$ with
$S_\alpha\overset d=cS_\alpha'+R_c$, proving self-decomposability. This is
the classical criterion: a law on $[0,\infty)$ whose L\'evy density has the
form $k(t)/t$ with $k$ nonincreasing is self-decomposable
\citep[Corollary~15.11]{sato1999levy}. For the xi function see also
\citet[v24, Theorem~5.6(ii)]{polson2026hilbert}.

\begin{lemma}[Density, support, and moments]\label{lem:Sregularity}
For every $\alpha>0$, $S_\alpha$ has a $C^\infty(\reals)$ density $\theta_\alpha$, zero on $(-\infty,0]$ and positive on $(0,\infty)$. Its Laplace transform is holomorphic on $\real s>-q$. It has exponential moments of all orders $r<q$, no exponential moment of order $r\ge q$, and finite negative moments of all orders.
\end{lemma}
The proof extracts gamma convolution factors of arbitrarily high order from
the L\'evy measure. This gives smoothness and full positive support;
the remaining moment estimates are in Appendix~\ref{app:law-regularity}.
The exponential-moment threshold is also in
\citet[v24, Corollary~6.8(iii)]{polson2026hilbert} and
\citet[v9, Theorem~22(iii)]{polson2026duality}.

Let $G_d$ be standard Gaussian and independent of $S_\alpha$. The vector
$X_{\alpha,d}=\sqrt{2S_\alpha}G_d$ has characteristic function \eqref{eq:radial-main} and density
\begin{equation*}
 p_{\alpha,d}(x)=\int_0^\infty (4\pi u)^{-d/2}e^{-\|x\|^2/(4u)}\theta_\alpha(u)\dd u>0.
\end{equation*}
Its infinite divisibility follows by adding the parameter $\alpha$, and its self-decomposability follows from \eqref{eq:SD} with $c^2$ instead of $c$. This is the usual Gaussian-subordination mechanism \citep{schoenberg1938metric}; see also \citet[Proposition~I.8.1]{roynette2005wald}. For the xi function with $\alpha=d=1$ it is \citet[v24, Theorem~A]{polson2026hilbert}.
The Gaussian Gram matrix on distinct points is positive definite for each $u>0$: its quadratic form is the integral of the squared modulus of a nonzero exponential polynomial against a positive Gaussian density. Integration against $\theta_\alpha$ proves strictness for \eqref{eq:radial-main}.
By \eqref{eq:B-growth}, its characteristic function is integrable against $e^{A\|t\|}$ for every $A>0$. The inverse Fourier integral therefore defines an entire extension of $p_{\alpha,d}$ to $\C^d$.
This proves Theorem~\ref{thm:prob-main}.

Theorem~\ref{thm:prob-main} also removes the hypothesis from
\citet[eq.~(V.36)]{roynette2005wald}, which they prove under RH: for every
$\alpha>\tfrac12$ there is a random variable $\widetilde H_\alpha\ge0$ with
\begin{equation}\label{eq:RY-shift}
 \E e^{-\lambda^2\widetilde H_\alpha/2}=\frac{\xi(\alpha)}{\xi(\alpha+\lambda)}
 \qquad(\lambda\ge0).
\end{equation}
Indeed, $\widetilde H=2S_1$ satisfies $\E e^{-\beta^2\widetilde H/2}=c_0/\xi(\tfrac12+\beta)$ for
$\beta\ge0$, and their lemma \citep[Lemme~V.3.4]{roynette2005wald},
applied to $\widetilde H$ with shift $\alpha-\tfrac12>0$, gives \eqref{eq:RY-shift}.
For $\alpha\ge1$, $\widetilde H_\alpha/2$ can be taken to be a generalized
gamma convolution \citep[v9, Theorem~18]{polson2026duality}. The question
(V.35) of Roynette and Yor, which asks for the additional factor
$\exp(\lambda\xi'(\alpha)/\xi(\alpha))$ on the right of
\eqref{eq:RY-shift}, is not settled by this argument; see
\citet[v9, Remark~20]{polson2026duality}.

\subsection{The densities $\Theta$ and $\Lambda$}\label{sec:pf}
In the notation of \citet[Theorem~5]{grochenig2020schoenberg}, $\Xi(z)=\Xi_1(-z^2)$ with $\Xi_1=c_0h$.
Set
\begin{equation}\label{eq:ThetaLambda}
 \Theta(u)=\frac{\theta_1(u)}{c_0},\qquad
 \Lambda(x)=\frac{p_{1,1}(x)}{c_0}
 =\frac1{c_0}\int_0^\infty\frac{e^{-x^2/(4u)}}{\sqrt{4\pi u}}\theta_1(u)\dd u.
\end{equation}
Both integrate to $1/c_0$ on their respective domains.

\begin{corr}[Positive versions of the two Laplace representations]\label{cor:Groech}
The density $\Theta$ is smooth, supported in $[0,\infty)$, and strictly positive on $(0,\infty)$. The density $\Lambda$ is even, strictly positive on $\reals$, and extends to an entire function. They satisfy
\begin{align}
 \frac1{c_0h(s)}&=\int_0^\infty e^{-su}\Theta(u)\dd u,
       &&\real s>-\gamma_1^2,\notag\\
 \frac1{\Xi(z)}&=\int_\reals e^{-zx}\Lambda(x)\dd x,
       &&|\real z|<\gamma_1.\label{eq:bilateral}
\end{align}
In particular, for arbitrary $c_j\in\C$ and $t_j\in\reals$,
\begin{equation}\label{eq:PD-exact}
 \sum_{j,k}\frac{c_j\overline{c_k}}{\xi(\tfrac12+t_j-t_k)}
 =\int_\reals\left|\sum_jc_je^{it_jx}\right|^2\Lambda(x)\dd x\ge0.
\end{equation}
The inequality is strict for distinct $t_j$ and a nonzero coefficient vector.
\end{corr}
\begin{proof}
The first transform follows from \eqref{eq:S-laplace} and Lemma~\ref{lem:Sregularity}; equality on the full half-plane follows by holomorphic continuation. For real $|z|<\gamma_1$, conditioning the Gaussian mixture gives
\[
 \int_\reals e^{-zx}\Lambda(x)\dd x
 =c_0^{-1}\E e^{z^2S_1}=\frac1{c_0h(-z^2)}.
\]
The exponential moments of the resulting symmetric variable give holomorphy of the bilateral transform throughout the strip, so the identity theorem proves \eqref{eq:bilateral} there. Take $z=it$ and use the functional equation and evenness to obtain the Fourier transform $1/\xi(\tfrac12+t)$. This proves \eqref{eq:PD-exact}. Strictness follows because $\Lambda>0$ and a nonzero exponential polynomial with distinct frequencies cannot vanish identically.
\end{proof}
Consequently the exact inverse Fourier transform specified in \citet[Theorem~4]{grochenig2020schoenberg} is
\[
 \frac1{2\pi}\int_\reals\frac{e^{-itx}}{\xi(\tfrac12+t)}\dd t=\Lambda(x)>0.
\]
The Gaussian mixture in \eqref{eq:ThetaLambda} gives a positive
representation of this inverse transform and proves the positive-definiteness
property discussed after Theorem~5 of \citet{grochenig2020schoenberg}.
\citet[Section~4.1]{konstantopoulos2024dantzig} ask whether $F$ belongs
to their subclass $\mathcal D_P$ of van Dantzig functions. By
Theorem~\ref{thm:prob-main}, $F$ has the properties that
\citet[Theorem~4.1]{konstantopoulos2024dantzig} derive for members of
$\mathcal D_P$: $\log F(it)=B(t^2)$, and the L\'evy density $\ell$ of $B$
satisfies $(t\ell(t))'=H'(t)<0$, so $B$ lies in their class
$\mathcal B_J$. These properties therefore cannot exclude $F$ from
$\mathcal D_P$. The same observation
is made in \citet[v24, Section~1.2]{polson2026hilbert}, and
\citet[v24, Section~13]{polson2026hilbert} further restricts the possible
representations.

\section{Proof of Theorem~\ref{cor:Lambda-PF6}: Log-concavity and finite PF orders}\label{sec:boundary}\label{sec:finite-translation}
Log-concavity follows from exponential tilting and the low verified zeros.
The higher orders follow by comparing an exact real exponential product
with a small signed correction. Rank comparison under exponential tilting
and a moving polynomial basis give a global bound that reduces all
required determinant signs to scalar inequalities. A separate
analytic theorem transfers the same finite order from $\Theta$ to $\Lambda$.

\subsection{Log-concavity from the low verified zeros}
\label{sec:Theta-PF2}\label{sec:Lambda-PF2}
The one-sided density has a boundary at zero, where a small additive
error need not preserve logarithmic curvature. Exponential tilting
avoids that difficulty. A second low real zero keeps the tilted
residual self-decomposable; its unimodality and its shape before the
mode then control the required primitive on the whole positive axis.
The same tilting argument also gives bilateral log-concavity.
Only the zeros through height $100$ are used in this subsection.

Log-concavity is additional information: the self-decomposable law
$\GammaLaw(\alpha,b)$, with $0<\alpha<1$ and $b>0$, has density
$g(x)=Cx^{\alpha-1}e^{-bx}$ and $(\log g)''(x)=(1-\alpha)/x^2>0$.
Here the work is to verify that the required exponential tilts of the
laws of $S_1$ and $X_{1,1}$ remain self-decomposable, so that classical shape theorems apply.
These reciprocal densities also differ from the Riemann Fourier kernel
$\Phi$, whose $\PF_2$ property is recorded in
\citet[Corollary~2.11]{csordas1998convexity}, and from the xi Taylor
coefficient sequence studied by \citet{katkova2005multiple}.

We use two facts from the theory of self-decomposable laws. Every such
law on $\reals$ is unimodal \citep{yamazato1978unimodality}. A drift-free
positive self-decomposable law with a smooth positive density has
\emph{Yamazato's property}: its density is log-concave before its mode
and nonincreasing after it. This is Yamazato's result, as recalled by
\citet[pp.~50--51]{watanabe1992yamazato}, whose Theorem~1.2 extends it
to unimodal one-sided L\'evy processes without drift.
For a variable $V$ with a finite exponential moment at $r$, its
\emph{exponential tilt} is the probability law
$e^{rx}\Pr(V\in\dd x)/\E e^{rV}$. Tilting a positive infinitely divisible
law with L\'evy density $k(t)/t$ replaces $k(t)$ by $e^{rt}k(t)$.

The next elementary lemma is a localized strict version of the classical
passage from a log-concave density to its distribution function; compare
\citet[Theorem~1 and Corollary~1]{bagnoli2005logconcave}.
They credit Pr\'ekopa for the integration principle and Dierker for
the elementary proof. We give the precise version needed here.

\begin{lemma}[Log-concavity of a primitive]\label{lem:primitive-LC}
Let $g>0$ be $C^1$ on $(0,\infty)$ and locally integrable at zero.
Suppose that, whenever $g'(x)>0$, $\log g$ is concave on $(0,x)$.
Then $G(x)=\int_0^xg(u)\dd u$ satisfies $(\log G)''(x)<0$ for every
$x>0$. In particular this holds for a positive smooth density with
Yamazato's property, or for a positive log-concave function.
\end{lemma}
\begin{proof}
If $g'(x)\le0$, then $G(x)g'(x)-g(x)^2<0$. Otherwise put
$r=g'(x)/g(x)>0$. The tangent inequality for $\log g$ gives
\[
 G(x)\le g(x)\int_0^x e^{r(u-x)}\dd u
       =g(x)\frac{1-e^{-rx}}r<\frac{g(x)}r.
\]
Again $Gg'-g^2<0$, which is exactly the asserted curvature sign.
\end{proof}

\begin{prop}[The one-sided reciprocal density]\label{cor:Theta-PF2}
The one-sided density satisfies
\[
 (\log\Theta)''(x)<0\qquad(x>0).
\]
Its extension by zero to $(-\infty,0]$ belongs to $\PF_2$.
\end{prop}
\begin{proof}
Put $q=\gamma_1^2\in(196,225)$ and let $m_1$ be the multiplicity of
the first ordinate. Choose a real ordinate $\eta\in(21,22)$; its
existence is certified in Appendix~\ref{app:zero-boxes}. We claim that
\begin{equation}\label{eq:Theta-tilted-spectral}
 J(t)=e^{qt}\bigl(H(t)-m_1e^{-qt}\bigr)>0,
 \qquad J'(t)<0\quad(t>0).
\end{equation}
To see this, retain the term at $\eta$ and compare it with
$e^{qt}H_{>100}(t)$. Put $\tau=\eta^2-q>216$. The shifted tail
exponents are $\lambda_\rho-q$. Both their real contributions and
negative derivatives are nonnegative when $|b_\rho|t\le\pi/2$.
A negative contribution therefore requires $t>\pi/(2\gamma)$, while
\[
 (a_\rho-q)-\tau=a_\rho-\eta^2\ge3\gamma^2/4,
 \qquad |\lambda_\rho-q|\le2\gamma^2\le\tau\gamma^2.
\]
The proof of Lemma~\ref{lem:anchor} now applies verbatim, giving
\[
 e^{-\tau t}+e^{qt}H_{>100}(t)\ge(1-\varepsilon_{100})e^{-\tau t},
 \quad
 -\frac{\mathrm d}{\mathrm dt}\bigl(e^{-\tau t}+e^{qt}H_{>100}(t)\bigr)
 \ge(1-\varepsilon_{100})\tau e^{-\tau t}.
\]
The other terms in $J$ are positive real exponentials with positive
rates, which proves \eqref{eq:Theta-tilted-spectral}.

Subtracting all $m_1$ first-zero factors gives
$S_1\overset d=\GammaLaw(m_1,q)+V$, with independence and with
$V$ having L\'evy density $(H(t)-m_1e^{-qt})/t$.
All remaining rates have real part greater than $q$, so
$M:=\E e^{qV}<\infty$. The tilted law $V^{(q)}$ has zero drift and
L\'evy density $J(t)/t$, hence is self-decomposable (a tilt in the
opposite direction appears in \citealp[v9, Theorem~22(ii)]{polson2026duality}). The gamma-extraction argument of
Appendix~\ref{app:law-regularity} gives it a smooth positive density:
removing finitely many factors and tilting preserve the small-time
divergence of $J$ and its positive exponential lower bound.

Let $F_q$ be its distribution function. Yamazato's property and
Lemma~\ref{lem:primitive-LC} give $(\log F_q)''<0$. The density of
$\Exp(q)+V$ is
\begin{equation*}
 f(x)=qM e^{-qx}F_q(x),\qquad x>0,
\end{equation*}
so it has the same strict logarithmic curvature. If $m_1>1$, convolve
with the other $m_1-1$ exponential factors. Each step replaces $f$ by
$q e^{-qx}\int_0^x e^{qu}f(u)\dd u$; since $e^{qu}f(u)$ is
log-concave, Lemma~\ref{lem:primitive-LC} preserves strict curvature.
This proves the claim for $\theta_1=c_0\Theta$, without assuming simplicity.
A log-concave function on its interval support, extended by zero, is
$\PF_2$. The support can force some ordered minors to vanish.
\end{proof}

\begin{prop}[The bilateral reciprocal density]\label{cor:Lambda-PF2}
The density $\Lambda$ satisfies $(\log\Lambda)''(x)<0$ on $\reals$.
Every ordered $\PF_2$ minor is strictly positive.
\end{prop}
\begin{proof}
Gaussian subordination \citep[Theorem~30.1]{sato1999levy} gives the
L\'evy density of $X_{1,1}$ as $A(|x|)/|x|$
(for the Xi function, cf.\ \citealp[v24, Section~2 and Proposition~12.1]{polson2026hilbert}), where
\begin{equation*}
 A(x)=\sum_{\rho\in\Zp}e^{-\gamma x}\cos(\delta_\rho x),\qquad x>0.
\end{equation*}
Indeed, the Gaussian integral applied to each term of $H$ gives
$e^{-(\gamma-i\delta_\rho)|x|}/|x|$. Absolute convergence for $x\ne0$
justifies the interchange: the integrals of the absolute values sum to
$|x|^{-1}\sum_\rho e^{-\sqrt{a_\rho}|x|}<\infty$.
Conjugation gives the displayed real sum.
Put $a=\gamma_1<15$. For every $|s|<a$, both $e^{sx}A(x)$ and
$e^{-sx}A(x)$ are positive and decreasing. Retain the same ordinate
$\eta\in(21,22)$. For either sign, the tail summands and their negative
derivatives are nonnegative for $x\le\pi$, because
$|\delta_\rho|\le1/2$. A negative term therefore has $x>\pi$ and
\[
 e^{-(\gamma-\eta)x}\le e^{-\pi(\gamma-22)},\qquad
 |\gamma\mp s-i\delta_\rho|\le2\gamma
              \le(\eta\mp s)\gamma^2.
\]
The reciprocal-square tail bound makes the negative tail contributions
in either estimate, relative to the retained term, at most
\[
 \varepsilon=100^4e^{-78\pi}<1.
\]
All remaining verified terms have positive coefficients and rates
$\gamma\mp s>0$. The tilt exists because
$\E e^{sX_{1,1}}=\E e^{s^2S_1}<\infty$ for $|s|<a$.
Its L\'evy density is $e^{sx}A(|x|)/|x|$, so the two monotonicity
statements show that it is self-decomposable
\citep[Theorem~15.10]{sato1999levy}, and hence unimodal.

Choose $\eta$ as the retained ordinate in Theorem~\ref{thm:PF-factors}.
This extracts a $\Lap(a)$ factor, with a remaining law having positive
density on the whole line. If $p=c_0\Lambda$, convolution with that
factor gives $|p'(x)|<a p(x)$ for every $x$: the positive contributions
from each side of $x$ make the triangle inequality strict. Were
$(\log p)''(x_0)>0$, the tilt with $s=-p'(x_0)/p(x_0)\in(-a,a)$
would have a strict local minimum at $x_0$, contradicting unimodality.
Therefore $\log p$ is concave. If it were affine on an interval,
analyticity would force $p(x)=Ce^{bx}$ everywhere, which cannot be a
probability density. Thus it is strictly concave, and all
$\PF_2$ minors are positive.
For the stronger pointwise curvature assertion, apply
Theorem~\ref{thm:Gaussian-PF-transfer} with $m=2$ to
$\theta=c_0\Theta$, using Proposition~\ref{cor:Theta-PF2} and the
regularity verified after that theorem. It gives
$\mathscr D_2(c_0\Lambda;x)>0$, hence $(\log\Lambda)''(x)<0$.
This step needs Theorem~\ref{thm:Gaussian-PF-transfer}: strict concavity
allows $(\log\Lambda)''$ to vanish at isolated points, as $-x^4$ does at $0$,
whereas the theorem writes $\mathscr D_2(c_0\Lambda;x)$ as a positive multiple
of a matrix density positive at $xI$.
\end{proof}
These arguments establish order two without the high verified height.
They do not by themselves give order three. For example, if
$0<\alpha<1$ and $b>q$, the $q$-tilt of a $\GammaLaw(\alpha,b)$
residual is self-decomposable. Yet the density $f$ of
$\Exp(q)+\GammaLaw(\alpha,b)$ satisfies $f(x)=Cx^\alpha(1+O(x))$
near zero, with a convergent power series in the last factor, and hence
\[
 -\det[f^{(i+j)}(x)]_{i,j=0}^2
 \sim2C^3\alpha^2(\alpha-1)x^{3\alpha-6}<0.
\]
Thus this density fails $\PF_3$. Higher orders need further information.
The next subsection supplies it through an exact signed correction.
Appendices~\ref{app:Theta-PF4-certificate} and~\ref{app:PF256-certificate}
keep direct certificates for $\Theta$ at order four and $\Lambda$ at order $256$
as independent checks.

\subsection{High finite orders from a signed perturbation}\label{sec:order-ten}
Order two used only the zeros through height $100$; the higher orders use the
full verified height. Replacing each unverified squared rate
$\lambda_\rho=(\gamma-i\delta_\rho)^2$ by the real rate $\gamma^2$ writes the law
of $S_1$ exactly as a real exponential product, which is $\PF_\infty$, convolved
with a signed correction (Lemma~\ref{lem:signed-factor}). The correction is small
in weighted total variation, but that alone does not preserve small determinants.
The key estimate is that, in a basis that moves with the evaluation point, a shift
by $u$ enlarges the Hankel matrix of derivatives of the selected product by at most
the factor $e^{Qu}$, where $Q$ bounds the selected rates uniformly in the order (Theorem~\ref{thm:moving-frame}). Its
proof rests on an exact Grassmannian integral and an order comparison for
exponentially tilted corner arrays. All determinant signs then reduce to scalar
inequalities.

For a verified height $100\le U\le U_{\rm PT}$ and value ${0\le a<U^2-1/4}$, we use the following:
\begin{align}
 \mathcal R_U&=\frac{\log(U/(2\pi))+1}{2\pi U}
       +\frac{\log(U+2)/4+8+1/16}{U^2},\nonumber\\
 c_U(a)&=1-\frac{a+1/4}{U^2},\qquad
 \epsilon_U(a)=
 \exp\!\left\{\frac{\mathcal R_U}{4}
       \left(1+\frac{c_U(a)+2}{c_U(a)^2}\right)\right\}-1.
 \label{eq:signed-scalar-definitions}
\end{align}
Also put $C(T)=M(T)-\log(T+2)/8-4$, the lower bound for $N(T)$
in \eqref{eq:count}.
\begin{theorem}[Scalar criterion for both reciprocal densities]
\label{thm:signed-PF-criterion}
Let $m,n$ be positive integers with $n\ge2m$ and
$\tfrac12(n-2m+1)\ge\sqrt{(m-1)(n-m+1)}$ (for instance $n\ge8m$), and let
$100\le T\le U\le U_{\rm PT}$. Put $Q=T^2$. Suppose $ C(T)\ge n$, $Q<U^2-\tfrac14$, and $\epsilon_U(Q)<1$.
Then, the zero extension of $\Theta$ belongs to $\PF_m$, with
ordered minors strictly positive exactly when $x_i>y_i$ for every $i$.
The density $\Lambda$ has strictly positive ordered minors through
order $m$ on the entire real line.
\end{theorem}
\begin{proof}
The first inequality selects $n$ verified exponential factors with rates at most
$Q$, counted with multiplicity; let $f$ be their product.
Lemma~\ref{lem:signed-factor} writes the law of $S_1$ as the real exponential
product convolved with a signed measure $\mathcal M_U$, and by the last two
inequalities $\mathcal M_U-\delta_0$ has weighted variation at most
$\epsilon_U(Q)<1$ with weight $e^{Qu}$. Theorem~\ref{thm:moving-frame}, whose
proof rests on Lemma~\ref{lem:Grassmann-identity} and
Theorem~\ref{thm:rank-increment}, then gives the minors of $f*\mathcal M_U$ on
all of $(0,\infty)$. Appendix~\ref{app:signed-closure} passes to $\Theta$ through
the support boundary and the remaining infinite convolution. The
Gaussian-mixture Theorem~\ref{thm:Gaussian-PF-transfer}, whose regularity
hypotheses hold for the reciprocal laws, then gives $\Lambda$ at the same order.
\end{proof}

\begin{proof}[Proof of Theorem~\ref{cor:Lambda-PF6}]
Propositions~\ref{cor:Theta-PF2} and~\ref{cor:Lambda-PF2} prove the
strict curvature assertions using only zeros through height $100$.
For the higher order, we take the exact profile
\eqref{eq:signed-certified-profile} at the full published height
$U_{\rm PT}$. The directed scalar inequalities in
\eqref{eq:signed-certified-errors} verify every hypothesis of
Theorem~\ref{thm:signed-PF-criterion}, including the exact factor budget.
The uniform error is below $0.999999764856$, against threshold one.
The theorem gives both global PF assertions, including their different
strictness conditions.
\end{proof}

The order-$k_*$ result also has an operator interpretation. Let $v(g)$ be
the largest number of sign changes of $g$ along increasing finite
sequences, with zero values discarded. For every bounded real Borel
function $g$, the finite-order variation-diminishing theorem
\citep[Chapter~5, Theorem~3.1]{karlin1968} gives
\[
 v(\Lambda*g)\le v(g)\qquad\text{if }v(g)\le k_*-1.
\]
The convolution converges absolutely because $\Lambda$ is integrable.
The bound is the finite-order counterpart of
\citet[Theorem~6]{grochenig2020schoenberg}.

There is also a necessary constraint on large orders. Theorem~\ref{pole:actual}
proves that the first nonreal poles obstruct PF after a precise number of
preceding real pole occurrences, including ties and multiplicities. By
Corollary~\ref{pole:height-corollary}, for $f\in\{\Theta,\Lambda\}$,
\[
 f\in\PF_k,\quad N(Y)\le k+1
 \quad\Longrightarrow\quad \text{RH through height }Y.
\]
For instance, if PF$_{10^{24}}$ held, it would imply RH through height $10^{23}$. %

\section{Proof of Theorem~\ref{thm:Laguerre-headline}: Laguerre inequalities on an interval}\label{sec:laguerre-interval}\label{sec:laguerre-bridge}\label{sec:ideas-laguerre}
This section proves Theorem~\ref{thm:Laguerre-headline} by expressing
each $L_n$ as a coefficient of a product with nonnegative coefficients.

\subsection{Coefficients of the squared modulus}
The definition \eqref{eq:Laguerre-definition} also applies to real
$C^{2n}$ functions. For example,
\begin{equation*}
 L_2[f]=\tfrac14(f'')^2-\tfrac13f'f'''+\tfrac1{12}ff^{(4)}.
\end{equation*}

\begin{lemma}[Generating and multiplication identities]\label{lem:Laguerre-generating}
For a real entire $f$ and real $x,y$,
\begin{equation*}
 |f(x+iy)|^2=f(x+iy)f(x-iy)
       =\sum_{n=0}^\infty L_n[f](x)y^{2n}.
\end{equation*}
For real entire $f,g$,
\begin{equation}\label{eq:Laguerre-product}
 L_n[fg](x)=\sum_{k=0}^n L_k[f](x)L_{n-k}[g](x).
\end{equation}
Local uniform convergence of real entire functions $f_j\to f$ implies
$L_n[f_j]\to L_n[f]$ locally uniformly, for each fixed $n$.
\end{lemma}
\begin{proof}
Multiply the two Taylor series at $x$. Odd powers of $y$ cancel, and the
coefficient of $y^{2n}$ is
$\sum_{j=0}^{2n}i^j(-i)^{2n-j}f^{(j)}(x)f^{(2n-j)}(x)/(j!(2n-j)!)$.
Since $i^j(-i)^{2n-j}=(-1)^{j+n}$, this is
\eqref{eq:Laguerre-definition}. Absolute convergence on compact sets
justifies the multiplication. Multiplying the generating identities for
$f$ and $g$ proves \eqref{eq:Laguerre-product}. Cauchy's integral formula
implies convergence of each fixed derivative under local uniform
convergence; the last statement follows from the finite sum in
\eqref{eq:Laguerre-definition}.
\end{proof}

\subsection{The zero-exclusion principle}
The Xi function satisfies the hypothesis of the next lemma by
Lemma~\ref{lem:product}.

\begin{lemma}[All-order zero-exclusion principle]\label{thm:Laguerre-exclusion}
Let $f\not\equiv0$ be a real entire function with a representation
\begin{equation}\label{eq:Laguerre-canonical-class}
 f(z)=C z^d e^{-a_0z^2+b_0z}
       \prod_{\zeta\in\mathcal Z_f^*}E_1(z/\zeta),
 \qquad E_1(w)=(1-w)e^w,
\end{equation}
where $C\in\reals\setminus\{0\}$, $d\in\N_0$, $a_0\ge0$, $b_0\in\reals$,
and $\mathcal Z_f^*$ is the multiset of nonzero zeros, invariant under
complex conjugation, with $\sum_{\zeta\in\mathcal Z_f^*}|\zeta|^{-2}<\infty$.
The product is grouped into real zeros and conjugate pairs.
Let $I\subset\reals$ be an interval such that
\begin{equation*}
 |x-\real\zeta|\ge|\imag\zeta|
 \quad(x\in I,\ \zeta\in\mathcal Z_f^*\setminus\reals).
\end{equation*}
Then $L_n[f](x)\ge0$ for every $n\in\N_0$ and every $x\in I$.

In particular, if every nonreal zero $\zeta=a+ib$ satisfies
$|a|>T$ and $|b|\le\Delta$, with $T>\Delta\ge0$, then
\begin{equation*}
 L_n[f](x)\ge0
 \qquad(n\in\N_0,\ |x|\le T-\Delta).
\end{equation*}
Likewise, if $|\imag\zeta|\le\Delta$ and no nonreal zero has real
part in $[x_-,x_+]$, where $x_+-x_->2\Delta$, the conclusion holds on
$[x_-+\Delta,x_+-\Delta]$.
\end{lemma}
\begin{proof}
Fix $x\in I$ and regard $Y=y^2$ as a formal variable. For a real zero
$r\ne0$,
\[
 E_1((x+iy)/r)E_1((x-iy)/r)
    =\frac{e^{2x/r}}{r^2}\bigl((x-r)^2+Y\bigr),
\]
which has nonnegative coefficients in $Y$. The factor $z^d$ similarly
contributes $(x^2+Y)^d$.

For a conjugate pair $\zeta=a+ib,\overline\zeta=a-ib$, $b\ne0$, set
$Q_{a,b}(z)=(z-a)^2+b^2$. Their canonical-product factors multiply to
\[
 E_1(z/\zeta)E_1(z/\overline\zeta)
   =\frac{Q_{a,b}(z)}{a^2+b^2}
       \exp\!\left(\frac{2az}{a^2+b^2}\right).
\]
With $u=x-a$, direct multiplication gives
\begin{equation*}
 Q_{a,b}(x+iy)Q_{a,b}(x-iy)
     =(u^2+b^2)^2+2(u^2-b^2)Y+Y^2.
\end{equation*}
Its coefficients are nonnegative by $|u|\ge|b|$. The exponential and
normalizing factors contribute the positive, $Y$-independent factor
$e^{4ax/(a^2+b^2)}/(a^2+b^2)^2$.

Finally, the leading factor in \eqref{eq:Laguerre-canonical-class}
contributes
\[
 C^2e^{-2a_0x^2+2b_0x}e^{2a_0Y},
\]
whose series in $Y$ has nonnegative coefficients because $a_0\ge0$.
Every finite grouped product therefore has nonnegative coefficients in
$Y$, by ordinary series multiplication. The canonical products converge
locally uniformly, and hence so do their products at $x+iy$ and $x-iy$
on compact complex $y$-disks. Cauchy's coefficient formula shows that
every coefficient converges to the corresponding coefficient of
$f(x+iy)f(x-iy)$. Those limiting coefficients are $L_n[f](x)$.
This proves the first assertion without division by $f(x)$ and so also
at real zeros, including multiple zeros.

For $|x|\le T-\Delta$,
$|x-a|\ge |a|-|x|>\Delta\ge|b|$, proving the symmetric-interval
assertion. For the final assertion, $a<x_-$ or $a>x_+$, so every
$x\in[x_-+\Delta,x_+-\Delta]$ satisfies $|x-a|>\Delta\ge|b|$.
\end{proof}

\begin{proof}[Proof of Theorem~\ref{thm:Laguerre-headline}]
A zeta zero $\rho=\beta+i\gamma$ corresponds to the Xi zero
$z_\rho=\gamma-i(\beta-\tfrac12)$. Thus the critical strip gives
$|\imag z_\rho|<1/2$, and the hypothesis gives
$|\real z_\rho|>T$ for every nonreal Xi zero. Apply
Lemma~\ref{thm:Laguerre-exclusion} with this $T$ and $\Delta=1/2$.
For the Xi function one can equivalently use the polynomial truncations
$c_0\prod_{\gamma\le Y}(1-z^2/\lambda_\rho)$ directly: real zeros
occur as pairs and nonreal zeros as conjugation-symmetric quartets.

To obtain strictness, fix $x$ in the stated interval and write
$\Xi(z)=(z-x)^{m_x}g(z)$, where $g(x)\ne0$. All nonreal zero pairs
satisfy the strict inequality $|x-\real\zeta|>|\imag\zeta|$.
Thus each remaining real-zero factor has two positive coefficients
in $Y$, and each nonreal pair has three. There are infinitely many
remaining factors, by zero counting. For any $k\ge0$, choose a finite
product $P$ whose polynomial in $Y$ has degree at least $k$; its
$k$th coefficient is positive. The remaining factor $g/P$ satisfies
Lemma~\ref{thm:Laguerre-exclusion} and is nonzero at $x$, so
$L_k[g](x)\ge L_k[P](x)(g/P)(x)^2>0$.
Finally $|\Xi(x+iy)|^2=y^{2m_x}|g(x+iy)|^2$ gives exactly the
stated equality cases, and positive normalization preserves them for $F$.
\end{proof}

\section{Related work}
\enlargethispage{2pt}\label{sec:related}
\subsection{The reciprocal law}
The reciprocal law was asked for, and constructed under RH, by
\citet[eq.~(0.29) and Th\'eor\`eme~V.3.2]{roynette2005wald}.
Polson's work connects reciprocal entire functions, Hadamard products,
and Thorin representations \citep[Section~2.1]{polson2018vandantzig}.
In \citet[v24, Theorems~A and~5.3]{polson2026hilbert}, whose $T$ is our
$S_1$, the L\'evy density $H(t)/t$, self-decomposability, and
Gaussian-mixture construction coincide with ours.
\citet[v9, Theorem~22]{polson2026duality} also gives the
self-decomposable law. The two proofs use the explicit formula near
zero and, at later times, verification through heights $1000$ and
$200$, respectively. Our proof extends their large-time comparison with
one low zero to every $t>0$ (Section~\ref{sec:heat}), so it needs neither
the explicit formula nor verification beyond height $100$. Positive powers and arbitrary dimension follow directly
from infinite divisibility and Gaussian subordination; the smoothness
and positivity of the densities, the negative moments of $S_\alpha$, and the entire
extension are not stated there.\footnote{Theorem~\ref{thm:prob-main},
the equivalence of RH with (ii)--(vi) in
Proposition~\ref{thm:RH-boundary}, and the
observation on \citet[Theorem~4.1]{konstantopoulos2024dantzig} in
Section~\ref{sec:pf}, were obtained independently of version~24 of \citet{polson2026hilbert}. They are in the
first version of this manuscript, recorded in commit
\texttt{eae8e7e0154b} of the author's repository on 11~September 2026
at 13:45~UTC; version~24 was submitted to arXiv at 16:43~UTC that day.}
Both papers also deduce unimodality of $S_1$ from Yamazato's theorem
(\citealp[v24, Corollary~6.8(iv)]{polson2026hilbert};
\citealp[v9, Theorem~22(ii)]{polson2026duality}), which
Theorem~\ref{cor:Lambda-PF6} strengthens to log-concavity.

\subsection{P\'olya frequency orders}
Gr\"ochenig's Xi specialization of Schoenberg's transform theorem
\citep{schoenberg1951polya,grochenig2020schoenberg} identifies
$\Theta,\Lambda\in\PF_\infty$ with RH. For $\Lambda$ this equivalence is
already in \citet[Th\'eor\`eme~V.3.2]{roynette2005wald}; see also
\citet[Corollary~2.4]{konstantopoulos2024dantzig}. Theorem~\ref{cor:Lambda-PF6}
establishes finite PF orders for these densities. Its
log-concavity argument combines Yamazato's one-sided shape property
\citep{yamazato1978unimodality,watanabe1992yamazato} with estimates
for exponential tilts. For higher orders, we separate a real exponential
product from a finite signed correction, and control the latter in
normalized derivative Grams. The classical orbital integral and the
uniform corner-array law \citep{olshanski2013orbital} identify the
finite-exponential determinants; we prove the comparison of their
separately tilted laws at successive ranks and the resulting moving-basis
bound. Classical Wronskian and composition arguments
\citep{karlin1966tchebycheff,karlin1968} give $\Theta$'s support-dependent
strictness. The Hermitian derivative-ensemble identities of
\citet{forster2021polynomial}, combined with a conditional Gaussian matrix,
give the same-order forward implication from $\Theta$ to $\Lambda$.

\subsection{Laguerre inequalities}
\citet[Lemma~4]{csordas1991laguerre} prove the Laguerre inequality on
intervals away from the nonreal zeros. Let
$f(z)=Ce^{-az^2+bz}\prod_\zeta(1-z/\zeta)e^{z/\zeta}$ be a real entire
function with $C\ne0$, $a\ge0$, and every nonreal zero in the strip
$|\imag z|<1$. If all zeros with real part in an interval $[x_-,x_+]$,
where $x_+-x_->2$, are real and simple, then $L_1[f]>0$ on
$[x_-+1,x_+-1]$; see also \citet[Proposition~2.2]{csordas2015fourier}.
Their proof uses the second logarithmic derivative. Our coefficient argument is the classical product argument
for \citet[Theorem~2.9]{csordas1990necessary}, applied to conjugate zero
pairs; it gives every $L_n$ on the same interval
(Lemma~\ref{thm:Laguerre-exclusion}) and, for $\Xi$, determines equality
from zero multiplicity.

\subsection{Equivalent forms of RH}
The equivalence of RH with complete monotonicity of $H$ is a special
case of \citet[v8, Lemma~3]{zhang2023genus}. Polson's RH equivalences
(\citealp[v24, Proposition~3.2 and Theorem~5.6]{polson2026hilbert};
\citealp[v9, Theorem~26]{polson2026duality})
identify complete monotonicity of $H$ with the generalized gamma
convolution property of $S_1$; Proposition~\ref{thm:RH-boundary} also
treats $H/t^\alpha$ for every finite $\alpha>0$.
\citet[v24, Theorem~7.1]{polson2026hilbert} shows moreover that Hankel
positivity of every order at a single time, with both Hankel blocks of
the Stieltjes moment problem, is equivalent to RH;
\citet[v25, Section~8.3]{polson2026hilbertv25} discusses the unshifted
test.

\subsection{The xi function in probability}
The xi function itself also occurs as a probability transform.
\citet{biane2001probability} survey
the probability laws whose transforms are $\xi$ and the Jacobi theta
function. \citet[Theorems~1.1--1.2]{nakamura2015complete} shows that
$\xi(\sigma-it)/\xi(\sigma)$ is a characteristic function of $t$ for
every real $\sigma$, and formulates an RH equivalence through
\emph{pretended infinite divisibility} for every $1/2<\sigma<1$.
Nakamura's signed L\'evy--Khintchine-type representations omit the usual
integrability condition on the total variation of the signed measure.
The quasi-infinitely divisible class studied by
\citet{lindner2018quasi} retains that condition.
\citet[Theorems~1.1--1.2]{nakamura2023infinitely} construct from the
explicit formula a function that is the characteristic function of an
infinitely divisible law if and only if RH holds. Under RH its L\'evy
measure is $\sum_\gamma m_\gamma\gamma^{-2}\delta_{-\gamma}$, where
$\gamma$ ranges over distinct signed zero ordinates and $m_\gamma$
is the multiplicity. For the reciprocal transform treated here, $1/h$
is infinitely divisible unconditionally, with the positive L\'evy density
$H(t)/t$; RH is instead the statement that this density is completely
monotone (Proposition~\ref{thm:RH-boundary}). The normalized pair
$(F(t),1/F(it))$ is a van Dantzig pair
\citep{konstantopoulos2024dantzig}.

\vskip12pt
\noindent\rule{\linewidth}{0.3pt}

\vskip12pt
\begin{aiusage}
The results in this paper were developed through extensive use of an AI
assistant (OpenAI's GPT~Pro). None of the claims were obtained zero-shot:
they emerged through sustained back-and-forth discussion, including failed
attempts, counterexamples, audits, and repairs.
The reciprocal-density revisions, verification programs, and attribution
review also used OpenAI Codex. Anthropic's Claude (Opus~5.5 and Fable~5.1)
was used extensively to assist with the writing. Assistant checks and
program reruns are not external peer review or formal verification.
\end{aiusage}

\aiacknowledgements

\appendix
\addtocontents{toc}{\protect\setcounter{tocdepth}{1}}
\begingroup
\renewcommand{\contentsname}{Appendix contents}
\tableofcontents
\endgroup

\section{Equivalent forms of the Riemann hypothesis}\label{sec:RH-equivalences}
Extended to every order, the positivity properties of
Sections~\ref{sec:prob} and~\ref{sec:boundary} are each equivalent to RH.
Schoenberg's characterization relates $\PF_\infty$ to real zeros of an
entire reciprocal transform \citep{schoenberg1951polya}. Its Xi
specialization (\citealp[Th\'eor\`eme~V.3.2]{roynette2005wald};
\citealp[Theorems~3--5]{grochenig2020schoenberg}) combines
with the theory of complete Bernstein functions to give the following
equivalences. A Bernstein function $B$ is \emph{complete}, written
$B\in\CBF$, when its L\'evy density is completely monotone. The
equivalence of (i) and (ii), and the continuation to the slit plane used
in the proof, are also in
\citet[v24, Propositions~3.2--3.3 and Theorem~5.6(iii)]{polson2026hilbert}
and \citet[v9, Lemma~8 and Theorem~26]{polson2026duality}; the
continuation argument goes back to \citet{polson2017hilbert}.

\begin{prop}[Equivalent forms of the Riemann hypothesis]\label{thm:RH-boundary}
The following assertions are equivalent:
\begin{enumerate}[label=(\roman*)]
\item RH;
\item $H$ is completely monotone;
\item $H(t)/t$ is completely monotone;
\item for one (equivalently, every) finite $\alpha>0$, $H(t)/t^\alpha$ is completely monotone;
\item $[H(u_i+u_j)]$ is positive semidefinite in every finite size for all $u_i>0$;
\item $B\in\CBF$;
\item $\Theta\in\PF_\infty$;
\item $\Lambda\in\PF_\infty$.
\end{enumerate}
Under these equivalent conditions the Hankel kernels of $H$ and $H/t^\alpha$ are strictly totally positive in every order for every $\alpha>0$.
\end{prop}

The central implication is analytic. Complete monotonicity gives a
positive Laplace representation, which extends $h'/h$ or $\log h$ to
the plane slit along the negative real axis. Such an extension excludes
every zero of $h$ off that axis, and hence every zeta zero off the
critical line.
\begin{proof}
We prove the implications in the order
\[
\begin{gathered}
 (\mathrm{i})\Rightarrow(\mathrm{ii})\Rightarrow
 \{(\mathrm{iii}),(\mathrm{iv})\},\qquad
 (\mathrm{ii})\Leftrightarrow(\mathrm{v}),\qquad
 (\mathrm{iii})\Leftrightarrow(\mathrm{vi})\Rightarrow(\mathrm{i}),\\
 (\mathrm{iv})\Rightarrow(\mathrm{i}),\qquad
 (\mathrm{i})\Leftrightarrow(\mathrm{vii})\Leftrightarrow(\mathrm{viii}).
\end{gathered}
\]

Under RH, $\lambda_\rho=\gamma^2>0$, so $H$ is a sum of positive real exponentials and is completely monotone. Products of completely monotone functions are completely monotone, and $t^{-\alpha}$ is completely monotone for every $\alpha>0$; hence (i) implies (ii)--(iv). Complete monotonicity immediately gives all-order Hankel positive semidefiniteness by the Gram representation. Conversely, the representation theorem for exponentially convex functions gives
$H(t)=\int_\reals e^{-tx}\nu(\dd x)$ with $\nu\ge0$; see \citet{widder1934doubly}. Since $H(t)\to0$ as $t\to\infty$, $\nu$ cannot charge $(-\infty,0]$. This makes $H$ completely monotone.

If $H/t$ is completely monotone, the L\'evy-density characterization of complete Bernstein functions \citep[Chapter~6]{schilling2012bernstein} makes $B$ complete Bernstein. Conversely a complete Bernstein function has a completely monotone L\'evy density, so uniqueness of the L\'evy representation gives complete monotonicity of $H/t$.

If $B$ is complete Bernstein, its Pick extension is holomorphic on
$\C\setminus(-\infty,0]$. The identity theorem extends $h=e^B$ from the
positive axis to that domain, where the exponential cannot vanish.
Every zero of the entire function $h$ must therefore lie on the removed
axis. Since $h(0)=1$, its zeros are negative real. For
$\rho=\beta+i\gamma\in\Zp$ the corresponding zero is
\[
 (\rho-\tfrac12)^2
 =\delta_\rho^2-\gamma^2+2i\gamma\delta_\rho.
\]
As $\gamma>0$, this is negative real exactly when $\delta_\rho=0$.
Thus completeness forces RH by combining slit-plane holomorphy with the
pre-existing entireness of $h$.

More generally, suppose $\ell_\alpha(t)=H(t)t^{-\alpha}$ is completely monotone for one finite $\alpha>0$. Write $\ell_\alpha(t)=\int_0^\infty e^{-rt}\mu(\dd r)$ with $\mu\ge0$. Then for $s>0$, Tonelli gives
\[
 B'(s)=\int_0^\infty e^{-st}H(t)\dd t
 =\Gamma(\alpha+1)\int_0^\infty(s+r)^{-\alpha-1}\mu(\dd r).
\]
At $s=1$ this gives $\int(1+r)^{-\alpha-1}\mu(\dd r)<\infty$.
For every compact $K\subset\C\setminus(-\infty,0]$,
$|z+r|\ge c_K(1+r)$ for $z\in K$, $r\ge0$, with $c_K>0$.
Thus the principal-branch integral defines a holomorphic function $J$
on the slit plane, as in the proof of \citet[v9, Theorem~26]{polson2026duality}. The identity theorem gives $h'=Jh$ there; a zero
of $h$ of multiplicity $m$ would give orders $m-1$ and at least $m$
on the two sides. Hence all zeros of $h$ are negative real, proving
(iv)$\Rightarrow$(i).

Under RH the genus-zero product \eqref{eq:product} becomes
\begin{equation*}
 h(s)=\prod_{\gamma>0}\left(1+\frac{s}{\gamma^2}\right),
 \qquad
 \Xi(z)=c_0\prod_{\gamma>0}\left(1-\frac{z^2}{\gamma^2}\right),
\end{equation*}
where ordinates are counted with multiplicity and
$\sum_{\gamma>0}\gamma^{-2}<\infty$. These are respectively the restricted
and even Laguerre--P\'olya products. Schoenberg's characterization applied
to $1/(c_0h)$ and $1/\Xi$ gives the PF densities; conversely, either PF
representation forces the corresponding product and hence RH
(\citealp[Th\'eor\`eme~V.3.2]{roynette2005wald};
\citealp[Theorems~3--5]{grochenig2020schoenberg}). Uniqueness of the Laplace
transforms identifies those densities with \eqref{eq:ThetaLambda}.
Under RH the representing measures of $H$ and of $H/t^\alpha$ have
infinite support, so every nonzero polynomial square has strictly positive
integral; Proposition~\ref{fin-differential-to-TP} then gives strict total
positivity in every finite order.
\end{proof}

\section{Estimates for \texorpdfstring{$H$}{H} and the reciprocal laws}\label{app:heat-details}
The reciprocal-law construction needs control of $H$ at both ends of
the time axis. Near zero, many high zeros contribute and zero counting
determines the leading scale; at large times, the first zero determines
the decay. These estimates justify the L\'evy integrals and moment bounds
used in the paper. We first check the counting constants and theta
normalization, then derive the asymptotics of $H$ at both ends, including their
fixed-order derivative versions; the last subsection completes the regularity
and moment statements of Lemma~\ref{lem:Sregularity}. Table~\ref{tab:inputs}
lists, for each result of the paper, the inputs beyond elementary calculus and
linear algebra.

\begin{table}[htb]\small
\begin{tabularx}{\textwidth}{@{}p{.33\textwidth}Y@{}}
\toprule
Assertion & Inputs beyond elementary calculus and linear algebra\\
\midrule
Positive $H$ and $-H'$; the Bernstein and self-decomposable laws & Critical strip, zero counting, the verified first zero, and all zeros through height $100$\\
Global log-concavity of $\Theta$ and $\Lambda$ & Reality through height $100$, the real zero enclosed near $21.022$ (Appendix~\ref{app:zero-boxes}), the shifted-exponent estimates, and Yamazato's shape theorems \citep{yamazato1978unimodality,watanabe1992yamazato}\\
PF convolution factors below $U$ & The estimate in Appendix~\ref{sec:kps-factors} and verified critical-line zeros below $U$\\
PF approximation in total variation & A retained zero in $(U/4,U/3]$, the reciprocal-square tail bound, the first exponential or Laplace factor, and the scalar check $225\,b_{\Tstar}<4.5\cdot10^{-9}$ in Corollary~\ref{cor:PF-approximation}\\
$\Theta,\Lambda\in\PF_{k_*}$ with the stated strictness & Verified height $U_{\rm PT}$, zero counting, signed factorization, the orbital integral, tilted rank comparison, moving-basis bound, support-boundary and Gaussian transfers, and directed scalar checks\\
Direct $\Theta\in\PF_4$ certificate & Verified height $\Tstar$, $64$ rigorous real-zero enclosures, uniform endpoint estimates, directed Fourier and residue bounds on closed intervals\\
Direct $\Lambda\in\PF_{256}$ certificate & Verified height $\Tstar$, residual mean and shared-cumulant bounds, $270$ rigorous real-zero enclosures, a closed interval cover, and an analytic tail bound\\
All $L_n[\Xi]\ge0$ on the interval determined by verified zeros & Canonical product, critical-strip width, and finite verification through $\Tstar$\\
\bottomrule
\end{tabularx}
\caption{Numerical inputs by result.}\label{tab:inputs}
\end{table}

\subsection{The counting constant}
For $Y\ge e$, \citet[Theorem~1.1]{bellotti2025improved} give the upper bound
\[
 0.11200\log Y+0.12567\log\log Y+3.77417.
\]
Put $u=\log Y\ge1$. The excess over $(\log Y)/8$ is
$0.12567\log u-0.013u+3.77417$, whose maximum at
$u=0.12567/0.013$ is less than $3.934<4$. At $Y=e$ the same
counting estimate gives $N(e)<4$; thus $N(Y)\le3$ and $-1<M(Y)\le0$
on $[0,e]$, which covers the remaining range. One-sided limits at zero
ordinates give the same inequality for either endpoint convention.
The counting theorem, including its analytic inputs and low-height
checks, is used as published.

\subsection{The theta normalization}
The classical theta formula, in our normalization, is
\begin{equation}\label{eq:theta}
 \Xi(z)=\int_\reals e^{izx}\Phi(x)\dd x,\qquad
 \Phi(x)=\sum_{n\ge1}(4\pi^2n^4e^{9x/2}-6\pi n^2e^{5x/2})
             e^{-\pi n^2e^{2x}}.
\end{equation}
This normalization follows by differentiating the Jacobi theta identity and
integrating by parts; see \citet[Section~2.16]{titchmarsh1986theory} and
\citet[eq.~(2.7)]{konstantopoulos2024dantzig}.
The kernel is even by the theta transformation, and each summand is positive
for $x\ge0$, since $2\pi n^2e^{2x}>3$. It has faster-than-exponential decay
at both ends. Thus $F$ is the characteristic function of the positive
density $\Phi/c_0$.

\subsection{Asymptotics of \texorpdfstring{$H$}{H} at both ends}
\begin{prop}[Asymptotics]\label{prop:heat-asymptotic}
As $t\downarrow0$,
\begin{equation}\label{eq:heat-asymptotic}
 H(t)=\frac{\log(1/t)-\gamma_E-2\log(4\pi)}{8\sqrt{\pi t}}
       +O(\log(1/t)).
\end{equation}
For each fixed nonnegative integer $j$,
\begin{equation}\label{eq:heat-mom-asymptotic}
 (-1)^jH^{(j)}(t)
 \sim \frac{\Gamma(j+\tfrac12)}{8\pi}
           t^{-j-1/2}\log(1/t).
\end{equation}
At infinity $H(t)=m_1e^{-qt}+O(e^{-(q+\varepsilon)t})$ for some $\varepsilon>0$, where $m_1$ is the multiplicity of the first zero. All fixed derivatives have analogous bounds.
\end{prop}
The leading term in \eqref{eq:heat-asymptotic} agrees with
\citet[v24, Corollary~5.5]{polson2026hilbert}, who also obtains the
constant term $7/8$ from the explicit formula.
\begin{proof}
First replace every $\lambda_\rho$ by $\gamma^2$. Stieltjes integration against \eqref{eq:count} gives the main integral
\[
 \frac1{2\pi}\int_0^\infty e^{-ty^2}\log(y/(2\pi))\dd y.
\]
The error is bounded by
$\int_0^\infty 2tye^{-ty^2}[\tfrac18\log(y+2)+4]\dd y
=O(1+\log(1/t))$.
Set $x=\sqrt t\,y$ and use
$\int_0^\infty e^{-x^2}\dd x=\sqrt\pi/2$ and
$\int_0^\infty e^{-x^2}\log x\dd x
=\sqrt\pi(-\gamma_E-2\log2)/4$ to obtain the stated main term.
For $0<t\le1$,
\[
 |e^{t\delta^2}\cos(2t\gamma\delta)-1|
 \le C(t+t^2\gamma^2).
\]
Zero counting then bounds the total replacement error by
$O(\sqrt t\,[1+\log(1/t)])$.
For \eqref{eq:heat-mom-asymptotic}, repeat the argument with
$y^{2j}e^{-ty^2}$. Its counting error is $O(t^{-j}\log(1/t))$; second-order Taylor expansion of
$z^{2j}e^{-tz^2}$ on the short vertical segments from $\gamma$ to $\gamma-i\delta$ gives an error
$O(t^{1/2-j}\log(1/t))$. More explicitly, its real first-order term is zero, and the second derivative is bounded by a constant depending on $j$ times
$(y^{2j-2}+ty^{2j}+t^2y^{2j+2})e^{-ty^2/2}$, with the negative power omitted at $j=0$. Integration of these bounds against $N$ proves the claimed order. The Gaussian main integral has leading constant $\Gamma(j+1/2)/(8\pi)$.
For the large-time statement, low zeros are discrete and real, and every tail real part exceeds $100^2-1/4>q$. There is therefore a positive gap above the minimum real part $q$. Local absolute convergence at $t=1$ controls the remaining exponentially decaying sum, also after any fixed number of differentiations.
\end{proof}

\subsection{Regularity and moments of the reciprocal laws}\label{app:law-regularity}
The Bernstein representation gives a probability measure, but the
later shape arguments need a smooth density that is positive in the
interior of its support. We complete that step by extracting gamma convolution
factors of arbitrarily high smoothness; the asymptotics of $H$ then control the
exponential moments.

\begin{proof}[Proof of Lemma~\ref{lem:Sregularity}]
For each integer $N\ge2$, \eqref{eq:heat-asymptotic} and the lower bound in Lemma~\ref{lem:anchor} allow $b_N>0$ such that
$\alpha H(t)\ge Ne^{-b_Nt}$ for every $t>0$. To see this explicitly, first choose $t_N>0$ so that $\alpha H\ge N$ on $(0,t_N]$, and then choose $b_N\ge q$ large enough that
$Ne^{-(b_N-q)t_N}\le\alpha(1-\varepsilon_{100})$.
Subtracting this L\'evy density gives independent variables with
\begin{equation}\label{eq:gamma-decompose}
 S_\alpha\overset d=\GammaLaw(N,b_N)+R_N.
\end{equation}
The gamma density, extended by zero to the negative half-line, has bounded continuous derivatives through order $N-2$. Convolution with the law of $R_N$ preserves them. As $N$ is arbitrary, $\theta_\alpha$ is smooth.
Every $R_N$ has zero drift: its Laplace exponent divided by $s$ tends to zero, by dominated convergence with the finite first L\'evy moment. For a nonnegative variable $Y$, the elementary Laplace bounds show
$\lim_{s\to\infty}-s^{-1}\log\E e^{-sY}=\inf\supp Y$.
Hence $\inf\supp R_N=0$, and convolution in \eqref{eq:gamma-decompose} makes $\theta_\alpha(x)>0$ for every $x>0$.

The large-time expansion of Proposition~\ref{prop:heat-asymptotic} shows
$\int_1^\infty(e^{rt}-1)H(t)\dd t/t<\infty$ exactly for $r<q$. The Poisson exponential formula, first with truncations and then by monotone convergence, proves the exponential-moment assertions. Holomorphy follows by domination on compact subsets.
Finally Stirling's formula in \eqref{eq:normalization} gives
\begin{equation}\label{eq:B-growth}
 B(s)=\tfrac14\sqrt s\log s+O(\sqrt s)\quad(s\to\infty).
\end{equation}
In particular $h(s)^{-\alpha}\to0$, so there is no atom at zero. For $p>0$, Tonelli's theorem yields
\[
 \E S_\alpha^{-p}=\frac1{\Gamma(p)}\int_0^\infty s^{p-1}h(s)^{-\alpha}\dd s<\infty,
\]
using \eqref{eq:B-growth}.
\end{proof}

\section{Higher P\'olya frequency orders}\label{app:reciprocal}
We first give the complete analytic proof of the scalar criterion in
Theorem~\ref{thm:signed-PF-criterion}: the signed factorization, the
finite-exponential reference, its exact orbital identity, comparison of
tilted corner arrays and rank increments, a global moving-basis bound,
and the two transfers to the actual reciprocal densities. The numerical specialization
is in Appendix~\ref{app:signed-numerics}. We then prove the necessary
first-pole obstruction and the exact scope of finite spatial stitching
in Appendices~\ref{pole:section} and~\ref{app:finite-PF-stitch}.
Subsequent subsections give
fixed-split bounds, positive-residual factorizations, approximation
results, and the direct
certificates as independent results. All accepted numerical inequalities
use directed interval arithmetic. The verified height and
first-zero isolation are published inputs.

\subsection{Wronskian transfers}\label{app:wronskian}\label{sec:tp}
The following is the classical Wronskian criterion for an extended complete
Chebyshev system \citep{karlin1966tchebycheff}. We include the short
proof to fix the signs used in the PF arguments of
this appendix and in the last assertion of
Proposition~\ref{thm:RH-boundary}.
\begin{lemma}\label{fin-evaluation}
Suppose real $C^{k-1}$ functions $f_1,\ldots,f_k$ on an interval
have positive initial Wronskians, where
$W(f_1,\ldots,f_j)=\det[f_b^{(a-1)}]_{a,b=1}^j$:
$W(f_1,\ldots,f_j)>0$ for $1\le j\le k$.
Then $\det[f_j(s_i)]_{i,j=1}^k>0$ whenever $s_1<\cdots<s_k$.
\end{lemma}
\begin{proof}
Induct on $k$. Since $f_1>0$, divide the evaluation row at $s_i$ by $f_1(s_i)$.
Set $g_j=(f_{j+1}/f_1)'$. Direct differentiation gives
\[
 W(g_1,\ldots,g_j)=\frac{W(f_1,\ldots,f_{j+1})}{f_1^{j+1}}>0.
\]
Subtract consecutive normalized rows, starting at the bottom, and expand in the first column. The remaining determinant is the integral, over
$[s_1,s_2]\times\cdots\times[s_{k-1},s_k]$, of the evaluation determinant of $g_1,\ldots,g_{k-1}$.
It is strictly positive in the interior by the induction hypothesis.
\end{proof}

Positive differential Hankel determinants pass to distinct arguments by
the same principle.
\begin{prop}\label{fin-differential-to-TP}
If a smooth real $f$ satisfies
$[(-1)^{i+j}f^{(i+j)}(t)]_{i,j=0}^{k-1}\succ0$ for all $t>0$, then $f(x+y)$ is strictly totally positive through order $k$.
\end{prop}
\begin{proof}
Every initial differential Hankel determinant
$\det[f^{(i+j)}(t)]_{i,j=0}^{r-1}$ is positive: the inserted row and column signs have product one.
Fixing $y$, these are the initial Wronskians of
$x\mapsto f(x+y),f'(x+y),\ldots,f^{(r-1)}(x+y)$.
Lemma~\ref{fin-evaluation} first proves positivity of
$\det[f^{(j-1)}(x_i+y)]_{i,j=1}^r$.
Fixing the $x_i$, these determinants are the transposes of the initial Wronskians of the functions $y\mapsto f(x_i+y)$.
A second application of the lemma proves the assertion.
\end{proof}

\subsection{Differential signs at the support boundary}\label{app:PF-support-boundary}
For a function with the required derivatives write
\[
 \mathscr D_k(f;x)=(-1)^{k(k-1)/2}
                       \det[f^{(i+j)}(x)]_{i,j=0}^{k-1}.
\]
\begin{lemma}[Differential criterion at a support boundary]
\label{lem:PF-support-boundary}
Let $m\ge1$ and let $f\in C^{2m-2}(\reals)$ vanish on $(-\infty,0]$.
Suppose that $\mathscr D_k(f;x)>0$ for every $x>0$ and $1\le k\le m$.
Then, for strictly increasing real tuples $x_1<\cdots<x_k$ and
$y_1<\cdots<y_k$, with $k\le m$,
\[
 \det[f(x_i-y_j)]_{i,j=1}^k\ge0,
 \qquad
 \det[f(x_i-y_j)]_{i,j=1}^k>0
 \quad\Longleftrightarrow\quad x_i>y_i\quad(1\le i\le k).
\]
\end{lemma}
\begin{proof}
We first extend the evaluation criterion used in Lemma~\ref{fin-evaluation}.
Suppose $a_1<\cdots<a_n$, each $g_j\in C^{n-1}(\reals)$ vanishes for
$t\le a_j$, and $W(g_1,\ldots,g_r)(t)>0$ for $t>a_r$, $r\le n$.
Then $\det[g_j(t_i)]\ge0$ at increasing nodes, with strict inequality
exactly when $t_i>a_i$ for all $i$. The case $n=1$ is immediate.
If $t_i\le a_i$, the first $i$ rows have nonzero entries in at
most $i-1$ columns. Otherwise $g_1>0$ on $[t_1,t_n]$. Put
$h_j=(g_{j+1}/g_1)'$ for $t>a_1$. These functions extend by zero to the left
of $a_{j+1}$, belong to $C^{n-2}$, and satisfy
\[
 W(h_1,\ldots,h_r)
   =\frac{W(g_1,\ldots,g_{r+1})}{g_1^{r+1}}>0
   \qquad(t>a_{r+1}).
\]
Dividing rows by $g_1(t_i)$ and taking consecutive row differences gives
\[
 \frac{\det[g_j(t_i)]_{i,j=1}^n}{\prod_i g_1(t_i)}
  =\int_{t_1}^{t_2}\cdots\int_{t_{n-1}}^{t_n}
       \det[h_j(v_i)]_{i,j=1}^{n-1}\dd v_{n-1}\cdots\dd v_1.
\]
Induction makes the integrand nonnegative. It is positive on the nonempty
open rectangle $\max(t_i,a_{i+1})<v_i<t_{i+1}$, proving the assertion.

Fix $n\le m$ and set $g_j(t)=f(t+x_{n+1-j})$, with
$a_j=-x_{n+1-j}$. For $t>a_r$, all arguments in the initial Wronskian
of order $r$ are positive. The functions
$f,-f',\ldots,(-1)^{r-1}f^{(r-1)}$ have positive initial Wronskians
$\mathscr D_j(f;\cdot)$ on $(0,\infty)$.
Lemma~\ref{fin-evaluation} therefore gives their positive evaluation
determinants. Reversing the arguments and removing the column signs
cancel the two factors $(-1)^{r(r-1)/2}$, so
$W(g_1,\ldots,g_r)(t)>0$. Apply the preceding assertion at
$t_i=-y_{n+1-i}$. Its evaluation matrix is the transpose of
$[f(x_i-y_j)]$, with both orders reversed, and has the same determinant.
Its strictness condition is precisely $x_i>y_i$ for every $i$.
\end{proof}

\subsection{A signed correction to an all-real product}\label{app:signed-factor}
Throughout the next subsections, total variation means the full variation
mass of a signed measure.
Fix a verified height $100\le U\le U_{\rm PT}$, and retain every occurrence of
each zero, including multiplicities. Define
\begin{align}
 K_0(u)&=\sum_{\rho\in\Zp}e^{-\gamma^2u},\nonumber\\
 \Delta K(u)&=\sum_{\gamma>U}e^{-\gamma^2u}
       \{e^{\delta_\rho^2u}\cos(2\gamma\delta_\rho u)-1\},
 &\Delta\nu(\dd u)&=\frac{\Delta K(u)}u\dd u.\label{eq:signed-kernel}
\end{align}
Conjugate occurrences are counted separately, so $H=K_0+\Delta K$ exactly.
The real rates $\gamma^2$ above $U$ are auxiliary; replacing
$\lambda_\rho$ by these rates does not assert that $\delta_\rho=0$.

\begin{lemma}[Finite signed factor]\label{lem:signed-factor}
Let $\mathcal R_U,c_U(a)$ and $\epsilon_U(a)$ be as in
\eqref{eq:signed-scalar-definitions}. The signed measure $\Delta\nu$ is
finite, with mass $d\in[-\mathcal R_U/4,0]$. The series
\begin{equation}\label{eq:signed-cp}
 \mathcal M_U=e^{-d}\sum_{j=0}^\infty\frac{\Delta\nu^{*j}}{j!},
 \qquad \Delta\nu^{*0}=\delta_0,
\end{equation}
defines a real signed measure of mass one, supported on $[0,\infty)$, and
\begin{equation}\label{eq:signed-weighted-TV}
 \int e^{au}|\dd(\mathcal M_U-\delta_0)(u)|\le\epsilon_U(a)
       \qquad(0\le a<U^2-1/4).
\end{equation}
If $V_0$ is the convergent sum of independent
$\Exp(\gamma^2)$ variables over all positive zero occurrences, then
\begin{equation}\label{eq:signed-exact-law}
                  \operatorname{Law}(S_1)
                    =\operatorname{Law}(V_0)*\mathcal M_U.
\end{equation}
\end{lemma}
\begin{proof}
Stieltjes integration of \eqref{eq:count} gives
\begin{equation}\label{eq:signed-R1}
 \sum_{\gamma>U}\gamma^{-2}\le \mathcal R_U
 =\frac{\log(U/(2\pi))+1}{2\pi U}
   +\frac{\log(U+2)/4+8+1/16}{U^2}.
\end{equation}
Indeed, writing $E(t)=\log(t+2)/8+4$, the error in integrating the
counting main term is at most
$E(U)U^{-2}+2\int_U^\infty E(t)t^{-3}\dd t$.
The inequality $\log(t+2)\le\log(U+2)+\log(t/U)$ evaluates this
to the second term in \eqref{eq:signed-R1}; integrating $M'$ gives
the first term. The endpoint convention for $N$ does not enlarge the bound.

For one occurrence, put $b=\gamma^2-\delta_\rho^2$. The inequalities
$e^v-1\le ve^v$ and $1-\cos v\le v^2/2$ give
\[
 \frac{|e^{-bu}\cos(2\gamma\delta_\rho u)-e^{-\gamma^2u}|}{u}
 \le\delta_\rho^2e^{-bu}
       +2\gamma^2\delta_\rho^2u e^{-bu}.
\]
Since $b-a\ge c_U(a)\gamma^2$ and $\delta_\rho^2\le1/4$,
integration and \eqref{eq:signed-R1} yield
\begin{equation}\label{eq:signed-levy-TV}
 \int e^{au}|\dd\Delta\nu(u)|
 \le \frac{c_U(a)+2}{4c_U(a)^2}\mathcal R_U=:V_U(a).
\end{equation}
This justifies absolute summation. Frullani's integral then gives
\[
 d=-\sum_{\gamma>U}\log(1+\delta_\rho^2/\gamma^2),
                 \qquad -\mathcal R_U/4\le d\le0.
\]
The series \eqref{eq:signed-cp} converges in every weighted variation
norm in \eqref{eq:signed-levy-TV}. Its Laplace transform is
\[
 \int e^{-su}\dd\mathcal M_U(u)
    =\exp\!\left(\int(e^{-su}-1)\dd\Delta\nu(u)\right).
\]
Separating the atom at zero, and using $e^{-d}\ge1$, bounds its
difference from $\delta_0$ in weighted variation by
$e^{-d}-1+e^{-d}(e^{V_U(a)}-1)\le e^{\mathcal R_U/4+V_U(a)}-1$.
This proves \eqref{eq:signed-weighted-TV}. The measure need not be
positive: if $d<0$, its atom has mass greater than one, so it is not a
probability law; if $d=0$, the correction is $\delta_0$.

The sum defining $V_0$ converges almost surely and in $L^1$, since
$\sum\gamma^{-2}<\infty$. Its Laplace exponent has kernel $K_0$.
Multiplying its Laplace transform by that of $\mathcal M_U$ gives
$\exp[-\int(1-e^{-su})H(u)\dd u/u]=h(s)^{-1}$.
Uniqueness for finite signed measures on $[0,\infty)$ proves
\eqref{eq:signed-exact-law}.
\end{proof}

\subsection{The full finite-exponential reference}\label{app:signed-Jacobi}
Here every norm is the Euclidean operator norm. The full reference has
dimension $n$, although only derivatives through $2m-2$ will be applied
to an actual zero-extended convolution.

\begin{lemma}[Finite-exponential Jacobi bound]\label{lem:signed-Jacobi}
Let $f=*_{\ell=1}^n\Exp(q_\ell)$ with $0<q_\ell\le Q$.
For each $x>0$, its entire positive-argument branch $\widetilde f$
has a lower triangular polynomial normalization $C_x$ with Gram
$J_n=\operatorname{diag}((-1)^i)_{i=0}^{n-1}$. In that basis,
multiplication by the derivative variable is a signed Jacobi matrix $B_x$
satisfying
\begin{equation}\label{eq:signed-Jacobi-bound}
 \|B_x\|_2\le\frac{2n}{x}+nQ,\qquad
 C_x[\widetilde f^{(i+j)}(x-u)]_{i,j<n}C_x^T=e^{-uB_x}J_n.
\end{equation}
Repeated rates are allowed.
\end{lemma}
\begin{proof}
Write $D_k(x)=\mathscr D_k(f;x)$, $D_0=1$.
We first prove a positive Laplace representation, valid for $k\le n$:
\begin{equation}\label{eq:signed-network-representation}
 D_k(x)=d_kx^{k(n-k)}\int e^{-xA}\dd\nu_k(A),\qquad
 d_k>0,\quad\supp\nu_k\subset[kq_{\min},kq_{\max}],
\end{equation}
where $\nu_k$ is a probability measure independent of $x$.

For one exponential kernel its ordered minor is
$q^ke^{-q\sum(u_i-v_i)}\det[\mathbf1_{u_i>v_j}]$.
The step determinant is zero or one, and is one on the interlacing
region $v_1<u_1<v_2<u_2<\cdots<v_k<u_k$.
Iterated continuous Cauchy--Binet therefore represents a minor of $f$
as a positive integral over a bounded interlacing network. Take its
endpoint tuples to be $x\delta\xi_i$ and $x(1+\delta\xi_i)$,
with the same increasing $\xi_i$ at both ends, and scale every
intermediate coordinate by $x$. The domain is independent of $x$.
Every increment is nonnegative and their sum is $k$, so the exponent
is $-xA$ with $kq_{\min}\le A\le kq_{\max}$.

Normalize Lebesgue measure on this scaled domain and push it forward
by $A$. Its volume is nonzero for sufficiently small $\delta$:
replacing every exponential by the step kernel gives
$h(x)=x_+^{n-1}/(n-1)!$, and direct differentiation gives
\[
 \mathscr D_k(h;x)
 =\frac{x^{k(n-k)}}{((n-1)!)^k}
       \prod_{j=0}^{k-1}j!(n-1)_{\underline j}>0.
\]
The ratio of the two minors is the Laplace transform of the
normalized pushforward, times $(\prod q_\ell)^k$.
Choose one weakly convergent subsequence as $\delta\downarrow0$;
the measures have common compact support and are independent of $x$.
Confluence of the endpoint minors yields
\eqref{eq:signed-network-representation}, with the same limiting measure
for every $x>0$ and
$d_k=(\prod q_\ell)^k\prod_{j<k}j!(n-1)_{\underline j}/((n-1)!)^k$.
This proves strict positivity.

Compact support permits differentiating the representation. With
expectation and variance under its tilt by $e^{-xA}$,
\begin{align}
 (\log D_k)'&=\frac{k(n-k)}x-\E_{k,x}A,\nonumber\\
 (\log D_k)''&=-\frac{k(n-k)}{x^2}+\operatorname{Var}_{k,x}(A).
 \label{eq:signed-network-derivatives}
\end{align}
All leading Hankel determinants through size $n$ are nonzero, so exact
orthogonalization of $P\mapsto P(D)f(x)$ gives normalized polynomials
of degrees $0,\ldots,n-1$ with squared norms $(-1)^i$.
Multiplication has diagonal entries
$\alpha_i=(\log(D_{i+1}/D_i))'$ and off-diagonal entries
$B_{i,i+1}=b_i$, $B_{i+1,i}=-b_i$. The last row is well defined:
$\widetilde f$ solves $\prod_\ell(D+q_\ell)\widetilde f=0$,
so multiplication takes place modulo $\prod_\ell(z+q_\ell)$.
This also covers repeated rates.

Equation~\eqref{eq:signed-network-derivatives} and the supports of the
probability measures give
\[
 |\alpha_i|\le\frac{n-1}{x}+nQ.
\]
Here both expectations are nonnegative, so the absolute difference
between the ones of orders $i+1$ and $i$ is at most $(i+1)Q$.
The signed Hankel Sylvester identity gives, for $1\le k<n$,
\[
 b_{k-1}^2=\frac{D_{k+1}D_{k-1}}{D_k^2}
       =-(\log D_k)''\le\frac{k(n-k)}{x^2}\le\frac{n^2}{4x^2}.
\]
The entrywise absolute matrix is symmetric, and its row sums are at
most $(2n-1)/x+nQ$. This proves the norm bound.

Freeze the normalized polynomial basis at $x$. Quotient multiplication
shows that the derivatives of the two expressions in
\eqref{eq:signed-Jacobi-bound} agree at $u=0$ at every order. Both sides
are entire exponential polynomials, so the translation identity follows.
\end{proof}

\subsection{The Grassmannian law of a finite-exponential minor}
\label{app:Grassmann-identity}
The compactly supported measure in
\eqref{eq:signed-network-representation} has an exact orbital
interpretation. Identifying the measure, rather than only its support,
will allow comparison of its separately tilted laws at successive ranks.

\begin{lemma}[Exact Grassmannian identity]\label{lem:Grassmann-identity}
Let $f=*_{\ell=1}^n\Exp(q_\ell)$, where $q_\ell>0$ and repeated
rates are allowed. Put $D_k(x)=\mathscr D_k(f;x)$, $D_0=1$, and
\[
 \mathcal C_0=1,\qquad
 \mathcal C_j=\prod_{\ell=0}^{j-1}\ell!,\qquad
 d_k=\Bigl(\prod_{\ell=1}^nq_\ell\Bigr)^k
                        \frac{\mathcal C_k\mathcal C_{n-k}}{\mathcal C_n}.
\]
Let $V$ have normalized Haar law on $U(n)$, and let $P_k$ be the
diagonal rank-$k$ projection with eigenvalues $0^{n-k},1^k$.
For $x>0$ and $0\le k\le n$,
\begin{equation}\label{eq:Grassmann-identity}
 \frac{D_k(x)}{d_kx^{k(n-k)}}
 =\E\exp\!\left[-x\operatorname{tr}
           \bigl(\operatorname{diag}(q_1,\ldots,q_n)VP_kV^*\bigr)\right].
\end{equation}
The constant $d_k$ agrees with that in
\eqref{eq:signed-network-representation}, and $\nu_k$ is exactly
the law of
\[
 A_k=\operatorname{tr}(\operatorname{diag}(q_1,\ldots,q_n)VP_kV^*).
\]
\end{lemma}
\begin{proof}
First assume the rates are pairwise distinct and write $a_j=-q_j$,
relabeling them so that $a_1<\cdots<a_n$. Partial fractions give
\[
 f^{(r)}(x)
 =\Bigl(\prod_{\ell=1}^nq_\ell\Bigr)
       \sum_{j=1}^n
       \frac{a_j^re^{xa_j}}{\prod_{\ell\ne j}(a_j-a_\ell)}
       \qquad(x>0,\ r\ge0).
\]
Cauchy--Binet applied to the resulting finite moment sum expresses
the Hankel determinant as a sum over $k$-element subsets $I$ of
$\{1,\ldots,n\}$. In each term,
\[
 \prod_{\substack{i,j\in I\\i\ne j}}(a_i-a_j)
       =(-1)^{k(k-1)/2}\Delta(a_I)^2.
\]
Thus the internal Vandermonde factors cancel, including the
differential sign, and
\begin{equation}\label{eq:Grassmann-subset}
 \begin{split}
 D_k(x)&=\Bigl(\prod_{\ell=1}^nq_\ell\Bigr)^k S_k(a;x),\\
 S_k(a;x)&=
  \sum_{\substack{I\subset\{1,\ldots,n\}\\|I|=k}}
   \frac{\exp(x\sum_{i\in I}a_i)}
        {\prod_{i\in I,\,j\notin I}(a_i-a_j)}.
 \end{split}
\end{equation}
Empty products give $S_0=1$ and $S_n=e^{x\sum_j a_j}$.
Individual terms of this subset sum need not be positive.

For distinct real vectors $a,b$, the normalized Harish--Chandra
integral is
\begin{equation}\label{eq:Grassmann-HCIZ}
 \int_{U(n)}
 e^{x\operatorname{tr}(\operatorname{diag}(a)V
                      \operatorname{diag}(b)V^*)}\dd V
 =\mathcal C_n x^{-n(n-1)/2}
       \frac{\det[e^{xa_ib_j}]_{i,j=1}^n}{\Delta(a)\Delta(b)},
\end{equation}
where $\Delta(a)=\prod_{i<j}(a_j-a_i)$; see
\citet[Eq.~(12)]{olshanski2013orbital}.
All unitary integrals here use probability Haar measure.
One may equivalently obtain the real-exponential form from the
Fourier form by analytic continuation, since the compact integral
is entire in the spectral arguments.

Conflate the first $n-k$ entries of $b$ to zero and the last $k$
entries to one. Within each group, Taylor confluence replaces its
columns by derivatives of orders $0,1,\ldots$, divided by the
corresponding factorials. Cross-group Vandermonde factors tend to one.
Define the $n\times n$ matrix
\[
 M^{(k)}_{ij}
 =\begin{cases}
      a_i^{\,j-1},&1\le j\le n-k,\\
      e^{xa_i}a_i^{\,j-(n-k)-1},&n-k<j\le n.
   \end{cases}
\]
The determinant quotient in \eqref{eq:Grassmann-HCIZ} consequently
tends to
\[
 \frac{x^{\binom{n-k}{2}+\binom{k}{2}}}
             {\mathcal C_{n-k}\mathcal C_k}
       \frac{\det M^{(k)}}{\Delta(a)}.
\]
Laplace expansion in the last $k$ columns gives
$\det M^{(k)}/\Delta(a)=S_k(a;x)$: its cross-Vandermonde signs
are precisely the signs in \eqref{eq:Grassmann-subset}. Since
$\binom n2-\binom{n-k}2-\binom k2=k(n-k)$, we obtain
\[
 \E e^{x\operatorname{tr}(\operatorname{diag}(a)VP_kV^*)}
 =\frac{\mathcal C_n}{\mathcal C_k\mathcal C_{n-k}}\,
                  x^{-k(n-k)}S_k(a;x).
\]
Together with \eqref{eq:Grassmann-subset}, this proves
\eqref{eq:Grassmann-identity}. The constants coincide because
\[
 \frac{\mathcal C_k\mathcal C_{n-k}}{\mathcal C_n}
 =\prod_{j=0}^{k-1}\frac{j!}{(n-1-j)!}
 =\frac{\prod_{j=0}^{k-1}j!(n-1)_{\underline j}}{((n-1)!)^k}.
\]

To allow repeated rates, approximate a positive rate vector by
pairwise distinct positive vectors. The entire positive branch of
$f$ solves $\prod_\ell(D+q_\ell)\widetilde f=0$, with initial
derivatives zero through degree $n-2$ and derivative
$\prod_\ell q_\ell$ at degree $n-1$. Its companion-matrix
exponential therefore depends continuously on the rates, as do all
fixed derivatives on the positive axis, even when roots coalesce.
The orbital integral is continuous by dominated convergence on
$U(n)$. Hence \eqref{eq:Grassmann-identity} survives the limit;
no simplicity assumption on the selected zeros is needed.

Finally $A_k\in[kq_{\min},kq_{\max}]$. The identity just proved
and \eqref{eq:signed-network-representation} give equality of the
Laplace transforms of its law and $\nu_k$ for every $x>0$.
These measures have compact support, so they are equal by
Laplace-transform uniqueness. Equivalently, their entire transforms
give equality of all moments, and polynomials are dense among
continuous functions on their common compact support. Thus their
tilted expectations are also identical:
\begin{equation}\label{eq:Grassmann-tilted-mean}
 \E_{k,x}A
   =\frac{\E[A_ke^{-xA_k}]}{\E[e^{-xA_k}]}.
\end{equation}
\end{proof}

\subsection{Order comparison for tilted corner arrays}\label{app:GT-rank}
We need a comparison between separately normalized exponential tilts,
not merely a coupling of the untilted projection matrices.
The comparison below uses the full array of corner eigenvalues.

For a weakly increasing vector
$\lambda=(\lambda_1,\ldots,\lambda_n)$, put
$X_\lambda=V\operatorname{diag}(\lambda)V^*$, with $V$ Haar distributed
on $U(n)$. Denote the increasingly ordered eigenvalues of its leading
$r\times r$ corner by $z^{(r)}_1,\ldots,z^{(r)}_r$, and write
$\mathbf z(X_\lambda)$ for the array of proper corners. They interlace:
\[
 z^{(r+1)}_j\le z^{(r)}_j\le z^{(r+1)}_{j+1},
 \qquad z^{(n)}=\lambda.
\]
Put $S_r(X_\lambda)=\sum_{j=1}^r z^{(r)}_j$ and $S_0=0$.
For fixed real fields $c_1,\ldots,c_{n-1}$, let $\E_{\lambda,c}$
denote the probability law obtained from this orbital law by the
positive weight $\exp(\sum_{r<n}c_rS_r)$, divided by its expectation.
This definition applies also to repeated entries of $\lambda$.

\begin{lemma}[Comparison of tilted corner arrays]\label{lem:GT-tilted-order}
Fix the same real fields $c_1,\ldots,c_{n-1}$ for two weakly increasing
top rows $\lambda,\mu$, and suppose $\lambda_j\le\mu_j$ for every $j$.
For every continuous coordinatewise nondecreasing function $\Phi$
of the proper-corner array,
\[
 \E_{\lambda,c}\Phi(\mathbf z)
       \le \E_{\mu,c}\Phi(\mathbf z).
\]
In particular, $\E_{\lambda,c}S_r\le\E_{\mu,c}S_r$ for $r<n$.
The fields may have either sign.
\end{lemma}
\begin{proof}
For $n=1$ there are no proper-corner coordinates, so there is nothing
to prove. Assume $n\ge2$.

First let both top rows be strictly increasing and rational.
For a strictly increasing boundary, the full proper-corner array of
Haar orbital measure has constant density on the interlacing polytope
$\mathrm{GT}(\lambda)$; the density is
$\prod_{j=0}^{n-1}j!/\prod_{i<j}(\lambda_j-\lambda_i)$
by \citet[Corollary~3.2]{olshanski2013orbital}.
The tilted law therefore has density proportional to
\begin{equation}\label{eq:GT-row-field}
 \exp\left(\sum_{r<n}c_r\sum_{j=1}^r z^{(r)}_j\right)
                \mathbf1_{\mathrm{GT}(\lambda)}(\mathbf z).
\end{equation}
The proportionality constant depends on the boundary, but the
fields in the exponent do not.

We prove the comparison on a finite mesh before taking a continuous
limit. Choose $h>0$ such that both top rows lie in $h\mathbb Z^n$,
and restrict every proper-corner coordinate to $h\mathbb Z$.
The admissible set is finite and nonempty. Give each array the positive
weight in \eqref{eq:GT-row-field}. Conditional on all coordinates
except $z^{(r)}_j$, its possible values form a nonempty lattice interval
$[L,U]\cap h\mathbb Z$. Here $L$ is the maximum of its lower
interlacing neighbors and $U$ is the minimum of its upper neighbors,
including fixed boundary coordinates when appropriate. Both endpoints
are nondecreasing functions of the other coordinates and of the top
row. Its conditional weight at $t$ is $e^{c_rt}$.

For any fixed positive weights $w(t)$, their normalized restrictions
to two lattice intervals $[L,U]$ and $[L',U']$ are stochastically
ordered when $L\le L'$ and $U\le U'$. Indeed, outside the common
overlap the CDF comparison is immediate. For a test point in the
overlap, write $A$ for the weight below $L'$, $C$ for the weight
from $L'$ through the test point, $D$ for the weight above it through
$U$, and $E$ for the weight above $U$ through $U'$. The first CDF
is $(A+C)/(A+C+D)$ and the second is $C/(C+D+E)$.
Cross multiplication leaves
$AD+(A+C)E\ge0$. Thus common-uniform inverse-CDF sampling gives
ordered conditional updates. In particular, $c_r$ may have either sign.

For completeness, the single-site update chain is irreducible.
The coordinatewise minimum admissible array for boundary $\lambda$ is
$z^{(r)}_j=\lambda_j$. From any other array there is a legal decrease
of one coordinate by $h$ towards this minimum. If a coordinate above
its target cannot decrease by $h$, some lower neighbor has the same
current value, because all coordinates lie on the same mesh.
That lower neighbor is itself above its target: the target array
satisfies the ordering constraint. A fixed boundary value cannot block
such a decrease, since it would force the target of the original
coordinate to be at least its current value.

Follow a blocking lower neighbor until reaching a coordinate that can
decrease. This terminates: the directed interlacing edges are
\[
 (r+1,j)\longrightarrow(r,j)\longrightarrow(r+1,j+1),
\]
and the integer $2j-r$ increases by one along each edge.
Following lower neighbors therefore follows an acyclic graph in the
opposite direction. Performing the resulting decrease reduces the
finite integer sum of all mesh distances to the minimum by one.
Iteration reaches the minimum. Every such move, and its reversal,
has positive conditional probability. Hence random-scan conditional
resampling is irreducible. It has positive self-loops, so is aperiodic.
The weighted law is stationary by conditional resampling; the finite
chain consequently converges to that unique stationary law.

Start the two boundary chains at their respective minimum arrays,
which are ordered because $\lambda\le\mu$. Use the same random
coordinate and the same inverse-CDF uniform variable at each update.
The monotonicity of the interval endpoints and of the conditional
laws preserves the order of the two arrays at every step.
Taking the finite-chain stationary limit proves the desired
expectation comparison for the two lattice laws.

Now let $h$ tend to zero through common divisors of the rational
boundaries. For a strictly increasing top row the interlacing
polytope is bounded, has nonempty interior in the full space of
$n(n-1)/2$ proper-corner coordinates, and has boundary of Lebesgue
measure zero. Riemann sums therefore converge for both the positive
normalizing integral in \eqref{eq:GT-row-field} and the same integral
with the continuous factor $\Phi$. The arrays lie in a common compact
set, so this also covers functions $\Phi$ not globally bounded.
Dividing the two convergent sums proves the comparison for strictly
increasing rational top rows under the continuous orbital laws.

It remains to justify arbitrary real and repeated top entries.
The tilted orbital law is weakly continuous in the top row.
Use one Haar matrix $V$ for a sequence of top rows tending to $\lambda$.
Then $V\operatorname{diag}(\lambda)V^*$ and the ordered eigenvalues of
each Hermitian corner vary continuously. On a bounded set of top rows
the corner arrays have common compact support, and their exponential
weights are bounded. Dominated convergence therefore applies to the
weighted test integrals and their normalizing integrals; the limiting
normalizing integral is strictly positive.

To approximate ordered weak boundaries while preserving their order,
take $h\downarrow0$ through positive rational values and put
\[
 \lambda^{(h)}_j=h\lfloor\lambda_j/h\rfloor+jh,\qquad
 \mu^{(h)}_j=h\lfloor\mu_j/h\rfloor+jh.
\]
Both rows are rational and strictly increasing,
$\lambda^{(h)}\le\mu^{(h)}$, and they converge to the respective
original boundaries. The comparison already proved passes to the
limit by the orbital continuity just established.
In particular the conclusion applies to top rows consisting of zeros
and ones. No ambient Lebesgue density on a collapsed polytope has
been used.
\end{proof}

\Needspace{12\baselineskip}
\begin{theorem}[Rank increments under exponential tilting]
\label{thm:rank-increment}
Let $D=\operatorname{diag}(q_1,\ldots,q_n)$ be real, and put
$q_-=\min_jq_j$, $Q=\max_jq_j$. Let $P_k$ be an orthogonal
projection of rank $k$, and define, with normalized Haar measure,
\[
 Z_k(x)=\int_{U(n)}
       e^{-x\operatorname{tr}(DVP_kV^*)}\dd V,\qquad
 \Gamma_k(x)=-\frac{\dd}{\dd x}\log Z_k(x).
\]
For every real $x$ and $0\le k<n$,
\begin{equation}\label{eq:rank-increment}
                 q_-\le\Gamma_{k+1}(x)-\Gamma_k(x)\le Q.
\end{equation}
Repeated rates are allowed. The assertion compares the separately
normalized tilted laws at ranks $k$ and $k+1$.
\end{theorem}
\begin{proof}
Compactness of the unitary group justifies differentiation under the
integral for every real $x$. Thus $\Gamma_k(x)$ is the expected
value of $\operatorname{tr}(DVP_kV^*)$ under the indicated normalized
tilt. In particular $\Gamma_0=0$ and $\Gamma_n=\operatorname{tr}D$.

For any Hermitian matrix $X$, let $S_r(X)$ be the trace of its
leading $r\times r$ corner. Since $X_{rr}=S_r-S_{r-1}$,
\begin{equation}\label{eq:rank-energy}
 \operatorname{tr}(DX)
    =q_nS_n+\sum_{r=1}^{n-1}(q_r-q_{r+1})S_r.
\end{equation}
On an orbit with fixed top eigenvalues, $S_n$ is constant and its
factor cancels in the normalization of the exponential tilt.
The fields of the proper-corner array are therefore
\[
                         c_r=x(q_{r+1}-q_r),
\]
the same at every top boundary.

The increasingly ordered top row
$(0^{n-k},1^k)$ is coordinatewise at most
$(0^{n-k-1},1^{k+1})$. Apply Lemma~\ref{lem:GT-tilted-order}
with these boundaries and the same fields. It gives
\[
 \Delta_r:=\E_{k+1,x}S_r-\E_{k,x}S_r\ge0
                    \quad(r<n),\qquad \Delta_n=1,
\]
where each expectation uses its own normalized tilt.
The repeated top eigenvalues, including ranks zero and $n$, are
covered by the limiting part of that lemma.

By Haar invariance, $Z_k$ and $\Gamma_k$ are unchanged when the
$q_j$ are permuted. First list them increasingly. Then $q_n=Q$
and $q_r-q_{r+1}\le0$, so \eqref{eq:rank-energy} gives
\[
 \Gamma_{k+1}-\Gamma_k
       =Q+\sum_{r<n}(q_r-q_{r+1})\Delta_r\le Q.
\]
For the lower bound, list the same rates decreasingly.
The comparison lemma still applies, because it allows arbitrary real
fields. Now $q_n=q_-$ and $q_r-q_{r+1}\ge0$, giving
$\Gamma_{k+1}-\Gamma_k\ge q_-$.
These two arguments can use different couplings: their scalar
left-hand sides agree by unitary invariance. This proves
\eqref{eq:rank-increment}.
\end{proof}

Theorem~\ref{thm:rank-increment} applies to the \emph{same} tilted
means that occur in the finite-exponential reference.
Indeed Lemma~\ref{lem:Grassmann-identity} identifies the normalized
determinant in \eqref{eq:signed-network-representation} with $Z_k$,
including its constant. Differentiation of this positive compact
integral, together with \eqref{eq:signed-network-derivatives}, therefore
identifies the expectations in Lemma~\ref{lem:signed-Jacobi} with
$\Gamma_k$. Its diagonal coefficients satisfy the sharper bound
\begin{equation}\label{eq:rank-alpha-bound}
 \alpha_i(y)
    =\frac{n-2i-1}{y}-(\Gamma_{i+1}(y)-\Gamma_i(y))
    \ge\frac{n-2i-1}{y}-Q,\qquad y>0.
\end{equation}
This conclusion uses the exact orbital identity and the comparison
\emph{after} tilting. A nested coupling before tilting alone would
not imply \eqref{eq:rank-increment}.

\subsection{Uniform signed smoothing in a moving basis}
\label{app:moving-frame}
We now use the rank-increment bound to control actual positive
shifts uniformly up to the support boundary. Write
\[
 \mathsf H_m[g](x)=[g^{(i+j)}(x)]_{i,j=0}^{m-1},
 \qquad J_m=\operatorname{diag}((-1)^i)_{i=0}^{m-1}.
\]
Here $\mathsf H_m[g]$ denotes a derivative Hankel matrix and is
unrelated to the scalar heat trace $H$ defined earlier.
Every matrix norm in this subsection is the Euclidean operator norm.

\begin{theorem}[Moving-frame signed smoothing]\label{thm:moving-frame}
Let $m,n$ be positive integers such that
\begin{equation}\label{eq:moving-frame-budget}
 n\ge2m,\qquad
 \frac{n-2m+1}{2}\ge\sqrt{(m-1)(n-m+1)}.
\end{equation}
In particular, $n\ge8m$ suffices. Let
$f=*_{\ell=1}^n\Exp(q_\ell)$, where $0<q_\ell\le Q$, with
repeated rates allowed. For each $x>0$ there is a lower triangular
matrix $C_x$ with positive diagonal such that
\begin{equation}\label{eq:moving-frame-shift}
 C_x\mathsf H_m[f](x)C_x^T=J_m,\qquad
 \|C_x\mathsf H_m[f](x-u)C_x^T\|_2\le e^{Qu}
                  \quad(u\ge0).
\end{equation}
The shifted matrix in this formula uses the actual zero extension of
$f$, including at $u=x$.

Consequently, if $\eta$ is a finite real signed measure supported on
$[0,\infty)$ and
\begin{equation}\label{eq:moving-frame-error}
             \varepsilon_\eta
               :=\int_{[0,\infty)}e^{Qu}\dd|\eta|(u)<1,
\end{equation}
then $p=f*(\delta_0+\eta)$ satisfies
$\mathscr D_k(p;x)>0$ for every $x>0$ and $1\le k\le m$.
Its zero extension is $\PF_m$, and its ordered minors through
order $m$ are strictly positive exactly when $x_i>y_i$ for
every diagonal pair.
\end{theorem}
\begin{proof}
For $y>0$, let $\mathcal L_y(P)=P(D)f(y)$.
Lemma~\ref{lem:signed-Jacobi} gives normalized polynomials
$Q_i(y;z)$ of degree $i$, with positive leading coefficient, such that
\[
 \mathcal L_y(Q_iQ_j)=s_i\delta_{ij},\qquad s_i=(-1)^i.
\]
Their coefficients are smooth in $y$: exact orthogonalization
expresses the monic coefficients in the analytic moments with
nonvanishing leading determinants, and each normalization is the
positive square root of a nonzero absolute pivot.
Let $C_y$ be the lower triangular coefficient matrix of the first
$m$ polynomials. Then $C_y\mathsf H_m[f](y)C_y^T=J_m$.

We first differentiate this basis. This differentiation is with
respect to its evaluation parameter $y$, not its polynomial
variable $z$. Put $K_y=C_y'C_y^{-1}$. Since a parameter derivative
does not increase polynomial degree, $K_y$ is lower triangular.
Differentiating orthogonality and using
$\partial_y\mathcal L_y(P)=\mathcal L_y(zP)$ for a frozen
polynomial gives
\[
 0=\mathcal L_y((\partial_yQ_i)Q_j)
    +\mathcal L_y(Q_i(\partial_yQ_j))
    +\mathcal L_y(zQ_iQ_j).
\]
In the signed multiplication matrix of
Lemma~\ref{lem:signed-Jacobi}, the entries are
$B_{ii}=\alpha_i$, $B_{i,i+1}=b_i>0$, and
$B_{i+1,i}=-b_i$. For $j<i$ the second term in the preceding
display vanishes, so $(K_y)_{ij}=-B_{ij}$.
For $j=i$ the two basis-derivative terms are equal.
We have therefore obtained the exact closed system
\begin{equation}\label{eq:moving-frame-generator}
 C_y'=K_yC_y,\qquad
 (K_y)_{ii}=-\alpha_i(y)/2,\qquad
 (K_y)_{i,i-1}=b_{i-1}(y),
\end{equation}
with every other entry zero. The first $m$ rows are closed under this
parameter derivative.

Let $\Gamma_k(y)=\E_{k,y}A$, with $\Gamma_0=0$.
Lemma~\ref{lem:Grassmann-identity} identifies these means with
the tilted orbital means, so Theorem~\ref{thm:rank-increment} gives
$\Gamma_{i+1}(y)-\Gamma_i(y)\le Q$.
Equation~\eqref{eq:signed-network-derivatives} and the Jacobi
identities consequently yield
\begin{align}
 \alpha_i(y)
   &=\frac{n-2i-1}{y}-(\Gamma_{i+1}(y)-\Gamma_i(y))
       \ge\frac{n-2i-1}{y}-Q,\nonumber\\
 b_{i-1}(y)^2
   &=\frac{D_{i+1}(y)D_{i-1}(y)}{D_i(y)^2}
       =-(\log D_i(y))''
       \le\frac{i(n-i)}{y^2}\qquad(1\le i<n).
 \label{eq:moving-frame-coefficients}
\end{align}
Write $S_y$ for the strictly lower part of $K_y$.
A lower bidiagonal shift has norm equal to its largest
subdiagonal magnitude. Since $n\ge2m$,
\[
 \|S_y\|_2\le\frac{\sqrt{(m-1)(n-m+1)}}{y}.
\]
The diagonal entries of $K_y$ are at most
$Q/2-(n-2m+1)/(2y)$. Thus
\begin{equation}\label{eq:moving-frame-lognorm}
 \begin{split}
 \lambda_{\max}\!\left(\frac{K_y+K_y^T}{2}\right)
 &\le\frac Q2-\frac{n-2m+1}{2y}
                   +\left\|\frac{S_y+S_y^T}{2}\right\|_2\\
 &\le\frac Q2-\frac{n-2m+1}{2y}
                   +\frac{\sqrt{(m-1)(n-m+1)}}y
 \le\frac Q2.
 \end{split}
\end{equation}
The last step is precisely \eqref{eq:moving-frame-budget}.
For $n\ge8m$, its first term satisfies
$(n-2m+1)/2\ge3n/8$, while
$\sqrt{(m-1)(n-m+1)}\le\sqrt{mn}\le n/\sqrt8<3n/8$,
which proves the asserted sufficient factor budget.

Fix $x>0$. For $0\le u<x$ define
$R_x(u)=C_xC_{x-u}^{-1}$. The order of its factors matters:
differentiation of \eqref{eq:moving-frame-generator} gives
\[
 R_x(0)=I,\qquad R_x'(u)=R_x(u)K_{x-u}.
\]
For an arbitrary row vector $v_0$, its evolved row
$v(u)=v_0R_x(u)$ therefore satisfies
\[
 \frac{\dd}{\dd u}\|v(u)\|_2^2
   =v(u)(K_{x-u}+K_{x-u}^T)v(u)^T
   \le Q\|v(u)\|_2^2.
\]
Gronwall yields $\|R_x(u)\|_2\le e^{Qu/2}$ for the actual nonautonomous
right evolution. Since
\[
 C_x\mathsf H_m[f](x-u)C_x^T
        =R_x(u)J_mR_x(u)^T,
\]
the required norm bound follows for $u<x$.

The finite exponential convolution has a $C^{n-2}$ zero extension,
and all its derivatives through degree $n-2$ vanish at zero.
They are bounded on the real line, being exponential polynomials
on the positive axis and zero on the negative axis. As
$2m-2\le n-2$, the actual matrix
$\mathsf H_m[f](x-u)$ is zero for $u\ge x$.
This proves \eqref{eq:moving-frame-shift} at and beyond the support
boundary.

The same bounded derivatives justify differentiation of
$p=f*(\delta_0+\eta)$ under integration against the finite
variation measure $|\eta|$. Thus
\[
 C_x\mathsf H_m[p](x)C_x^T=J_m+E_x,\qquad
 \|E_x\|_2
 \le\int_{[0,\infty)}e^{Qu}\dd|\eta|(u)
 =\varepsilon_\eta<1.
\]
Each leading compression $J_k+(E_x)_{[k]}$ has the same strict
error bound. The symmetric path $J_k+t(E_x)_{[k]}$,
$0\le t\le1$, remains nonsingular because the least singular
value of $J_k$ is one. Its inertia, and hence its determinant
sign, therefore equals that of $J_k$.
The lower triangularity of $C_x$ makes this leading compression
the congruence of $\mathsf H_k[p](x)$ by the leading $k$ block
of $C_x$; its determinant multiplier is positive.
It follows that $\mathscr D_k(p;x)>0$ for all $k\le m$.

Finally the zero extension of $p$ belongs to $C^{2m-2}(\reals)$
and vanishes on $(-\infty,0]$. Also $p\in L^1(\reals)$ because
$f\in L^1(\reals)$ and $\eta$ has finite variation.
Lemma~\ref{lem:PF-support-boundary}
therefore proves the global PF assertion and the exact strictness
pattern. If $\delta_0+\eta$ has mass one, the $k=1$ conclusion
also shows that $p$ is a positive probability density.
\end{proof}

\subsection{Closure for \texorpdfstring{$\Theta$}{Theta} and the infinite complement}\label{app:signed-closure}
We now prove the $\Theta$ assertion of Theorem~\ref{thm:signed-PF-criterion}.
The counting inequality selects the first $n$ verified zero occurrences
at most $T$, with $q_\ell\le Q=T^2$. Let $f$ be their exponential
convolution and put $p=f*\mathcal M_U$. By
Lemma~\ref{lem:signed-factor}, $\eta=\mathcal M_U-\delta_0$ satisfies
$\int e^{Qu}\dd|\eta|(u)\le\epsilon_U(Q)<1$.
Theorem~\ref{thm:moving-frame} therefore gives the required differential
signs for this single selected core at every $x>0$.

The actual $p$ is $C^{n-2}$ on the line, vanishes on the nonpositive
half-line, and has all required strict differential signs.
Lemma~\ref{lem:PF-support-boundary} makes it $\PF_m$, with strict
ordered minors exactly when $x_i>y_i$ for every $i$.
In particular $p>0$ on the positive half-line. It integrates to one.

Let $\nu_0$ be the infinite convolution of all remaining real
$\Exp(\gamma^2)$ factors, including the auxiliary ones above $U$.
To justify its PF property, retain two of these factors, whose convolution
$g$ is bounded, uniformly continuous and positive on $(0,\infty)$.
Write $v(x)=\E g(x-W)$ for the remaining sum $W\ge0$.
Finite partial sums $W_j$ converge in $L^1$, so uniform continuity
gives $\E g(x-W_j)\to v(x)$ uniformly in $x$.
Each finite convolution is $\PF_\infty$, and passing to the limit
in each fixed determinant proves this property for $v$.
Moreover $\Pr(W<\varepsilon)>0$ for every $\varepsilon>0$:
choose a tail of mean less than $\varepsilon/4$, use Markov's
inequality on that tail, and use independence and the positive
small-value probabilities for the finite initial sum.
Thus $v(x)>0$ exactly for $x>0$.

Equation~\eqref{eq:signed-exact-law} gives $c_0\Theta=p*v$.
Continuous Cauchy--Binet shows that its ordered minor is the integral
over increasing $z_i$ of
$\det[p(x_i-z_j)]\det[v(z_i-y_j)]$, a nonnegative integrand.
Absolute convergence follows by expanding the determinants into products
of finite one-dimensional convolution integrals. If $x_i>y_i$ for all
$i$, take
\[
 y_i<z_i<\min(x_i,y_{i+1})\ (i<k),\qquad y_k<z_k<x_k.
\]
On this nonempty open box the first determinant is positive, and the
second is lower triangular with positive diagonal. The integral is
strictly positive. If some $x_i\le y_i$, the support gives a zero block
with $i$ rows and $k-i+1$ columns, forcing the determinant to vanish.
This proves the full support-dependent assertion for $\Theta$.

\subsection{Gaussian variance mixtures preserve the finite order}\label{app:Gaussian-PF}
The forward implication used for $\Lambda$ is analytic and preserves the
same order.

\begin{theorem}[Forward Gaussian transfer]\label{thm:Gaussian-PF-transfer}
Let $\theta$ be a Schwartz probability density on $\reals$, zero on
$(-\infty,0]$, whose positive random variable $V$ has every negative
moment. Define
\[
 \lambda(x)=\int_0^\infty\frac{e^{-x^2/(4u)}}{\sqrt{4\pi u}}
                         \theta(u)\dd u.
\]
If $\theta\in\PF_m$, then $\lambda$ has strictly positive ordered
translation minors through order $m$ on $\reals$. Also
$\mathscr D_r(\lambda;x)>0$ for all $x\in\reals$, $1\le r\le m$.
\end{theorem}
\begin{proof}
We use the Fourier convention $\widehat f(t)=\int e^{itx}f(x)\dd x$
in this proof. Fix $r\le m$, write
$\Delta(a)=\prod_{i<j}(a_j-a_i)$ and $Z_r=r!\prod_{j=0}^{r-1}j!$.
The derivative-ensemble density
\begin{equation}\label{eq:Gaussian-PF-ensemble}
 P_\theta(a)=\frac{\Delta(a)}{Z_r}
          \det[(-D)^{j-1}\theta(a_i)]_{i,j=1}^r
\end{equation}
is nonnegative by confluence of the second nodes of a PF minor.
It integrates to one: Andreief's identity and integration by parts give
the triangular moment matrix
\[
 \int a^{i-1}(-D)^{j-1}\theta(a)\dd a
 =\begin{cases}(i-1)!\E V^{i-j}/(i-j)!,&i\ge j,\\0,&i<j.
 \end{cases}
\]
All polynomial boundary terms vanish by the Schwartz hypothesis.
The density in \eqref{eq:Gaussian-PF-ensemble} vanishes if an $a_i<0$,
and its eigenvalues are strictly positive almost surely. Choose an
independent Haar unitary $U_r$ and set $S=U_r\operatorname{diag}(a)U_r^*$.

The Harish--Chandra integral followed by Andreief gives, for Hermitian
$T$ of eigenvalues $t_j$,
\begin{equation}\label{eq:Gaussian-PF-matrix-transform}
               \E e^{i\operatorname{tr}(TS)}
                         =\prod_{j=1}^r\widehat\theta(t_j).
\end{equation}
These derivative-ensemble identities are recorded in
\citet[Theorem~2.9(1), Corollary~3.4(1), Remark~3.5]{forster2021polynomial}.
To fix the constants, the unitary integral at distinct $t_j$ is
\[
 \left(\prod_{j=0}^{r-1}j!\right)
 \frac{\det[e^{ia_it_j}]}{\Delta(a)\Delta(it)}.
\]
After integration, the resulting determinant has entries
$(it_j)^{k-1}\widehat\theta(t_j)$, which cancel its Vandermonde
and the normalizer $Z_r$. Repeated $t_j$ follow by continuity.
For diagonal $T$, \eqref{eq:Gaussian-PF-matrix-transform} proves that
the diagonal entries of $S$ are independent with density $\theta$.
Unitary invariance consequently gives
\begin{equation}\label{eq:Gaussian-PF-matrix-Laplace}
 \E e^{-\operatorname{tr}(AS)}=\prod_j L_\theta(a_j)
                    \qquad(A\succeq0),
\end{equation}
where $a_j$ are the eigenvalues of $A$ and $L_\theta$ is the scalar
Laplace transform.

Let $G$ be independent complex Ginibre, with entry density
$\pi^{-1}e^{-|z|^2}$, and define
$X=\sqrt S\,G+G^*\sqrt S$. Gaussian conditioning gives
\[
 \E_G e^{i\operatorname{tr}(TX)}=e^{-\operatorname{tr}(ST^2)},
 \qquad
 \E e^{i\operatorname{tr}(TX)}
       =\prod_jL_\theta(t_j^2)=\prod_j\widehat\lambda(t_j).
\]
Here $\operatorname{tr}(TX)=2\Re\operatorname{tr}(T\sqrt S G)$
and $\|T\sqrt S\|_F^2=\operatorname{tr}(ST^2)$, which fix the
Gaussian normalization. This law is unitarily invariant. Since $S>0$,
its conditional Gaussian is nondegenerate on the real space of Hermitian
matrices and has an everywhere positive density.

All positive and negative moments of $V$ imply that $\lambda$ is
Schwartz: weighted $L^1$ norms of its Gaussian derivatives are bounded
by constants times $\E V^{(p-j)/2}$, and weighted suprema by
$\E V^{(p-j-1)/2}$. Thus \eqref{eq:Gaussian-PF-ensemble} with
$\theta$ replaced by $\lambda$ defines a finite signed spectral measure.
The same transform calculation identifies its unitary lift with the law
of $X$, by uniqueness of Fourier transforms of finite signed measures.

This identification is pointwise, not just almost everywhere. Negative
moments give $\widehat\lambda(t)=\E e^{-Vt^2}\le C_N(1+t^2)^{-N}$
for every $N$. Weyl integration adds only $\Delta(t)^2$, so the
matrix Fourier transform is integrable and its inverse $p_X$ is
continuous. Some event $\varepsilon I\preceq S\preceq MI$ has
positive probability; on it the conditional Gaussian densities have a
common positive lower bound on every compact set. Therefore $p_X$ is
strictly positive everywhere. Away from repeated eigenvalues, Weyl's
formula and continuity imply
\begin{equation}\label{eq:Gaussian-PF-positive-ratio}
 \frac{\det[(-D)^{j-1}\lambda(a_i)]}{\Delta(a)}
             =b_r p_X(\operatorname{diag}a)>0,
                         \qquad b_r>0.
\end{equation}
The exact positive Weyl constant is immaterial to the signs.
For increasing $a_i$, the derivative evaluation determinant is positive.
Fix increasing $x_i$ and apply Lemma~\ref{fin-evaluation} to
$y\mapsto\lambda(x_i-y)$: each initial Wronskian is the transpose
of such a positive determinant at $x_i-y$. This proves all strict
translation minors.

Finally let all $a_i$ tend to $x$ in
\eqref{eq:Gaussian-PF-positive-ratio}. Taylor confluence gives
\[
 \frac{\mathscr D_r(\lambda;x)}{\prod_{j=0}^{r-1}j!}
                         =b_r p_X(xI)>0.
\]
The nonzero limit follows from the positive matrix density at $xI$.
\end{proof}

For the present laws take $\theta=c_0\Theta$ and $\lambda=c_0\Lambda$;
their mixture relation is \eqref{eq:ThetaLambda}. The hypotheses follow
from Lemma~\ref{lem:Sregularity}. More explicitly, its decompositions
$S_1=\GammaLaw(N,b_N)+R_N$ for arbitrarily large integers $N$ give
weighted derivative bounds for $\theta$: choose $N>j+1$, convolve
the $j$th derivative of the gamma density, and use $R_N\le S_1$
in that independent realization to bound all moments of $R_N$.
The corresponding weighted suprema prove that $\theta$ is Schwartz.
All negative moments were proved in Appendix~\ref{app:law-regularity}.
Theorem~\ref{thm:Gaussian-PF-transfer} now completes the $\Lambda$ assertion
of Theorem~\ref{thm:signed-PF-criterion}. It also gives the strict
pointwise logarithmic curvature for $\Lambda$ from
$\Theta\in\PF_2$, using only the low-height inputs.

\subsection{Exact-budget scalar specialization and its scope}\label{app:signed-numerics}
Use the full published verified height
\[
 U=U_{\rm PT}=3\,000\,175\,332\,800,
\]
from \citet[Theorem~1]{platt2021riemann}. The round value
$3\cdot10^{12}$, used by the direct certificates, lies slightly below
that theorem's endpoint. Take the exact profile
\begin{equation}\label{eq:signed-certified-profile}
 \begin{gathered}
 k_*=1\,810\,540\,931\,075,\qquad
 n_*=12\,363\,146\,804\,211,\\
 T_{\rm sel}=\frac{6\,000\,347\,566\,081}{2},\qquad
 Q_*=T_{\rm sel}^2.
 \end{gathered}
\end{equation}
Directed evaluation gives
\begin{equation}\label{eq:signed-certified-errors}
 \begin{gathered}
 C(T_{\rm sel})-n_*>0.889309597244,\\
 \epsilon_U(Q_*)<0.999999764856<1,\qquad Q_*<U^2-\tfrac14,\\
 (n_*-2k_*+1)^2-4(k_*-1)(n_*-k_*+1)
                  =3\,780\,596\,175\,292>0.
 \end{gathered}
\end{equation}
Also $n_*\ge2k_*$. These are all the scalar hypotheses of
Theorem~\ref{thm:signed-PF-criterion}. The weighted error leaves
a margin of about $2.4\cdot10^{-7}$; the replays below check it at
higher precision.
The uniform moving-frame estimate applies on the entire positive axis.

We give the exact capacity calculation for this sufficient criterion.
For fixed $n\ge2m$, squaring the nonnegative sides of
\eqref{eq:moving-frame-budget} shows that its second condition is
\[
 8m^2-(8n+12)m+n^2+6n+5\ge0.
\]
In the range $n\ge2m$, this is equivalent to
\begin{equation}\label{eq:exact-budget-capacity}
 m\le \frac{2n+3-\sqrt{2n^2-1}}4.
\end{equation}
Equivalently, for positive integer $m$, the least permissible factor count is
\begin{equation}\label{eq:exact-budget-minimum}
 n_{\min}(m)=\max\left\{2m,
 \left\lceil4m-3+2\sqrt{2m^2-3m+1}\right\rceil\right\}.
\end{equation}
This is an exact integer calculation: the ceiling of the square root of
$4(2m^2-3m+1)$ is determined by an integer square root and one comparison.
The required count increases with $m$.

Set $R=\mathcal R_U$, $r=R/4$ and $H=\log2-r>0$, and define
\begin{equation}\label{eq:signed-weight-threshold}
 c_* =\frac{r+\sqrt{r^2+8rH}}{2H},\qquad
 A_*=U^2(1-c_*)-\tfrac14.
\end{equation}
For $0\le a<U^2-1/4$, the function $\epsilon_U(a)$ is strictly increasing.
Writing $c=c_U(a)$ shows that $\epsilon_U(a)<1$ is equivalent to
$Hc^2-rc-2r>0$, and hence exactly to $a<A_*$. For the present height,
$0<A_*<U^2-1/4$. Since
\[
 C'(T)=\frac1{2\pi}\log\frac{T}{2\pi}-\frac1{8(T+2)}>0
          \qquad(T\ge100),
\]
every admissible selection height satisfies
$C(T)<C(\sqrt{A_*})$. Outward bounds give
\begin{equation}\label{eq:signed-capacity-errors}
 12\,363\,146\,804\,212.45185057
 < C(\sqrt{A_*})
 <12\,363\,146\,804\,212.45185058.
\end{equation}
Thus the gate $C(T)\ge n$ permits no integer factor count exceeding
$12\,363\,146\,804\,212$. Equation~\eqref{eq:exact-budget-capacity}
then gives at most $k_*$, and the accepted profile above attains that order.
More directly,
\begin{equation}\label{eq:signed-capacity-exclusion}
 \begin{gathered}
 n_{\min}(k_*+1)=12\,363\,146\,804\,218,\\
 C(\sqrt{A_*^{\rm up}})-n_{\min}(k_*+1)<-5.54814942,
 \end{gathered}
\end{equation}
where $A_*^{\rm up}\ge A_*$ is computed using an outward lower bound
for $R$ and outward upper endpoints thereafter. Decreasing $R$ permits
a larger weight, so this exclusion uses the optimistic endpoint.
Larger orders require still more factors and are excluded as well.
This proves exact maximality only for the displayed moving-frame budget,
the fixed height $U_{\rm PT}$, $Q=T^2$, the count gate $C(T)\ge n$,
and the stated variation majorant. Allowing smaller verified heights
cannot improve this maximum: $\mathcal R_U$ is strictly decreasing for
$U\ge100$, whereas $c_U(a)$ increases for fixed admissible $a$.
Since $1+(c+2)/c^2$ decreases for $c>0$, $\epsilon_U(a)$ decreases
with $U$. Thus every profile feasible at a smaller height remains feasible
at $U_{\rm PT}$. Even exploiting the integrality of $N(T)$
would not change the integer order here: it permits at most one additional
factor, still fewer than $n_{\min}(k_*+1)$.

For comparison, the simpler budget $n\ge8m$ at the round height
$U=3\cdot10^{12}$ gives the maximum $1\,545\,299\,548\,368$, and the exact
budget at that height gives $1\,810\,431\,035\,011$; the full published
endpoint accounts for the remaining increase to
\eqref{eq:signed-certified-profile}. The round profile
$m=1\,500\,000\,000\,000$, $n=8m$, $T=2\,950\,000\,000\,000$,
$U=3\cdot10^{12}$ has the much wider error margin
$\epsilon_U(T^2)<6.886657\cdot10^{-10}$.

All evaluations are scalar and independent of the order; they check the
hypotheses of Theorem~\ref{thm:signed-PF-criterion}. The accepted profile
was checked at $100$ and $160$ decimal digits and replayed by an
independent implementation at $180$ digits; the capacity calculation was
also repeated at $220$ digits. All scalar bounds use outward interval
arithmetic \citep{johansson2017arb}; the verified height is the published
input \citep{platt2021riemann}.

The value $k_*$ is the largest order this criterion certifies; a
different analytic estimate, for example a sharper matrix norm bound,
could improve it with the same zero information.
The asymptotic scales are $\mathcal R_U=O(\log U/U)$ and
$C(T)\asymp T\log T$. For fixed $0<\beta<1$, selection at $T=\beta U$
gives $\epsilon_U(T^2)=O_\beta(\log U/U)$ and permits
$m\asymp U\log U$; the exact factor budget changes its constant.

\subsection{Necessary PF orders and the first nonreal poles}
\label{pole:section}
The sufficient bounds above have a converse constraint: a high PF order
forces a long prefix of real poles. We prove the constraint for each
actual reciprocal density, with arbitrary multiplicities and coincident
decay rates.

By \eqref{eq:bilateral}, their transforms are
\[
 G_\Theta(s):=\frac1{c_0h(s)}
       =\frac1{\xi(\tfrac12+\sqrt s)},\qquad
 G_\Lambda(s):=\frac1{\Xi(s)}.
\]
These identities hold on their respective initial convergence domains,
including the imaginary axis. Corollary~\ref{cor:Groech} gives the actual
continuous densities; their positive probability normalization does not
change the PF property. We continue to use the zero extension for
$\Theta$ and the non-strict PF convention.

Write $\rho=\tfrac12+\delta_\rho+i\gamma_\rho$ for a zeta zero with
$\gamma_\rho>0$, and let $m_\rho$ denote its multiplicity. Thus
$|\delta_\rho|<1/2$. The rates of the poles in the left half-plane are
\begin{equation}\label{pole:rates}
 \lambda_\Theta(\rho)=(\gamma_\rho-i\delta_\rho)^2,
 \qquad
 \lambda_\Lambda(\rho)=\gamma_\rho-i\delta_\rho.
\end{equation}
The corresponding pole is $s=-\lambda_f(\rho)$.
Each pole has multiplicity $m_\rho$. We justify these facts, including
the absence of a factor of two in the first multiplicity, in the contour argument below.
Let $N_0(t-)$ and $N_0(t)$ count critical-line zero occurrences of positive
ordinate strictly below, and at most, $t$, respectively.

\begin{theorem}[First nonreal pole, with multiplicities]\label{pole:actual}
Suppose an off-critical zero exists. Set
\[
 \Gamma=\min_{\delta_\rho\ne0}\gamma_\rho,
 \qquad A_\Theta=\min_{\delta_\rho\ne0}
                   (\gamma_\rho^2-\delta_\rho^2),
 \qquad \tau=\sqrt{A_\Theta}.
\]
These minima are attained, and $\tau<\Gamma$.
For $\Theta$, let $r_\Theta$ be the critical-line multiplicity at ordinate
$\tau$ (zero if absent), and let $\mu_\Theta$ be the largest multiplicity
among off-critical zeros whose squared rate has real part $A_\Theta$.
For $\Lambda$, let $r_\Lambda$ be the critical-line multiplicity at
ordinate $\Gamma$ (zero if absent), and let $\mu_\Lambda$ be the largest
off-critical multiplicity at that ordinate. Put
\begin{align}
 d_\Theta&=N_0(\tau-)+\max(0,r_\Theta-\mu_\Theta+1),\label{pole:dTheta}\\
 d_\Lambda&=N_0(\Gamma-)+\max(0,r_\Lambda-\mu_\Lambda+1).
 \label{pole:dLambda}
\end{align}
Then $\Theta\notin\PF_{d_\Theta+1}$ and
$\Lambda\notin\PF_{d_\Lambda+1}$. In particular,
\[
 \Theta\in\PF_k\ \Longrightarrow\ k\le d_\Theta\le N_0(\Gamma-),
 \qquad
 \Lambda\in\PF_k\ \Longrightarrow\ k\le d_\Lambda\le N_0(\Gamma).
\]
The last count for $\Lambda$ includes a possible critical-line zero tied
at the first off-critical ordinate. When that tied zero is absent,
$d_\Lambda=N_0(\Gamma-)$.
\end{theorem}
This is a necessary condition. If RH holds, there is no first nonreal pole and the
conditional upper bounds do not apply.

\begin{corr}[PF order forces zero verification]\label{pole:height-corollary}
Let $N(Y)$ count all nontrivial zeta zero occurrences with
$0<\gamma\le Y$, including multiplicity. For each $f\in\{\Theta,\Lambda\}$,
\begin{equation}\label{pole:height-implication}
 f\in\PF_k\quad\text{and}\quad N(Y)\le k+1
 \quad\Longrightarrow\quad
 \beta=\tfrac12\ \text{for every zero with }0<\gamma\le Y.
\end{equation}
It suffices to have a rigorous real upper bound $N_{\rm up}(Y)<k+2$.
\end{corr}
\begin{proof}[Proof from Theorem~\ref{pole:actual}]
If there were an off-critical zero through $Y$, its first ordinate
$\Gamma\le Y$ would contain a reflected pair
$\tfrac12\pm\delta+i\Gamma$, of equal positive multiplicity.
There are therefore at least two noncritical occurrences at that height.
For either density, the theorem implies
\[
 k\le d_f\le N(\Gamma)-2\le N(Y)-2\le k-1,
\]
a contradiction. The strict upper-bound formulation follows because
$N(Y)$ is an integer. This also covers $Y$ equal to a zero ordinate.
\end{proof}

\subsubsection{Finite translates under weak total positivity}
We give the finite-order ingredient explicitly, including rank-deficient
collocation matrices. A strict sign change means an alternation along an
increasing finite sequence after zero values have been deleted.

\begin{lemma}\label{pole:translates}
If $f\in\PF_k$ and $n\le k$, every real linear combination of $n$ distinct
translates of $f$ has at most $n-1$ strict sign changes.
\end{lemma}
\begin{proof}
For increasing real $a_j$, the initial Wronskians of $e^{a_jx}$ are
positive exponentials times positive Vandermonde products. Thus
Lemma~\ref{fin-evaluation} makes every ordered evaluation determinant
of these exponentials strictly positive.

Here is the finite matrix consequence needed below. Let $A$ be a totally
nonnegative $(n+1)\times n$ matrix. Suppose a vector $y$ in its column
space has $n+1$ nonzero alternating entries. The rank $r$ of $A$ is positive and
at most $n$. Choose $r$ independent columns $A_C$, with the original
column order. All their maximal minors are nonnegative and at least one
is positive. For $\varepsilon>0$ set
\[
 (B_\varepsilon)_{ij}=e^{-(i-j)^2/\varepsilon},\qquad1\le i,j\le n+1.
\]
This matrix tends to the identity. It is strictly totally positive:
positive row and column factors reduce its minors to those of
$e^{2ij/\varepsilon}$, and the exponential evaluation fact applies.
Cauchy--Binet now makes every $r$-row minor of $B_\varepsilon A_C$
strictly positive, since the same positive maximal minor of $A_C$
contributes positively to each sum.

For sufficiently small $\varepsilon$, $B_\varepsilon y$ still has all
$n+1$ alternating signs and belongs to the column space of
$B_\varepsilon A_C$. Restrict to the first $r+1$ rows. The alternating
maximal cofactors of this $(r+1)\times r$ matrix form a left null vector;
all their magnitudes are positive. Its dot product with the first $r+1$
entries of $B_\varepsilon y$ has summands of one nonzero sign, so cannot
vanish. This contradiction proves the matrix consequence. Rank zero is
immediate.

Finally, $n$ strict sign changes of
$g(x)=\sum_{j=1}^n c_jf(x-y_j)$ would give $n+1$ increasing evaluation
points with alternating nonzero values. Order the distinct $y_j$
increasingly. The resulting $(n+1)\times n$ matrix
$[f(x_i-y_j)]$ is totally nonnegative, since every possible minor has
size at most $n\le k$. The preceding consequence excludes its column-space
vector of alternating values. This proves the lemma for the non-strict
PF convention, without assuming full rank or strictly positive minors.
\end{proof}
This is a special case of the classical finite matrix variation-diminishing
theorem \citep{karlin1968}. The proof above supplies exactly
the version needed for finite translates.

\subsubsection{The pole obstruction with equal decay rates}
\begin{theorem}\label{pole:general}
Suppose a real function has, as $x\to+\infty$, an expansion
\begin{equation}\label{pole:finiteexpansion}
 f(x)=\sum_{\lambda\in\mathcal S}e^{-\lambda x}P_\lambda(x)
                      +O(e^{-Bx}),
\end{equation}
where $\mathcal S$ is a finite set of distinct rates, closed under
conjugation, each $P_\lambda$ has exact degree $m_\lambda-1$ with
$m_\lambda\ge1$, and conjugate rates have conjugate polynomials.
Assume that $\mathcal S$ contains a nonreal rate, let $A$ be the least
real part of such a rate, and suppose every rate in $\mathcal S$ has
real part at most $A$, with $B>A$. Thus all nonreal rates in the sum
have real part $A$. Let
\[
 M_< =\sum_{\substack{q\in\mathcal S\cap\reals\\q<A}}m_q,
 \qquad r=\begin{cases}m_A,&A\in\mathcal S,\\0,&A\notin\mathcal S,\end{cases}
 \qquad \mu=\max_{\substack{\lambda\in\mathcal S\setminus\reals\\\real\lambda=A}}
                  m_\lambda.
\]
Then, with $d=M_<+\max(0,r-\mu+1)$, one has $f\notin\PF_{d+1}$.
\end{theorem}
\begin{proof}
Choose $h>0$ small enough that
$0<h|\imag\lambda|<2\pi$ for every nonreal rate in the group at $A$.
Set $E_hf(x)=f(x+h)$, $\ell=\max(0,r-\mu+1)$, and
\begin{equation}\label{pole:annihilator}
 P(z)=\prod_{\substack{q\in\mathcal S\cap\reals\\q<A}}
               (z-e^{-qh})^{m_q}\,(z-e^{-Ah})^\ell.
\end{equation}
It has degree $d$, so $P(E_h)f$ is a combination of $d+1$ translates.
For a polynomial $p$,
\[
 (E_h-e^{-qh})[e^{-qx}p(x)]
       =e^{-q(x+h)}\bigl(p(x+h)-p(x)\bigr).
\]
The factor lowers polynomial degree by one and kills constants.
It follows that all real-rate terms below $A$ are killed. The tied real
term, when present, has remaining degree at most $\mu-2$, or is killed:
if $r\ge\mu$, its degree drops from $r-1$ by $r-\mu+1$;
if $r<\mu$, its original degree already is at most $\mu-2$.
This includes $r=0$ and $\mu=1$.

For a nonreal rate $\lambda$ at $A$, the operator multiplies the leading
coefficient by $P(e^{-\lambda h})\ne0$. The roots associated with $q<A$
have different moduli, and the possible root at $A$ has a different
phase by the choice of $h$. Hence each nonreal term of maximal degree
$\mu-1$ survives, with nonzero leading coefficient. Distinct such rates
have distinct nonzero frequencies. Fixed finite translates preserve the
exponential gap in the remainder, and therefore
\begin{equation}\label{pole:oscillatory}
 P(E_h)f(x)=e^{-Ax}x^{\mu-1}\bigl(T(x)+o(1)\bigr),
\end{equation}
where $T(x)=\sum_\nu a_\nu e^{i\nu x}$ is a nonzero real trigonometric
polynomial, with distinct nonzero frequencies and
$a_{-\nu}=\overline{a_\nu}$.

Direct integration gives its mean zero and positive mean square:
\[
 \lim_{L\to\infty}\frac1L\int_0^L T(x)\dd x=0,
 \qquad
 \lim_{L\to\infty}\frac1L\int_0^L T(x)^2\dd x
                     =\sum_\nu|a_\nu|^2=:s^2>0.
\]
Let $K\ge\max(1,\sup|T|)$ and choose $0<\epsilon<s^2/K$.
If $T\ge-\epsilon$ eventually, then
$(T+\epsilon)(T-K)\le0$ eventually. Taking long-time means gives
$s^2\le K\epsilon$, a contradiction. Apply the same argument to $-T$.
Thus values above $\epsilon$ and below $-\epsilon$ both occur arbitrarily
far out. They can be selected successively in an alternating increasing
sequence. The $o(1)$ term in \eqref{pole:oscillatory} is eventually smaller
than $\epsilon/2$, so $P(E_h)f$ has arbitrarily many strict sign changes.
Lemma~\ref{pole:translates} contradicts this if $f\in\PF_{d+1}$.
No rational independence of frequencies is needed.
\end{proof}

\subsubsection{Finite contour shifts for the actual transforms}\label{pole:contours}
The rates in \eqref{pole:rates} are discrete, have positive real parts,
and have only finitely many occurrences with real part below any fixed
bound. For $\Theta$, use
$\real\lambda_\Theta=\gamma^2-\delta^2$ and $|\delta|<1/2$;
for $\Lambda$, use $\real\lambda_\Lambda=\gamma$.
Positive real parts follow also from the known first-zero isolation.
All repeated occurrences at a rate are collected as its multiplicity.

Evenness of $\xi(\tfrac12+z)$ makes
$\xi(\tfrac12+\sqrt s)$ entire and single-valued. Near a nontrivial
zero, the map $z\mapsto z^2$ has nonzero derivative. Consequently its
zero multiplicity is unchanged by this map. The two zeros corresponding
to opposite values of $z$ give the same squared zero; choosing positive
zeta ordinate selects one of them, without doubling its multiplicity.
There is no branch-cut contribution or pole at zero. For $\Lambda$, the
map defining $\Xi$ is linear; evenness supplies the negative-half-plane
pole $s=-\lambda_\Lambda$ with the same multiplicity. Reflection of zeta
zeros in the critical line gives the conjugate rate, with equal multiplicity.

Fix $B>0$ so that no pole lies on $\real s=-B$. We require uniform
large-$|t|$ estimates for $s=\sigma+it$, $-B\le\sigma\le0$.
For $\Theta$, choose the principal square root solely for this estimate.
Uniformly on the fixed strip,
\[
 \sqrt{\sigma+it}
 =\sqrt{|t|}\,e^{i\operatorname{sgn}(t)\pi/4}+O_B(|t|^{-1/2}).
\]
Thus $w=\tfrac12+\sqrt s$ has real part tending to infinity, and $w/2$
stays in a closed sector within $|\arg z|<\pi$. The Dirichlet series gives
\[
 |\zeta(w)-1|\le2^{-\real w}
                  +\frac{2^{1-\real w}}{\real w-1}=o(1).
\]
Apply sector Stirling \citep[Section~5.11]{dlmf} to
$\xi(w)=\tfrac12w(w-1)\pi^{-w/2}\Gamma(w/2)\zeta(w)$. It gives
\begin{equation}\label{pole:Theta-decay}
 \log|G_\Theta(\sigma+it)|
 =-\frac{\sqrt{|t|}}{4\sqrt2}\log|t|+O_B(\sqrt{|t|}).
\end{equation}
Indeed $\real(w/2)=\sqrt{|t|}/(2\sqrt2)+O(1)$ and
$\log|w/2|=\tfrac12\log|t|-\log2+o(1)$; all other real parts
in Stirling's formula contribute $O_B(\sqrt{|t|})$.

For $\Lambda$, the functional equation gives
\[
 \Xi(\sigma+it)=\xi(\tfrac12+|t|\mp i\sigma),
\]
using the minus sign when $t>0$ and the plus sign when $t<0$.
The same Dirichlet-series estimate and Stirling formula, now with bounded
imaginary part, give uniformly on the strip
\begin{equation}\label{pole:Lambda-decay}
 \log|G_\Lambda(\sigma+it)|
 =-\tfrac12|t|\log|t|
     +\tfrac12(1+\log(2\pi))|t|+O_B(\log|t|).
\end{equation}
In particular, for either transform and every fixed integer $j\ge0$,
\[
 \int_\reals|(-B+it)^jG_f(-B+it)|\dd t<\infty.
\]
The same polynomial multiples tend to zero on horizontal edges uniformly
in $\sigma\in[-B,0]$. All constants here may depend on the fixed cutoff
and derivative order.

\begin{lemma}[Actual finite residue expansion]\label{pole:residues}
For either $f=\Theta$ or $f=\Lambda$, every such $B$ gives
\begin{equation}\label{pole:actualexpansion}
 f^{(j)}(x)=\sum_{\real\lambda<B}
       \left(\frac{\dd}{\dd x}\right)^j
          [e^{-\lambda x}P_\lambda(x)]
                 +R_{B,j}(x),\qquad x\ge0,
\end{equation}
where $P_\lambda$ has exact degree $m_\lambda-1$ and
$|R_{B,j}(x)|\le C_{B,j}e^{-Bx}$ for a finite constant $C_{B,j}$.
The sum is finite and conjugate rates have conjugate polynomials.
\end{lemma}
\begin{proof}
Fourier inversion on the imaginary axis is absolutely convergent by
\eqref{pole:Theta-decay}--\eqref{pole:Lambda-decay}. It recovers the actual
continuous density, since its transform there is already known. The same
bounds justify every fixed number of derivatives. Shift this upward
Bromwich contour to the upward line $\real s=-B$. The bounding rectangle
has positive orientation when its right edge is traversed upward;
the two horizontal integrals vanish. If
\[
 G_f(s)=\sum_{a=1}^{m_\lambda}
        \frac{c_{\lambda,-a}}{(s+\lambda)^a}
           +\text{holomorphic},\qquad c_{\lambda,-m_\lambda}\ne0,
\]
the residue of $e^{sx}G_f(s)$ is
\[
 e^{-\lambda x}P_\lambda(x),\qquad
 P_\lambda(x)=\sum_{a=1}^{m_\lambda}
                   c_{\lambda,-a}\frac{x^{a-1}}{(a-1)!}.
\]
This gives its exact degree and nonzero leading coefficient. Multiplying
the integrand by $s^j$ differentiates each finite residue term. The remaining
vertical integral has absolute value at most
\[
 \frac{e^{-Bx}}{2\pi}
       \int_\reals|(-B+it)^jG_f(-B+it)|\dd t.
\]
It supplies the asserted uniform remainder. Only finitely many residues
are crossed.
\end{proof}

\begin{proof}[Completion of Theorem~\ref{pole:actual}]
For either density choose $B$ strictly above the least nonreal decay rate
and strictly below the next larger rate real part. Local finiteness makes
this possible. Lemma~\ref{pole:residues} then supplies precisely the
hypotheses of Theorem~\ref{pole:general}. For $\Theta$, the real rates are
$\gamma^2$ with critical-line multiplicities, giving \eqref{pole:dTheta}.
For $\Lambda$, they are the critical-line ordinates themselves, giving
\eqref{pole:dLambda}.

An off-critical zero at the first ordinate $\Gamma$ has $\delta\ne0$,
so $A_\Theta\le\Gamma^2-\delta^2<\Gamma^2$. This proves $\tau<\Gamma$,
even if the squared rate attaining $A_\Theta$ belongs to a different,
slightly higher ordinate. Also $\max(0,r-\mu+1)\le r$ since $\mu\ge1$.
Consequently $d_\Theta\le N_0(\tau)\le N_0(\Gamma-)$ and
$d_\Lambda\le N_0(\Gamma)$. The theorem follows.
\end{proof}

\subsubsection{What a quadratic order would imply}
The counting inequality \eqref{eq:count}, derived from
\citet[Theorem~1.1]{bellotti2025improved}, gives
\begin{equation}\label{pole:countup}
 N(Y)\le N_{\rm up}(Y):=
 \frac{Y}{2\pi}\left(\log\frac{Y}{2\pi}-1\right)
                       +\frac18\log(Y+2)+4.
\end{equation}
Directed scalar evaluation gives
\begin{equation}\label{pole:countprofiles}
 N_{\rm up}(10^{23})<8\cdot10^{23},\qquad
 N_{\rm up}(10^{24})<8.4\cdot10^{24}.
\end{equation}
Hence Corollary~\ref{pole:height-corollary} proves the following implications for
each density separately:
\begin{gather*}
 f\in\PF_{10^{24}}\ \Longrightarrow\ \text{RH through height }10^{23},\\
 f\in\PF_{9\cdot10^{24}}\ \Longrightarrow\ \text{RH through height }10^{24}.
\end{gather*}
The checks in \eqref{pole:countprofiles} are repeated at $100$ and $180$
decimal digits with outward Arb arithmetic.

More generally $N_{\rm up}(Y)\sim Y\log Y/(2\pi)$.
The inverse counting scale is $Y\sim2\pi k/\log k$:
for any fixed $0<c<2\pi$, taking
$Y=ck/\log k$ satisfies the sufficient count inequality for all sufficiently
large $k$. Thus $k\asymp T^2$ forces reality of zeros through heights of
order $T^2/\log T$, far beyond $T$. A new analytic proof of that PF
claim would simultaneously prove the stronger zero information.

The present sufficient construction, by comparison, uses the factor budget \eqref{eq:moving-frame-budget},
whose minimal admissible factor count satisfies $n\sim(4+2\sqrt2)k$,
and $n\le N(T)$, giving order of size $T\log T$.
The pole obstruction explains why replacing this count scale
by the squared-rate scale $T^2$ cannot be justified solely as a numerical
optimization of the existing certificate.

\subsection{What a finite spatial sweep proves}\label{app:finite-PF-stitch}
The variable in a spatial sweep of the reciprocal density is distinct
from the ordinate of a zeta zero. We write $L$ for a spatial cutoff and
$Y$ for a zero-height cutoff. There is a valid local-to-global passage
on a fixed spatial interval, but it preserves the certified order and
the domain of differences.

\begin{lemma}[Finite-window differential criterion]\label{local:finite-window}
Let $m\ge1$ and let $I=(a,b)$ be an interval. Suppose
$f\in C^{2m-2}(I)$ and
$\mathscr D_j(f;t)>0$ for $t\in I$ and $1\le j\le m$.
Then, for increasing tuples of size $k\le m$,
\begin{equation}\label{eq:finite-PF-window}
 \det[f(x_i-y_j)]_{i,j=1}^k>0
 \quad\text{if}\quad a<x_1-y_k\quad\text{and}\quad x_k-y_1<b.
\end{equation}
There is also a support-boundary version. Let $L>0$, suppose
$f\in C^{2m-2}(( -\infty,L))$ vanishes on $(-\infty,0]$, and assume
$\mathscr D_j(f;t)>0$ for $0<t<L$, $j\le m$. Then
\begin{equation}\label{eq:finite-PF-support-window}
 \det[f(x_i-y_j)]_{i,j=1}^k\ge0\qquad(x_k-y_1<L),
\end{equation}
and this determinant is strictly positive exactly when
$x_i>y_i$ for every $i$.
\end{lemma}
\begin{proof}
Fix $k\le m$. The functions
$f,-f',\ldots,(-1)^{r-1}f^{(r-1)}$, $r\le k$, have positive
initial Wronskians on $I$, precisely the signs $\mathscr D_j$.
Lemma~\ref{fin-evaluation} therefore gives their positive evaluation
determinants. Put $g_j(t)=f(t+x_{k+1-j})$.
On the common interval $(a-x_1,b-x_k)$ the initial Wronskians of the
$g_j$ are positive: reversing the evaluation arguments and removing
the derivative column signs cancels two factors
$(-1)^{r(r-1)/2}$. Under the conditions in
\eqref{eq:finite-PF-window}, all nodes
$t_i=-y_{k+1-i}$ lie in this common interval.
A second use of Lemma~\ref{fin-evaluation} proves positivity.
The resulting evaluation matrix is the transpose of the desired
matrix with both orders reversed, so its determinant is unchanged.

For the support version use the same functions on
$t<L-x_k$, with support thresholds $a_j=-x_{k+1-j}$.
Their initial Wronskians are positive for
$a_r<t<L-x_k$: the arguments used are positive and less than $L$.
The support-boundary evaluation argument in the proof of
Lemma~\ref{lem:PF-support-boundary} is local on the interval containing
its nodes. Explicitly, if some $t_i\le a_i$, support gives the zero
block forcing the determinant to vanish. Otherwise divide by $g_1$,
subtract consecutive rows, and integrate the derived evaluation
determinant over
$[t_1,t_2]\times\cdots\times[t_{k-1},t_k]$.
Every integration variable stays below $t_k<L-x_k$.
The induction for the derived functions therefore remains inside this
same domain, and is strictly positive on
$\max(t_i,a_{i+1})<v_i<t_{i+1}$.
No argument of $f$ at or above $L$ is used. The strictness condition
$t_i>a_i$ is exactly $x_i>y_i$ after reversing the orders.
\end{proof}

Consequently, charts covering $(0,L)$ can be combined when every chart
certifies all the strict differential signs through the \emph{same}
order $m$ for the \emph{same actual function} $f$.
Their union supplies the hypotheses of
Lemma~\ref{local:finite-window}; its conclusion includes minors whose
differences cross several charts. The bases used to bound those same
signs may differ from chart to chart. Exact coverage, including chart
endpoints, is necessary.

A practical sweep at a fixed order $m$ can therefore combine an analytic
endpoint bound on $(0,a]$ with rigorous bounds on closed intervals
covering $[a,L]$. Every interval must prove all signs through $m$ for the
actual kernel, with certified joins. Positive
scalar normalizations, or lower triangular congruences with positive
squared determinant multipliers, can be chosen separately in each chart:
they preserve the signs of the same actual leading determinants.
For $\Lambda$, the first part of the lemma likewise applies to a cover
of any interval of differences. The Gaussian-transfer theorem used in
this paper assumes global PF order of $\Theta$; a finite $\Theta$
window does not meet that premise, since its variance integration ranges
over all $u>0$. A finite-window transfer would require additional control
of that unbounded integration domain.

For example, the finite-translate argument used in the pole obstruction
can then be applied only where every matrix entry is in the certified
domain. For the support version, a combination
$g(x)=\sum_{j=0}^d c_jf(x+jh)$, $d<m$, has at most $d$ sign changes
on $x<L-dh$. Indeed every collocation matrix there has maximum
cross difference less than $L$, and all its minors are covered by
\eqref{eq:finite-PF-support-window}; the same finite-matrix argument
applies. For the interval version with shifts $y_1<\cdots<y_n$,
the corresponding domain is $a+y_n<x<b+y_1$.
These statements provide no control of the arbitrarily late values
used by the pole obstruction.

\begin{lemma}[Finite-window positivity can precede a pole obstruction]
\label{lem:finite-window-counterexample}
For every integer $m\ge2$ and every finite $L>0$, there is a positive
probability density $f$ on $(0,\infty)$ whose zero extension is
$C^{2m-2}$, which satisfies all conclusions of
\eqref{eq:finite-PF-support-window} through order $m$, but is not
$\PF_m$ globally. Its Laplace transform has real poles,
followed by a nonreal conjugate pair, with arbitrarily small nonzero
imaginary parts, and then further real poles.
\end{lemma}
\begin{proof}
Choose distinct positive rates
\[
 0<q_1<\cdots<q_{m-1}<A<r_1<\cdots<r_{m-1}.
\]
For $b>0$ define the normalized transform
\begin{equation}\label{eq:finite-window-model}
 F_b(s)=
 \prod_{j=1}^{m-1}\frac{q_j}{s+q_j}\,
 \frac{A^2+b^2}{(s+A)^2+b^2}\,
 \prod_{j=1}^{m-1}\frac{r_j}{s+r_j}.
\end{equation}
Its inverse is positive. To see this before convolving the remaining
real exponentials, put $a=A-q_1>0$. The inverse of the first real
factor and the conjugate pair is
\[
 \frac{q_1(A^2+b^2)}{a^2+b^2}e^{-q_1x}
 \left[1-e^{-ax}\left(\cos bx+\frac ab\sin bx\right)\right]
 \quad(x>0).
\]
The bracket starts at zero and has derivative
$(a^2+b^2)e^{-ax}\sin(bx)/b$.
Its positive-time local minima occur at $2\pi\ell/b$, where its value
is $1-e^{-2\pi\ell a/b}>0$. It is therefore positive for every $x>0$.
Convolution with the remaining positive exponential densities preserves
positivity, and $F_b(0)=1$ gives unit mass.

There are $N=2m$ factors, counting the conjugate pair twice.
The positive branch $f_b$ solves the corresponding monic constant-
coefficient ODE of order $N$, whose coefficients depend polynomially
on $b^2$. Its initial derivatives vanish through degree $N-2$, and
\[
 f_b^{(N-1)}(0)=K_b:=(A^2+b^2)\prod_{j=1}^{m-1}q_jr_j>0.
\]
These facts follow also by expanding \eqref{eq:finite-window-model}
at infinity. The companion-matrix solution is analytic jointly in
$x$ and $b^2$. Thus the actual zero extension is $C^{N-2}=C^{2m-2}$.

For $1\le j\le m$, define for $x>0$
\[
 E_j(x,b)=x^{-j(N-j)}\mathscr D_j(f_b;x).
\]
The power series of $f_b$ at zero shows that this expression has a
continuous, indeed analytic, extension in $x$ and $b^2$ at $x=0$.
To check the cancellation explicitly, factor
$x^{N-1-i-\ell}$ from the $(i,\ell)$ derivative entry; the total
power is $j(N-j)$, and the remaining entries are analytic power
series. Their constant determinant is the signed determinant of
the derivatives of $K_bx^{N-1}/(N-1)!$. Consequently
\[
 E_j(0,b)=\frac{K_b^j}{((N-1)!)^j}
              \prod_{a=0}^{j-1}a!(N-1)_{\underline a}>0.
\]
At $b=0$, all $N$ rates are real, with the rate $A$ repeated twice.
The positive representation \eqref{eq:signed-network-representation}
therefore gives $\mathscr D_j(f_0;x)>0$ for every $x>0$ and $j\le m$.
Each $E_j(x,0)$ is positive on the compact interval $[0,L]$,
including at zero by the displayed formula. Joint continuity and the
finite number of indices imply $E_j(x,b)>0$ on this entire interval
for all sufficiently small $b>0$.
Lemma~\ref{local:finite-window} now gives the asserted finite-window
minor signs and their support-dependent strictness.

Nevertheless, the first nonreal rates of $F_b$ are $A\pm ib$.
There are exactly $m-1$ preceding real poles, no tied real pole,
and all these poles are simple. Theorem~\ref{pole:general} therefore gives
$f_b\notin\PF_m$. The extra real rates above $A$ supply the required
boundary regularity without changing that obstruction.
\end{proof}

This example shows why a finite spatial certificate alone
cannot replace the global PF premise in the zero-height implication.
The conjugate rates $A\pm ib$ can be written
$(\gamma_b\pm i\delta_b)^2$, with
$\gamma_b\to\sqrt A$ and $\delta_b\to0$ as $b\downarrow0$.
Thus the obstruction can lie at a bounded zero-height scale while its
failure of PF escapes every prescribed finite spatial window.
A finite-window zero-exclusion argument for the actual densities would
need additional quantitative estimates placing enough alternating
values of the annihilated tail inside the certified window. The
qualitative residue expansion supplies no such uniform window:
nonreal frequencies can approach zero, and the constants and gaps in
the remainder also require control.

\subsection{Fixed-split bounds}\label{app:signed-local}
The following local and bulk estimates give a fixed-split
criterion, which we keep as an independent bound.

\begin{lemma}[Local signed smoothing]\label{lem:signed-local}
Let $n\ge8m$, $0<q_\ell\le Q$, and
$f=*_{\ell=1}^n\Exp(q_\ell)$. Let $\eta$ be a finite real signed
measure supported on $[0,\infty)$, put $p=f*(\delta_0+\eta)$,
and set $\varepsilon=\|\eta\|_{\rm TV}$. Define
\[
 r_m=\sqrt{\frac{m-1}{n-m+1}},\qquad c_m=\frac{r_m}{1-r_m}.
\]
Then $\mathscr D_k(p;x)>0$ for every $k\le m$ and every $x>0$ satisfying
\begin{equation}\label{eq:signed-local-error}
 e^{(1+2c_m)Qx}-1+e^{2c_mQx}\varepsilon<1.
\end{equation}
In particular, on $0<x\le1/(16Q)$ it suffices that
$e^{7/48}-1+e^{1/12}\varepsilon<1$.
\end{lemma}
\begin{proof}
Put $N=n-1$, $C=\prod q_\ell$, and $h(x)=x_+^N/N!$.
At a fixed $x>0$, the normalized polynomial functional
$\mathcal L(z^j)=x^jh^{(j)}(x)/h(x)=N_{\underline j}$
has signed orthonormal polynomials
\begin{equation}\label{eq:signed-Laguerre-basis}
 Q_i(z)=\frac{L_i^{(-n)}(-z)}{\sqrt{\binom Ni}},\qquad
 \mathcal L(Q_iQ_j)=(-1)^i\delta_{ij},\quad 0\le i,j<m.
\end{equation}
These are finite polynomial identities. For completeness, the Laguerre
generating function and $\mathcal L(e^{tz})=(1+t)^N$ give
\[
 \mathcal L\!\left[
 (1-u)^N e^{zu/(1-u)}(1-v)^N e^{zv/(1-v)}\right]
                         =(1-uv)^N;
\]
coefficient comparison proves \eqref{eq:signed-Laguerre-basis}.
All polynomial operator coefficients
below are frozen at the chosen $x$; $z=xD$ is not an Euler operator
acting on those coefficients.

Let $L$ represent polynomial differentiation, $Q'(z)=LQ(z)$.
Its strictly lower triangular entries are
\[
 L_{ij}=\sqrt{\frac{i!\,N_{\underline j}}
                         {j!\,N_{\underline i}}}\le r_m^{i-j}
       \quad(j<i).
\]
Both row and column sums are at most $c_m$, so
$\|L\|_2\le c_m$. Hence, for $v\ge0$,
\begin{equation}\label{eq:signed-Laguerre-translation}
 Q(z-v)=T(v)Q(z),\qquad T(v)=e^{-vL},\qquad
                         \|T(v)\|_2\le e^{c_mv}.
\end{equation}
The polynomial Euler operator $zQ'(z)=AQ(z)$ is instead bidiagonal:
$A_{ii}=i$, $A_{i,i-1}=\sqrt{i(n-i)}$. Thus
\[
 \lambda_{\max}((A+A^T)/2)
 \le(m-1)+\sqrt{(m-1)N}<N/2,
\]
because $N\ge8(m-1)$ and $1/8+1/\sqrt8<1/2$; for $m=1$
the assertion is immediate. Differentiating squared Euclidean norms gives
\begin{equation}\label{eq:signed-dilation-contract}
                   \|e^{t(A-NI/2)}\|_2\le1\quad(t\ge0).
\end{equation}

Let $J_m=\operatorname{diag}((-1)^i)_{i=0}^{m-1}$. For a shift
$0\le u<x$, put $r=(x-u)/x$ and $t=-\log r$. The normalized
Gram matrix of the shifted monomial is exactly
\begin{align}
 G(r)&=\left[\frac{Q_i(xD)Q_j(xD)h(x-u)}{h(x)}\right]_{ij}
      =r^N e^{tA}J_m e^{tA^T},\nonumber\\
 \|G(r)\|_2&\le1.\label{eq:signed-shift-contraction}
\end{align}
Indeed its moments are $N_{\underline j}r^{N-j}$, which replaces
each polynomial argument by $z/r$ and multiplies by $r^N$.
For $u\ge x$ its actual zero-extended derivative matrix is zero,
including $u=x$, since $2m-2\le n-2$.

The convolution-simplex formula reads
\[
 f(y)=C h(y)\,\E e^{-\mathcal A y},\qquad
 \mathcal A=\sum q_\ell V_\ell\in[0,Q],
\]
where $V$ is uniform on the standard simplex and its law is independent
of $y$. Normalize derivative Grams by the fixed basis $Q_i(xD)$ and
the positive scalar $Ch(x)$. For a fixed value $a\in[0,Q]$, the
unshifted Gram of $Ch(y)e^{-ay}$ is
$e^{-ax}T(ax)J_mT(ax)^T$. The matrix operator
$X\mapsto-X-LX-XL^T$ has norm at most $1+2c_m$, so its
difference from $J_m$ has norm at most $e^{(1+2c_m)ax}-1$.
Averaging gives the first term of \eqref{eq:signed-local-error}.

For $0\le u<x$, the corresponding shifted Gram is
\[
 e^{-a(x-u)}T(ax)G((x-u)/x)T(ax)^T.
\]
Its norm is at most $e^{2c_mQx}$ by
\eqref{eq:signed-Laguerre-translation}--\eqref{eq:signed-shift-contraction}.
The order of these matrices matters.
For $u\ge x$ the Gram is zero. Integrating against $\eta$ proves
the second term of \eqref{eq:signed-local-error}.
Differentiation under that integral is legitimate: the zero-extended
$f$ has bounded continuous derivatives through order $n-2$.

An error of norm less than one from $J_m$ preserves the inertia of
every leading block along the straight-line homotopy. The basis change
is lower triangular with positive diagonal, hence preserves all leading
determinant signs. Finally $n\ge8m$ gives $r_m<2/5$ and $c_m<2/3$,
which proves the stated uniform specialization.
\end{proof}

\begin{lemma}[Bulk signed smoothing]\label{lem:signed-bulk}
Let $n\ge2m$, $f=*_{\ell=1}^n\Exp(q_\ell)$ with $0<q_\ell\le Q$,
and $p=f*(\delta_0+\eta)$ for a finite real signed measure supported
on $[0,\infty)$. If $x_0>0$ and
\[
 a_0=2n/x_0+nQ,\qquad \int e^{a_0u}|\dd\eta(u)|<1,
\]
then $\mathscr D_k(p;x)>0$ for every $x\ge x_0$ and $k\le m$.
\end{lemma}
\begin{proof}
At a fixed $x\ge x_0$, use the full normalization of
Lemma~\ref{lem:signed-Jacobi}. In its leading $m$-block, the difference
between the actual Grams of $p$ and $f$ is the integral of the
actual shifted Gram of $f$ against $\eta$. For $u>x$ this is zero.
For $0\le u\le x$ it is the leading block of $e^{-uB_x}J_n$.
At $u=x$ the leading block is zero as well, because only derivatives
through $2m-2\le n-2$ are involved. Thus the error norm is at most
\[
 \int_{[0,x]}e^{u\|B_x\|_2}|\dd\eta(u)|
             \le\int e^{a_0u}|\dd\eta(u)|<1.
\]
Compression and the signed-inertia argument in
Lemma~\ref{lem:signed-local} prove every leading determinant sign.
There is no crossing error: the comparison is directly with the unit
atom, and the actual shifted kernel is zero beyond its support.
\end{proof}

\subsection{Explicit PF convolution factors}\label{sec:kps-factors}
For this factorization we use a variant of Lemma~\ref{lem:anchor}.
Its bounds hold with $q$ replaced by any real $q_*$ satisfying
$2\le q_*\le U^2/4-1/4$, and with
$K_U(t)=e^{-q_*t}+H_{>U}(t)$: the same proof applies because
$a_\rho-q_*\ge3\gamma^2/4$ and $|\lambda_\rho|\le q_*\gamma^2$.

For $100\le U\le U_{\rm PT}$, choose a verified ordinate $\gamma_*$ with
$2\le\gamma_*^2\le U^2/4-1/4$; the first ordinate is always a possible
choice. Let $\mathcal A_U$ be the finite multiset of ordinates
$\gamma\le U$, with \emph{one occurrence} of $\gamma_*$ deleted, and define
\[
 P_U(s)=\prod_{\gamma\in\mathcal A_U}(1+s/\gamma^2).
\]
Repeated zeros are retained with their multiplicities.

\begin{theorem}[Extraction of verified-zero PF factors]\label{thm:PF-factors}
There is a positive self-decomposable random variable $R_U$ such that
\begin{align}
 S_1&\overset d=\sum_{\gamma\in\mathcal A_U}E_\gamma+R_U,
       &E_\gamma&\sim\Exp(\gamma^2),\label{eq:exp-factors}\\
 X_{1,1}&\overset d=\sum_{\gamma\in\mathcal A_U}L_\gamma+Y_U,
       &L_\gamma&\sim\Lap(\gamma),\label{eq:lap-factors}
\end{align}
where $\Lap(\gamma)$ has density $\gamma e^{-\gamma|x|}/2$ and $Y_U\overset d=\sqrt{2R_U}\,G_1$ with $G_1$ independent of $R_U$. The summands on each right-hand side are independent. The density of either finite sum is a $\PF_\infty$ function. The residual $R_U$ has a smooth positive density on $(0,\infty)$; $Y_U$ is self-decomposable and has an entire positive density.
\end{theorem}
Under RH, $1/h(s)=\prod_\gamma(1+s/\gamma^2)^{-1}$ is the Laplace
transform of a convergent sum of independent exponential variables, a
generalized gamma convolution in Thorin's sense
\citep[Chapter~3]{bondesson1992generalized}, and correspondingly
$X_{1,1}=\sum_\gamma L_\gamma$; for xi these factorizations are written in
\citet[Section~2.1]{polson2018vandantzig} and
\citet[v9, Theorem~12(i)]{polson2026duality}. The point here is the
unconditional extraction of finitely many verified factors with a
self-decomposable residual.
\begin{proof}
Set $C_U(s)=B(s)-\log P_U(s)$. By the product and integral formulas,
\begin{equation}\label{eq:residual-exponent}
 C_U(s)=\int_0^\infty(1-e^{-st})\frac{K_U(t)}t\dd t,
 \qquad \E e^{-sR_U}=e^{-C_U(s)}=\frac{P_U(s)}{h(s)}.
\end{equation}
Here $K_U(t)=e^{-\gamma_*^2t}+H_{>U}(t)$, which is positive and decreasing
by the variant of Lemma~\ref{lem:anchor} stated above. Thus the Poisson and self-decomposability constructions from Section~\ref{sec:prob} apply. Multiplication of Laplace transforms proves \eqref{eq:exp-factors}. On the Fourier axis each factor becomes $(1+t^2/\gamma^2)^{-1}$, the characteristic function of $\Lap(\gamma)$, proving \eqref{eq:lap-factors}.
The finite exponential and centered Laplace densities are $\PF_\infty$
by Schoenberg's product characterization \citep{schoenberg1951polya}.
Equivalently, their transform reciprocals have only real zeros of the
required sign. Removing finitely many exponentials leaves the small-time divergence of $H$ unchanged and changes $B(s)$ by $O_U(\log(2+s))$. The density and entire-extension arguments in Section~\ref{sec:prob} therefore apply to both residuals.
\end{proof}
The densities in \eqref{eq:ThetaLambda} therefore have explicit finite
$\PF_\infty$ convolution factors. A higher retained zero makes the
residual small and yields quantitative approximations.

\begin{corr}[Approximation by verified-zero PF laws]\label{cor:PF-approximation}
For $1000\le U\le\Tstar$, the retained ordinate can be chosen in
$(U/4,U/3]$. Put
\[
 b_U=\frac{16}{U^2}+\frac{2(\log U+1)}U,
 \qquad A_U=\sum_{\gamma\in\mathcal A_U}E_\gamma,
 \qquad Z_U=\sum_{\gamma\in\mathcal A_U}L_\gamma.
\]
Then $\E R_U\le b_U$. If $d_{\mathrm{TV}}$ denotes total variation distance
between probability laws, with values in $[0,1]$, then
\begin{equation}\label{eq:PF-approximation}
 d_{\mathrm{TV}}(S_1,A_U)\le q b_U,\qquad
 d_{\mathrm{TV}}(X_{1,1},Z_U)\le q b_U,\qquad q=\gamma_1^2<225.
\end{equation}
Both bounds are less than $4.5\cdot10^{-9}$ at $U=\Tstar$.
The approximating laws have $\PF_\infty$ densities.
\end{corr}
\begin{proof}
The counting bound gives
\[
 N(U/3)-N(U/4)
 \ge\frac{U}{24\pi}\log\frac{U}{8\pi}
       -\frac14\log(U+2)-8>0\qquad(U\ge1000).
\]
Indeed, the first term exceeds $U/32$; the resulting lower bound is
positive at $1000$ and increasing thereafter. This supplies $\gamma_*$,
and its square satisfies the hypothesis of Lemma~\ref{lem:anchor}.
Differentiating \eqref{eq:residual-exponent} at zero gives
\[
 \E R_U=\gamma_*^{-2}+\sum_{\gamma>U}\real\frac1{\lambda_\rho}
 \le\frac{16}{U^2}+\sum_{\gamma>U}\gamma^{-2}\le b_U,
\]
since $\real(\lambda_\rho^{-1})
=(\gamma^2-\delta_\rho^2)/(\gamma^2+\delta_\rho^2)^2\le\gamma^{-2}$.

The finite sums include the first-zero factors. Translation by $r\ge0$
changes the law $\Exp(q)$ in total variation by $1-e^{-qr}\le qr$.
Convolution contracts total variation, so conditioning on $R_U$ proves
the first bound in \eqref{eq:PF-approximation}.
For the second, let $l(x)=\gamma_1e^{-\gamma_1|x|}/2$. Its second
distributional derivative is $D^2l=q l-q\delta_0$, of total variation
$2q$. Consequently
\[
 \|l(\cdot-y)+l(\cdot+y)-2l\|_1\le2q y^2.
\]
This follows by integrating translates of $D^2l$ against the triangular
weight $(|y|-|v|)_+$. The residual $Y_U$ is symmetric and
$\E Y_U^2=2\E R_U$. Averaging the displayed inequality and dividing
the resulting $L^1$ bound by two gives
$d_{\mathrm{TV}}(L_{\gamma_1}+Y_U,L_{\gamma_1})\le q\E R_U$.
Convolving the other Laplace factors proves the claim. The final
numerical bound follows from $q<225$ and $U=\Tstar$.
\end{proof}
These estimates concern the normalized densities $c_0\Theta$ and
$c_0\Lambda$.

\subsection{Residual moments for the higher-order certificate}
\label{sec:PF-residual-moment}
Put $U=\Tstar$ and retain $\gamma_*\in(U/4,U/3]$ as in
Corollary~\ref{cor:PF-approximation}. Write $R=R_U$ and
$Y=\sqrt{2R}G_1$. For order $256$, we use sharper bounds that retain
the common residual mean across all matrix entries. The counting
estimate supplies the following enclosures.
Set
\begin{align*}
 T&=\frac{\log(U/(2\pi))+1}{2\pi U},&
 C&=\frac{\log(U+2)/4+129/16}{U^2},&D_0&=T+C,\\
 B&=D_0+\frac{16}{U^2},&
 V_2&=\frac{D_0}{U^2}+\frac{256}{U^4},&
 A(s)&=1-\frac{16(s+1/4)}{U^2}.
\end{align*}
\begin{lemma}\label{lem:PF-residual-moment}
The residual mean $\mu=\E R$ and variance satisfy
\begin{equation}\label{eq:PF256-moments}
 \mu_-:=T-C+\frac9{U^2}-\frac{3D_0}{4U^2}
 \le\mu\le\mu_+:=B,\qquad \operatorname{Var}(R)\le V_2.
\end{equation}
For $s\ge0$ with $A(s)>0$,
\begin{equation}\label{eq:PF256-mgf}
 \log\E e^{sR}\le\frac{Bs}{A(s)},\qquad
 0\le\log\E e^{sR}-\mu s\le\frac{s^2V_2}{2A(s)^2}.
\end{equation}
In particular, with $\Omega=3\cdot10^7$ and
$\kappa=\exp(B\Omega^2/A(\Omega^2))-1$,
$\E\cosh(\Omega Y)-1\le\kappa$.
\end{lemma}
\begin{proof}
Stieltjes integration of the counting error, including its boundary
term at $U$, gives
$\sum_{\gamma>U}\gamma^{-2}\in[T-C,T+C]$.
For $\lambda=(\gamma+i\delta)^2$, $|\delta|\le1/2$,
\[
 0\le\gamma^{-2}-\real\lambda^{-1}\le\frac3{4\gamma^4}.
\]
Combine this with $9/U^2\le\gamma_*^{-2}\le16/U^2$ in
\eqref{eq:residual-exponent} to obtain the mean interval.
The second logarithmic derivative at zero is bounded by
$\gamma_*^{-4}+\sum_{\gamma>U}|\lambda|^{-2}\le V_2$.
For $0\le v\le s$, each resolvent denominator in the first two
logarithmic derivatives loses at most a factor $A(s)$ relative to its
value at zero. Equivalently, integrate the positive residual kernel
against $e^{vt}$ and bound each zero term by its absolute value.
The derivatives are bounded by $B/A(s)$ and $V_2/A(s)^2$,
respectively. Integration and convexity prove
\eqref{eq:PF256-mgf}. Gaussian conditioning gives the cosh bound.
All domain inequalities are checked with directed arithmetic.
\end{proof}

Higher cumulants retain the same dependence. If
$h_j=\kappa_j(R)/(j-1)!$, then
\[
 h_j=\gamma_*^{-2j}+\sum_{\gamma>U}\real\lambda_\rho^{-j}.
\]
Writing $q=2j-1$, put
\[
 T_j=\frac{\log(U/(2\pi))+1/q}{2\pi qU^q},\qquad
 C_j=\frac{\log(U+2)/4+8+1/(16j)}{U^{2j}}.
\]
The same integration and the inequality
$0\le\gamma^{-2j}-\real\lambda^{-j}
\le j(2j+1)\gamma^{-2j-2}/4$ yield
\begin{equation}\label{eq:PF256-cumulants}
 T_j-C_j+\frac{9^j}{U^{2j}}-\frac{j(2j+1)D_0}{4U^{2j}}
 \le h_j\le T_j+C_j+\frac{16^j}{U^{2j}}.
\end{equation}
These bounds use only the published verified height and the counting
estimate.

\subsection{Differential determinants and strict composition}
\label{sec:Lambda-PF6}
For a sufficiently smooth density $p$, define
\[
 \mathscr D_k(p;x)=(-1)^{k(k-1)/2}
       \det[p^{(i+j)}(x)]_{i,j=0}^{k-1}.
\]
\begin{lemma}[Differential criterion for PF minors]
\label{lem:PF-differential}
If $p\in C^{2m-2}(\reals)$ and $\mathscr D_k(p;x)>0$ for every
$x\in\reals$ and $1\le k\le m$, then $p(x-y)$ has strictly positive
ordered minors through order $m$.
\end{lemma}
\begin{proof}
Fix $y$. The initial Wronskians of the functions
$x\mapsto(-1)^jp^{(j)}(x-y)$, $0\le j<m$, are precisely
$\mathscr D_k(p;x-y)$. Lemma~\ref{fin-evaluation} makes their ordered
evaluation determinants positive. For fixed increasing $x_i$, these
determinants are the transposes of the initial Wronskians in $y$ of
$y\mapsto p(x_i-y)$. A second application of that lemma proves the claim.
\end{proof}

Choose the $270$ distinct verified ordinates $a_0<\cdots<a_{269}<U/4$
specified in Appendix~\ref{app:PF256-certificate}, and let
$f=\Lap(a_0)*\cdots*\Lap(a_{269})$ denote their convolution density.
Theorem~\ref{thm:PF-factors} gives
\begin{equation}\label{eq:PF256-factorization}
 c_0\Lambda=p*\nu,\qquad p=\E H_Rf,\qquad H_t=e^{tD^2},
\end{equation}
where $\nu$ is the convolution of the remaining extracted Laplace
laws. The positive Fourier transform
$\widehat f(t)=\prod_l a_l^2/(a_l^2+t^2)$ has integrable moments
through degree $538$, so $f$ and $p$ have bounded continuous
derivatives through that degree. This exceeds the $510$ derivatives
needed for $m=256$.

\begin{lemma}[Certified differential signs]\label{lem:PF256-signs}
For the density $p$ in \eqref{eq:PF256-factorization},
$\mathscr D_k(p;x)>0$ for every $x\in\reals$ and
$1\le k\le256$.
\end{lemma}
\begin{proof}
The interval calculation in Appendix~\ref{app:PF256-certificate}
certifies all leading determinant signs uniformly on $[0,3]$,
using the residual bounds just proved. Its analytic tail estimate
certifies the same signs on $[3,\infty)$. Evenness gives the
negative half-line: the derivative matrices at $-x$ and $x$ are
congruent by $\operatorname{diag}((-1)^i)$.
\end{proof}

Lemma~\ref{lem:PF-differential} now gives strict ordered minors
through size $256$ for $p$. Convolution with each remaining Laplace
density preserves strictness. The continuous Cauchy--Binet formula
\citep{karlin1968} integrates the product of a strictly positive
minor and a nonnegative Laplace minor over $u_1<\cdots<u_k$.
The latter is positive near $u_i=y_i$: the matrix
$[e^{-a|y_i-y_j|}]$ is positive definite for distinct $y_i$, by its
strictly positive Fourier density. Thus the integral is strictly
positive. Iterating over the remaining finite set of factors in
$\nu$ proves the assertion for $c_0\Lambda$, and hence for $\Lambda$.

\subsection{Fixed-order positivity near the one-sided boundary}
\label{sec:Theta-local-PF}
We use the exponential-tilting (saddlepoint) method of
\citet{daniels1954saddlepoint} to obtain relative determinant estimates
near the left endpoint of $\Theta$.
Under a sufficiently strong negative
exponential tilt, the one-sided law becomes Gaussian after centering
and scaling. To control determinants, weak convergence is insufficient:
we prove convergence of every fixed number of density derivatives.

\begin{prop}[Local PF positivity of $\Theta$]
\label{prop:Theta-local-PF}
For each fixed integer $m\ge1$, there is $\varepsilon_m>0$ such that
\[
 \det[\Theta(x_i-y_j)]_{i,j=1}^k>0\qquad(1\le k\le m)
\]
whenever $x_1<\cdots<x_k$, $y_1<\cdots<y_k$, and
$0<x_i-y_j<\varepsilon_m$ for every $i,j$.
More precisely, put $x_s=B'(s)$ and $v_s=-B''(s)$. As $s\to\infty$,
\begin{equation}\label{eq:Theta-local-determinant}
 \mathscr D_k(\Theta;x_s)
 \sim \Theta(x_s)^k v_s^{-k(k-1)/2}\prod_{j=0}^{k-1}j!
 \qquad\text{for every fixed }k.
\end{equation}
\end{prop}
\begin{proof}
We use the branch of $B(z)=\log h(z)$ that is real on the positive
axis. Stirling's expansion \citep[Section~5.11]{dlmf}, together with
$\log\zeta(w)=O(2^{-\real w})$ for $\real w\to\infty$, gives
\begin{equation}\label{eq:Theta-B-sector}
 B(z)=\frac{\sqrt z}{2}\left(\log\frac{\sqrt z}{2\pi}-1\right)
       +\frac78\log z+C+O(|z|^{-1/2})
\end{equation}
uniformly in the right half-plane as $|z|\to\infty$. Indeed,
$|\arg\sqrt z|\le\pi/4$ and
$\real\sqrt z\ge\sqrt{|z|/2}$ there, so the gamma expansion and
Euler-product estimate are uniform. Cauchy's estimates permit
termwise differentiation near the positive axis. In particular,
\begin{equation}\label{eq:Theta-saddle-scales}
 x_s\sim\frac{\log s}{8\sqrt s},\qquad
 v_s\sim\frac{\log s}{16s^{3/2}},\qquad
 B'''(z)=O\!\left(\frac{\log s}{s^{5/2}}\right)
 \quad(|z-s|\le s/2).
\end{equation}
For large $s$, $v_s>0$ and $x_s$ decreases continuously to zero.

The tilted density is $f_s(u)=e^{B(s)-su}\theta_1(u)$, with mean $x_s$
and variance $v_s$. Write $\sigma_s=\sqrt{v_s}$ and let
$g_s(z)=\sigma_s f_s(x_s+\sigma_s z)$ be its centered, standardized
density. Its characteristic function is
\[
 \phi_s(t)=\exp\left\{B(s)-B(s-it/\sigma_s)-itx_s/\sigma_s\right\}.
\]
Taylor expansion using \eqref{eq:Theta-saddle-scales} gives
$\phi_s(t)\to e^{-t^2/2}$ for each fixed $t$.

To pass to density derivatives, we need an integrable bound for all
frequencies. Write $\sqrt{1+iu}=a+ib$ with $a\ge1$, and set
$w=\sqrt s$, $L=\log(w/(2\pi))-1$. From
\eqref{eq:Theta-B-sector},
\[
 \real B(s+isu)-B(s)
 \ge\frac{wL}{2}(a-1)-\frac{C_1}{w}.
\]
Here $\real[(a+ib)\log(a+ib)]\ge0$: for $\eta=|b|/a<1$, its quotient
by $a$ is $\operatorname{arctanh}(\eta^2)-\eta\arctan\eta$, whose derivative
is $\eta/(1-\eta^2)-\arctan\eta\ge0$. The logarithmic term
$\tfrac78\log|1+iu|$ is also nonnegative. Finally,
\[
 a-1=\frac{u^2}{2(a+1)(\sqrt{1+u^2}+1)}
 \ge c_1\min\{u^2,\sqrt{|u|}\}
\]
for an absolute $c_1>0$.
Put $A_s=s\sigma_s$, so $A_s^2\sim wL/8$ and $A_s\to\infty$.
Substitution of $u=t/A_s$ gives
$|\phi_s(t)|\le C_2\exp\{-c\min(t^2,A_s^{3/2}\sqrt{|t|})\}$,
with $c,C_2>0$ independent of large $s$. Since $A_s\ge1$ eventually,
\begin{equation*}
 |\phi_s(t)|\le C_2\exp\{-c\min(t^2,\sqrt{|t|})\}.
\end{equation*}
The right side is integrable after multiplication by any fixed power
of $|t|$. Dominated convergence in Fourier inversion therefore gives
$g_s^{(r)}(0)\to\varphi^{(r)}(0)$ for every fixed $r$, where
$\varphi(z)=(2\pi)^{-1/2}e^{-z^2/2}$.

The Gaussian moment determinant is
\[
 \mathscr D_k(\varphi;0)=\varphi(0)^k\prod_{j=0}^{k-1}j!>0.
\]
For example, this follows by replacing monomials with monic Hermite
polynomials in the Gaussian Gram matrix; their squared norms are $j!$.
For $M_f=(f^{(i+j)}(x))_{i,j=0}^{k-1}$, Leibniz's rule gives
$M_{e^{cx}f}=e^{cx}T_cM_fT_c^\mathsf T$, where $T_c$ is lower triangular
with $(T_c)_{ij}=\binom ij c^{i-j}$ for $j\le i$. Since $\det T_c=1$,
tilting multiplies the differential determinant by $e^{kcx}$.
Scaling then gives
\[
 \frac{\mathscr D_k(\Theta;x_s)}{\Theta(x_s)^k}
 =\sigma_s^{-k(k-1)}
   \frac{\mathscr D_k(g_s;0)}{g_s(0)^k},
\]
which proves \eqref{eq:Theta-local-determinant}. All initial signed
determinants through size $m$ are therefore positive on
$(0,\varepsilon_m)$ for some $\varepsilon_m>0$.
The proof of Lemma~\ref{lem:PF-differential} applies locally: both
Wronskian integrations remain in this interval when all the cross
differences do. This proves the stated $\PF_k$ minor inequalities.
\end{proof}

\subsection{The global order-four certificate for \texorpdfstring{$\Theta$}{Theta}}
\label{app:Theta-PF4-certificate}

For global PF positivity, the differential criterion must also handle
nonpositive cross differences, as in
Lemma~\ref{lem:PF-support-boundary}.

For $\Theta$, the required interior signs follow from uniform analytic
bounds near zero and infinity, joined by closed interval computations.
\begin{lemma}\label{lem:Theta-PF4-signs}
For every $x>0$ and $1\le k\le4$, $\mathscr D_k(\Theta;x)>0$.
\end{lemma}
\begin{proof}
Write $p=c_0\Theta=\theta_1$. The positive constant $c_0$ does not affect
the signs. The three domains used below are
\begin{equation}\label{eq:Theta-PF4-domains}
 (0,0.00211005],\qquad [0.00211004,0.01001],\qquad [0.008,\infty).
\end{equation}
All terminating decimals in this proof denote exact rationals.

\paragraph{The uniform interval at zero.}
We make the bounds in the proof of Proposition~\ref{prop:Theta-local-PF}
explicit. With principal square roots and logarithms, put
\[
 \mathcal F(z)=\frac{\sqrt z}{2}
       \left(\log\frac{\sqrt z}{2\pi}-1\right)+\frac78\log z.
\]
For $\real z\ge0$ and $|z|\ge10^4$, the remainder in
\eqref{eq:Theta-B-sector} satisfies
\begin{equation}\label{eq:Theta-PF4-sector-remainder}
 B(z)=\mathcal F(z)+C+E(z),\qquad |E(z)|\le |z|^{-1/2}.
\end{equation}
To check the constant, set $w=\sqrt z$ and $v=w/2+1/4$.
Binet's formula \citep[Eq.~5.9.10\_2]{dlmf} gives
\begin{align*}
 E(z)&=\log\left(1-\frac1{4w^2}\right)
   +\left(\frac w2-\frac14\right)\log\left(1+\frac1{2w}\right)
   -\frac14+R(v)+\log\zeta(w+\tfrac12),\\
 R(v)&=\int_0^\infty e^{-vt}
          \left(\frac12\coth\frac t2-\frac1t\right)\frac{\dd t}{t}.
\end{align*}
The partial-fraction expansion of $\coth$ bounds the bracket by $t/12$.
Thus $|R(v)|\le\sqrt2/(6|w|)$.
For $|w|\ge100$, the first logarithm is at most $(4|w|^2-1)^{-1}$
in modulus, and the second term together with $-1/4$ is at most
$3/[16|w|(1-1/(2|w|))]$. The Euler product and integral bound for
the zeta tail give, with $\sigma=\real(w+1/2)$,
\[
 |\log\zeta(w+\tfrac12)|\le\log\zeta(\sigma)
 \le 2^{-\sigma}\left(1+\frac2{\sigma-1}\right)
 \le\frac1{100|w|}.
\]
Here $\sigma\ge |w|/\sqrt2+1/2$; the last bound holds at $|w|=100$
and improves thereafter. The sum of these bounds is less than
$0.44/|w|$, proving \eqref{eq:Theta-PF4-sector-remainder}.

Retain $x_s=B'(s)$, $v_s=-B''(s)$, $\sigma_s=\sqrt{v_s}$ and the
standardized tilted density $g_s$ from the preceding subsection.
Put $w=\sqrt s$, $L=\log(w/(2\pi))-1$, $W=wL$, and $A_s=s\sigma_s$.
Cauchy's estimate on the disk of radius $s/2$ gives
$|E^{(r)}(s)|\le r!2^r\sqrt2\,s^{-r-1/2}$ and hence
\[
 A_s^2\in\left[\frac W8+\frac78-\frac{8\sqrt2}{w},
                \frac W8+\frac78+\frac{8\sqrt2}{w}\right].
\]
For $s\ge S:=275789$, these bounds imply $W/8\le A_s^2\le W/7$.
Taylor's formula for the characteristic function $\phi_s$ gives, when
$|t|\le T:=3113/500\le A_s/2$,
\[
 \log\phi_s(t)=-t^2/2+i b_s t^3+R_s(t),
 \qquad |b_s|\le a,\qquad |R_s(t)|\le d t^4,
\]
where the following expressions are evaluated at $s=S$:
\begin{align*}
 a&=\frac16\left(\frac8L\right)^{3/2}
     \left(\frac{3L-1}{16\sqrt w}+
       \frac7{4w^{3/2}}+\frac{48\sqrt2}{w^{5/2}}\right),\\
 d&=\frac1{24}\left(\frac8L\right)^2
     \left(\frac{15L+167/4}{32w}+\frac{21}{4w^2}
                         +\frac{384\sqrt2}{w^3}\right).
\end{align*}
These expressions decrease for $s\ge S$. They follow by differentiating
$\mathcal F$ three and four times and applying the same Cauchy bound;
on $s+i[-s/2,s/2]$ one uses $|\log(z/s)|\le1/2$.
With $c=1/2-dT^2>0$, the even part therefore satisfies
\[
 |\real\phi_s(t)-e^{-t^2/2}|
 \le \tfrac12 a^2t^6e^{-t^2/2}+dt^4e^{-ct^2}
 \qquad (|t|\le T).
\]

The sector estimate also bounds the entire Fourier tail.
The calculation preceding \eqref{eq:Theta-local-determinant} gives
\[
 |\phi_s(t)|\le e^{2/\sqrt s}
 \begin{cases}
  e^{-7t^2/24},& |t|\le A_s,\\
  \exp\{-\tfrac W{32}\sqrt{|t|/A_s}\},& |t|\ge A_s.
 \end{cases}
\]
Indeed, for $\alpha=\real\sqrt{1+iu}$, the identity used there gives
$\alpha-1\ge u^2/12$ when $|u|\le1$ and
$\alpha-1\ge\sqrt{|u|}/16$ when $|u|\ge1$.
Define
\[
 I_j(b,T)=\int_T^\infty t^je^{-bt^2}\dd t,
 \qquad
 J_j(W)=2(W/7)^{(j+1)/2}
        \frac{\Gamma(2j+2,W/32)}{(W/32)^{2j+2}}.
\]
The second part of the tail contributes at most $J_j(W)$.
For $0\le j\le6$, its finite incomplete-gamma expansion is decreasing
when $W\ge80$; also $W(S)>\max(80,32T^2)$.
Fourier inversion thus proves, for even $j\le6$ and every $s\ge S$,
\begin{align}
 |g_s^{(j)}(0)-\varphi^{(j)}(0)|&\le e_j,
 \qquad \varphi(z)=(2\pi)^{-1/2}e^{-z^2/2},\label{eq:Theta-PF4-even-errors}\\
 \pi e_j&=\tfrac12a^2 I_{j+6}(1/2,0)+dI_{j+4}(c,0)
     +I_j(1/2,T)\notag\\
 &\hspace{12mm}+e^{2/\sqrt S}
       \{I_j(7/24,T)+J_j(W(S))\}.\notag
\end{align}
The four Gaussian values are
$\varphi^{(j)}(0)=(2\pi)^{-1/2}(1,-1,3,-15)$ for $j=0,2,4,6$.
Outward evaluation of \eqref{eq:Theta-PF4-even-errors} proves the four
strict inequalities
\[
 g_s(0)>0,\qquad -g_s''(0)>0,\qquad
 g_s^{(4)}(0)-\frac{g_s''(0)^2}{g_s(0)}>0,\qquad
 -g_s^{(6)}(0)-\frac{g_s^{(4)}(0)^2}{-g_s''(0)}>0.
\]
Their smallest certified lower bound exceeds $1.5585\cdot10^{-6}$.
The even-index block of the derivative matrix is positive definite and
the odd-index block is negative definite. Eliminating the even block
leaves a negative definite Schur complement, whatever the mixed entries;
the same argument applies to every leading block. Hence
$\mathscr D_k(g_s;0)>0$ for $k\le4$.
The scaling identity in the preceding proof transfers these signs to
$\Theta$ at $x_s$. Finally,
\[
 x_S\ge\frac{L(S)+1}{4\sqrt S}+\frac7{8S}
           -\frac{2\sqrt2}{S^{3/2}}>0.00211005.
\]
Since $x_s'=-v_s<0$ and $x_s\to0$, this proves all signs on the first
domain of \eqref{eq:Theta-PF4-domains}.

\paragraph{The compact interval.}
For positive rational $s,\sigma$ and rational $x$, use the exact density
and characteristic function
\[
 g(z)=\sigma e^{B(s)-s(x+\sigma z)}p(x+\sigma z),
 \qquad
 \phi(t)=\frac{h(s)}{h(s-it/\sigma)}e^{-itx/\sigma}.
\]
At $U=3\cdot10^{12}$, Theorem~\ref{thm:PF-factors} and the retained
anchor from Corollary~\ref{cor:PF-approximation} extract the first
$64$ exponential rates $q_\ell=\gamma_\ell^2$ and leave a positive
independent residual law. Tilting changes these rates to $s+q_\ell$.
Consequently, with $r_\ell=(s+q_\ell)\sigma$,
\[
 |\phi(t)|\le P(t):=\prod_{\ell=1}^{64}(1+t^2/r_\ell^2)^{-1/2}.
\]
For $t\ge T_0:=64$, logarithmic differentiation gives
$P(t)\le P(T_0)(t/T_0)^{-\alpha}$, where
$\alpha=\sum_\ell T_0^2/(T_0^2+r_\ell^2)$. Thus, whenever
$\alpha>j+1$,
\begin{equation}\label{eq:Theta-PF4-Fourier-tail}
 \int_{T_0}^\infty t^j|\phi(t)|\dd t
 \le \frac{P(T_0)T_0^{j+1}}{\alpha-j-1}.
\end{equation}
Certified quadrature of
$\pi^{-1}\real\int_0^{T_0}(-it)^j\phi(t)\dd t$, together with
\eqref{eq:Theta-PF4-Fourier-tail}, encloses $v_j=g^{(j)}(0)$ for
$0\le j\le9$. The transform is evaluated by the exact gamma--zeta
formula \eqref{eq:normalization}.

Absolute Fourier integrals give bounds
$N_j\ge\sup_z|g^{(j)}(z)|$ for $4\le j\le10$.
On each interval $[u,v]\subset[0,64]$ of length $1/8$, use
$t^jP(t)\le v^jP(u)$ and add \eqref{eq:Theta-PF4-Fourier-tail}.
For $s\ge10^4$, the smaller of this bound and the decreasing majorant
\[
 \exp\left\{\frac2{\sqrt s}-\frac W2
 \left(\sqrt{\frac{\sqrt{1+(t/(s\sigma))^2}+1}{2}}-1\right)\right\}
\]
is used on each rectangle; it follows from
\eqref{eq:Theta-PF4-sector-remainder}. All comparisons use directed
upper endpoints. Taylor's theorem now encloses the whole derivative
matrix $H(z)=[g^{(i+j)}(z)]_{i,j=0}^3$ on $|z|\le r$:
\[
 H_{ij}(z)=\sum_{n=0}^3 v_{i+j+n}\frac{z^n}{n!}+R_{ij}(z),
 \qquad |R_{ij}(z)|\le N_{i+j+4}\frac{r^4}{4!}.
\]

Each box carries an exact rational lower triangular matrix $C$ with
positive diagonal. Apply this congruence to the Taylor coefficients
before enclosing $z$. Writing $CH(z)C^T=Q(z)+E(z)$, we have
\[
 |E_{ij}(z)|\le F_{ij}:=\frac{r^4}{4!}
     \sum_{u\le i,\,v\le j}|C_{iu}C_{jv}|N_{u+v+4}.
\]
Every leading block is transformed by the corresponding leading block
of $C$, so its determinant changes by a positive square.
One acceptance test bounds $\|CH(z)C^T-J\|_\infty<1$ uniformly, where
$J=\operatorname{diag}(1,-1,1,-1)$; the homotopy from $J$ preserves
every leading determinant sign. The other test forms each
$\det Q_k(z)$ as a polynomial in the same variable $z$ and bounds its
signed value below, using $z^n\in[-r^n,r^n]$ for odd $n$ and
$z^n\in[0,r^n]$ for positive even $n$.
If $M_{ij}\ge\sup_{|z|\le r}|Q_{ij}(z)|$, the complete determinant
error is bounded by
\[
 |\det(Q_k+E_k)-\det Q_k|
 \le\sum_{\pi\in S_k}\ \sum_{\varnothing\ne I\subset\{0,\ldots,k-1\}}
    \prod_{i\in I}F_{i,\pi(i)}\prod_{i\notin I}M_{i,\pi(i)}.
\]
The signed polynomial lower bound must exceed this error for every
$k\le4$. Both tests prove signs throughout the closed physical interval
$[x-r\sigma,x+r\sigma]$, because the exponential tilt and positive
scale preserve these signs.
The records certify $89$ boxes on $[0.00211004,0.006]$ and $40$ on
$[0.006,0.01001]$; $54$ use the norm test and $75$ use the polynomial
test. Exact rational endpoint comparisons verify the closed cover.

\paragraph{The residue tail.}
Extract the same $64$ factors and write $p=f_{64}*\operatorname{Law}(V)$,
where $V$ includes every omitted factor and the positive anchored high
residual. The verified zeros satisfy $q_{64}<\tau:=29470<q_{65}$,
and the certified zero count below $\sqrt\tau$ is $64$.
Every remaining finite exponential factor therefore has a finite
moment at $\tau$. Lemma~\ref{lem:PF-residual-moment} gives the same
conclusion for the high residual, since $\tau<40000$.
Having established $\E e^{\tau V}<\infty$, analytic continuation of
the cancelled transform gives
\[
 \kappa:=\E e^{\tau V}-1
   =\frac{c_0}{\Xi(\sqrt\tau)}
       \prod_{\ell=1}^{64}(1-\tau/q_\ell)-1.
\]
The finite convolution and the actual residues are
\[
 f_{64}(x)=\boldsymbol1_{\{x>0\}}
          \sum_{\ell=1}^{64}c_\ell e^{-q_\ell x},
 \quad c_\ell=\frac{\prod_hq_h}{\prod_{h\ne\ell}(q_h-q_\ell)},
 \quad d_\ell=c_\ell\E e^{q_\ell V}
          =-\frac{2\gamma_\ell c_0}{\Xi'(\gamma_\ell)}.
\]
The zero and derivative enclosures verify $(-1)^{\ell-1}d_\ell>0$.
Set $G(x)=\sum_{\ell=1}^{64}d_\ell e^{-q_\ell x}$.
For polynomials $P,Q$ of degree at most three, the part of the
convolution with $V\ge x$ gives
\begin{equation}\label{eq:Theta-PF4-crossing-error}
 |P(\partial_x)Q(\partial_x)(p-G)(x)|
 \le\frac{\kappa e^{-\tau x}}{1-e^{-\tau x}}
       \sum_{\ell=1}^{64}|c_\ell P(-q_\ell)Q(-q_\ell)|.
\end{equation}
Indeed, for $v\ge x$ and $q_\ell<\tau$,
$e^{q_\ell v}/(e^{\tau v}-1)
 \le e^{(q_\ell-\tau)x}/(1-e^{-\tau x})$.
Integrate against $V$ and use $\E(e^{\tau V}-1)=\kappa$.
Differentiation through order six is valid because $f_{64}$ and these
derivatives match zero at the support boundary.

Use the normalized Newton polynomials
\[
 P_i(z)=\frac{\prod_{h=1}^i(z+q_h)}
        {\sqrt{|d_{i+1}|}\prod_{h=1}^i(q_{i+1}-q_h)},
 \qquad 0\le i\le3.
\]
Their coefficient matrix is lower triangular and invertible.
After this congruence and a positive diagonal scaling with entries
$e^{q_{i+1}x/2}$, the residue matrix has entries
\[
 A^0_{ij}(x)=\sum_{\ell>\max(i,j)}
       w_{\ell ij}e^{-\beta_{\ell ij}x},\qquad
 w_{\ell ij}=d_\ell P_i(-q_\ell)P_j(-q_\ell),\quad
 \beta_{\ell ij}=q_\ell-\frac{q_{i+1}+q_{j+1}}2.
\]
The leading diagonal term is $(-1)^i$; all other surviving exponents
are positive. The actual transformed matrix $A(x)$ differs from
$A^0(x)$ by at most
\[
 E_{ij}(L)=\frac{\kappa
    e^{-(\tau-(q_{i+1}+q_{j+1})/2)L}}{1-e^{-\tau L}}
       \sum_\ell |c_\ell P_i(-q_\ell)P_j(-q_\ell)|
 \qquad(x\ge L>0).
\]
At $L=0.032$, the absolute row sums of $A-J$, bounded by the nonleading
residue terms and these errors, are less than $0.880786$.
Every term decreases for larger $x$, so the same homotopy argument
proves all four signs on $[0.032,\infty)$.

For a closed rational box $[L,H]\subset[0.008,0.032]$, put
$x_0=(L+H)/2$ and $r=(H-L)/2$. The fourth-degree Taylor coefficients
and the full remainder bound are
\begin{align*}
 (C_n)_{ij}&=\frac1{n!}\sum_\ell
    w_{\ell ij}(-\beta_{\ell ij})^ne^{-\beta_{\ell ij}x_0},
       \qquad 0\le n\le4,\\
 F_{ij}&=\frac{r^5}{5!}\sum_\ell
       |w_{\ell ij}|\beta_{\ell ij}^5e^{-\beta_{\ell ij}L}
       +E_{ij}(L).
\end{align*}
Each box supplies an invertible rational lower triangular matrix $T$.
Transform the coefficients to $TC_nT^T$ and the remainder to
$|T|F|T|^T$ before inserting the parity intervals for $(x-x_0)^n$.
Interval elimination on this enclosure verifies successive pivot signs
$+,-,+,-$, checking each strict sign before division.
Lower triangularity transfers these signs to the leading minors of
$A$, and hence to those of $p$.
The resulting $1041$ closed rational boxes cover $[0.008,0.032]$
exactly, completing the third domain of \eqref{eq:Theta-PF4-domains}.

\paragraph{Certificate records.}
All numerical comparisons use outward Arb balls
\citep{johansson2017arb}; all interval joins and basis entries use exact
rationals.
The compact quadratures were produced at $40$ decimal digits;
the returned balls, including their radii, enter every later test.
A global check replays the endpoint inequalities and all $129$ compact
boxes at $100$ digits, independently regenerates the residue inputs and
all $1041$ tail boxes, and verifies the monotone tail and exact overlaps.
Separate implementations replayed both finite covers at $120$ digits,
the compact one from the certified quadrature balls. The smallest
acceptance margin of the compact replay exceeds $0.0147355$, and the
smallest signed pivot of the tail replay exceeds $0.0100751$; these
margins refer to the recorded congruences. The checks verify all three
overlapping domains in \eqref{eq:Theta-PF4-domains}, proving the lemma.
\end{proof}

The density $\Theta$ is smooth and flat at zero, so its zero extension
belongs to $C^6(\reals)$. Lemmas~\ref{lem:Theta-PF4-signs}
and~\ref{lem:PF-support-boundary} therefore prove the global order-four assertion
for $\Theta$, including its exact strictness condition.

\subsection{Low-zero enclosures for log-concavity}\label{app:zero-boxes}
For the low-height arguments, define $I_i=[\bar a_i-10^{-40},\bar a_i+10^{-40}]$,
$i=1,2$, by the following exact decimal centers:
\begingroup\small
\begin{center}
\begin{tabular}{rl}
\toprule
$i$&$\bar a_i$\\\midrule
1 & 14.13472514173469379045725198356247027078425711569924318\\
2 & 21.02203963877155499262847959389690277733434052490278175\\
\bottomrule
\end{tabular}
\end{center}
\endgroup
Interval evaluation gives opposite signs for $\Xi$ at the two
endpoints of every $I_i$, so continuity gives a real zero in each box.
The boxes are positive, disjoint, and increasing. The numbering and
simplicity of these zeros are not needed.

The evaluation uses Euler--Maclaurin with $N=1000$, $p=64$, and
$s=1/2+it$:
\begin{align*}
 \zeta(s)={}&\sum_{n=1}^{N-1}n^{-s}
 +\frac{N^{1-s}}{s-1}+\frac12N^{-s}
 +\sum_{k=1}^p\frac{B_{2k}}{(2k)!}(s)_{2k-1}N^{-s-2k+1}+R_p(s),\\
 |R_p(s)|\le{}&
 \frac{4|(s)_{2p}|}
 {(2\pi)^{2p}(\real s+2p-1)N^{\real s+2p-1}}.
\end{align*}
Here $(s)_j=s(s+1)\cdots(s+j-1)$ and $B_{2k}$ are exact rational
Bernoulli numbers. The remainder follows by integration of the periodic
Bernoulli remainder and
$\|B_{2p}(\{\cdot\})\|_\infty/(2p)!\le4/(2\pi)^{2p}$.
The resulting disc is enclosed by a complex rectangle and multiplied
by $s(s-1)\pi^{-s/2}\Gamma(s/2)/2$, with the gamma factor also in
interval arithmetic, to enclose $\Xi(t)$.

\subsection{The order-256 certificate}\label{app:PF256-certificate}
\label{app:PF10-certificate}
We describe the finite inequalities used in Lemma~\ref{lem:PF256-signs}.
All zero enclosures, signed sums, matrix operations, and remainders use
outward ball arithmetic \citep{johansson2017arb}. The $270$ selected
zeros are specified by their increasing indices in the certificate;
the largest index is $19930$, with ordinate approximately $17991.0263$.
Each component regenerates ordinate balls using the rigorous
\texttt{acb.zeta\_zero} routine and checks the positivity, ordering,
and analytic cutoffs needed by its method. The final independent
audit also checks the exact critical-line identity, disjointness,
and the bounds $a_l<\min(U/4,\Omega)$ for every selected index.

\paragraph{Signed matrices and Fourier norms.}
Put $S=100$, $r_l=a_l/S$, and, for $x\ge0$,
\[
 f(x)=\sum_{l=0}^{269}c_l e^{-a_lx},\qquad
 c_l=\frac{\prod_h a_h^2}{2a_l\prod_{h\ne l}(a_h^2-a_l^2)}.
\]
The signs of $c_l$ alternate, starting positive. Let
$\phi(t)=\widehat f(t)>0$. For even $j\le538$,
\[
 F_j=\frac1{2\pi}\int_\reals(t/S)^j\phi(t)\dd t
     =(-1)^{j/2}\sum_l c_lr_l^j>0.
\]
For a fixed lower triangular polynomial basis $P_i$ of degree $i$,
define, whenever $2i+q\le538$ and $q$ is even,
\begin{equation}\label{eq:PF256-norms}
 R_{i,q}^2=\frac1{2\pi}\int_\reals(t/S)^q
                    |P_i(it/S)|^2\phi(t)\dd t.
\end{equation}
These are finite signed combinations of the $F_j$; each is enclosed
and proved positive before taking its square root. Fourier inversion
and Cauchy--Schwarz bound the derivative of order $q$ of the
transformed entry by $S^qR_{i,q}R_{j,q}$. The same bound applies to
$p$, since $|\widehat p(t)|=\phi(t)\E e^{-Rt^2}\le\phi(t)$.

The matrix to be enclosed is
$[p^{(i+j)}(x)/S^{i+j}]_{i,j=0}^{255}$. A lower triangular
congruence with nonzero diagonal preserves the signs of every leading
determinant. The numerical preconditioners are selected from midpoint factorizations, but
their coefficients are fixed dyadic numbers. Signed sums and congruences
are completed before uncertain common coefficients are inserted.
Every rounding used to store a completed coefficient is outward.

\paragraph{The central interval.}
Choose an exact dyadic center $\mu_0$ of the mean enclosure and a
directed bound $\Delta\ge|\mu-\mu_0|$. Set
\[
 g_L=\sum_{l=0}^L\frac{\mu_0^l}{l!}D^{2l}f,\qquad
 \eta=\mu-\mu_0.
\]
Taylor's formula in the heat parameter and contraction of $H_t$ give
\[
 p=g_L+\eta D^2g_{L-1}+\mathcal R.
\]
In the transformed $(i,j)$ entry the remainder is bounded by
\begin{equation}\label{eq:PF256-central-error}
 \frac{(V_2+\Delta^2)S^4}{2}R_{i,4}R_{j,4}
 +\left(\frac{\mu_0^{L+1}}{(L+1)!}
       +\frac{\Delta\mu_0^L}{L!}\right)
       S^{2L+2}R_{i,2L+2}R_{j,2L+2}.
\end{equation}
Indeed, expanding $H_Rf$ about $H_{\mu_0}f$ leaves a second-order
error bounded by $\E(R-\mu_0)^2/2\le(V_2+\Delta^2)/2$;
expanding $H_{\mu_0}f$ and $D^2H_{\mu_0}f$ gives the other terms.
The heat parameters along the first Taylor segment are nonnegative.

On a closed box, expand the two reference functions about a fixed
center $z$ through odd spatial degree $K$. If
$h=\max_{x\text{ in the box}}|x-z|$, their spatial remainder is at most
\begin{align}\label{eq:PF256-central-spatial}
 \frac{h^{K+1}}{(K+1)!}\Bigg\{&
 \sum_{l=0}^L\frac{\mu_0^lS^{2l+K+1}}{l!}
      R_{i,2l+K+1}R_{j,2l+K+1}\notag\\
 &+\Delta\sum_{l=0}^{L-1}\frac{\mu_0^lS^{2l+K+3}}{l!}
      R_{i,2l+K+3}R_{j,2l+K+3}\Bigg\}.
\end{align}
The common variable $\eta$ and spatial powers are retained until after
the signed congruences. Odd and even powers have their correct
symmetric and nonnegative enclosures. Intersecting termwise and
Horner enclosures can tighten a polynomial without changing its range
guarantee. The factor condition
\[
 n\ge m+L+(K+1)/2
\]
supplies every required Fourier moment. For $n=270$, $m=256$, the
certificate uses $(L,K)=(10,7)$ on $[0,1/16000]$ and $(8,11)$ on
$[1/16000,1/1000]$. Both use moments only through degree $538$.
These choices of $(L,K)$ shift the derivative budget toward the spatial
remainder.

\paragraph{Shared cumulants outside the center.}
For $x>0$ write
$g(x)=\sum_l c_lM_le^{-a_lx}$, where $M_l=\E e^{a_l^2R}$.
Independent interval choices for the $M_l$ lose the cancellation
needed at high order. Choose dyadic centers $h_j^0$ of
\eqref{eq:PF256-cumulants}, for $1\le j\le N$, and put
\[
 L_0(s)=\sum_{j=1}^N h_j^0s^j/j,\qquad
 \E e^{sR}e^{-L_0(s)}=\sum_{n\ge0}d_ns^n.
\]
Here $d_0=1$ and
$nd_n=\sum_{j=1}^n(h_j-h_j^0)d_{n-j}$ for $1\le n\le N$.
For $\rho<U^2/16$, define
\begin{align*}
 H&=16^{N+1}/U^{2N+2}+D_0/U^{2N},\\
 G&=\sum_{j=1}^N|h_j-h_j^0|\rho^j/j
       +\frac{H\rho^{N+1}}{(N+1)(1-16\rho/U^2)}.
\end{align*}
The reciprocal-zero series bounds the omitted logarithmic
coefficients by this geometric majorant. Cauchy's estimate therefore
gives, for $0\le s<\rho$,
\begin{equation}\label{eq:PF256-cauchy}
 \left|\E e^{sR}e^{-L_0(s)}-\sum_{n=0}^Nd_ns^n\right|
 \le e^G\frac{(s/\rho)^{N+1}}{1-s/\rho}.
\end{equation}
For each common $d_n$, form the signed derivative matrix first, and
then multiply by its interval enclosure. The radius and every
convergence inequality are checked in the component record.
Conservative far-region specializations retain only $d_0,d_1$,
moving the higher terms into the absolute error, or use the
independent enclosures supplied by \eqref{eq:PF256-mgf}:
\[
 M_l\in\exp\left(a_l^2[\mu_-,\mu_+]
          +\left[0,\frac{a_l^4V_2}{2A(a_l^2)^2}\right]\right).
\]

For any fixed basis $P_i$, the branch-crossing error between $p$ and
$g$ is bounded entrywise by
\begin{equation}\label{eq:PF256-branch}
 \frac{2\kappa e^{-\Omega x}}{(1-e^{-\Omega x})^2}
 \sum_l|c_l|\bigl(|P_i(r_l)P_j(r_l)|+
                       |P_i(-r_l)P_j(-r_l)|\bigr).
\end{equation}
Only $Y>x$ changes the exponential branch. On that event,
$\cosh(\Omega Y)-1\ge e^{\Omega Y}(1-e^{-\Omega x})^2/2$;
together with $a_l<\Omega$ this proves the bound. The finite
convolution is smooth enough to differentiate before splitting the
expectation.

The outer bases start from the Newton polynomials
\[
 Q_i(w)=\frac{\prod_{h<i}(w+r_h)}
              {\sqrt{|c_i|}\prod_{h<i}(r_i-r_h)}.
\]
Their zeros $Q_i(-r_l)=0$ for $l<i$ are structural identities. The components either freeze the basis and
its positive normalization at the box center, keep the normalization
$\operatorname{diag}(e^{a_ix/2})$ moving, or apply a fixed positive
rescaling $\operatorname{diag}(e^{\lambda_i(x-z)})$ to a frozen basis.
All are valid triangular congruences throughout the box.
For the Jacobi components, an interval LDL correction gives the
orthogonal polynomial basis for the signed discrete weights
$c_l\exp(L_0(a_l^2)-a_lz)$; its verified three-term recurrence
evaluates the signed matrices and Fourier norms. These rate-dependent
bases are enclosed throughout; only the auxiliary midpoint
preconditioners and chosen rescaling rates are fixed dyadics.
For a term with exponent $\beta$, the degree-$K$
spatial error is at most its absolute coefficient times
\[
 \frac{h^{K+1}|\beta|^{K+1}}{(K+1)!}
 \max\{e^{-\beta(z-h)},e^{-\beta(z+h)}\}.
\]
Both signs of $\beta$ are covered; a structurally constant term has
zero remainder. Absolute errors are propagated through a further
congruence $C$ by $|C|E|C|^{\mathsf T}$.

\paragraph{Fourier remainder in the near region.}
On $[1/1000,3/1000]$, use $N=8$ shared cumulants and $K=25$.
For a frozen basis followed by the positive rescaling above, put
$q=K+1$ and $u_{ij}=\lambda_i+\lambda_j$. The spatial remainder of
the \emph{true} transformed entry is at most
\begin{equation}\label{eq:PF256-near-fourier}
 \frac{h^q e^{|u_{ij}|h}}{q!}
 \sum_{d=0}^{q/2}\binom{q/2}{d}|u_{ij}|^{q-2d}
                       S^{2d}R_{i,2d}R_{j,2d}.
\end{equation}
This follows from $|u+it|^q=(u^2+t^2)^{q/2}$,
$|\widehat p|\le\phi$, and Cauchy--Schwarz. The largest moment
index is $510+26=536\le538$. Taylor coefficients are computed from
the shared-cumulant reference; their Cauchy and branch errors are
also included. Summing their absolute Taylor errors costs at most
$e^{(a_l+|u_{ij}|)h}$. Evaluating the corresponding error bounds at
$z-h$ and multiplying by $e^{|u_{ij}|h}$ covers this cost, since
$a_l<\Omega$. Thus \eqref{eq:PF256-near-fourier} holds for the actual
density $p$. Every side interval left by a successful centered box is checked
separately.

\paragraph{Acceptance and the unbounded tail.}
A finite box is accepted only after strict leading determinant
signs have been proved. The implementations use alternating
nonzero pivots, a strict Neumann bound relative to
$J=\operatorname{diag}((-1)^i)$, or positive definiteness of the
even-index block together with negative definiteness of the odd-index
block. The first gives the signs by elimination. If $\|A-J\|_\infty<1$,
each leading-block homotopy from $J$ is nonsingular, giving the same
signs. For the parity test, a Schur complement gives
$\lceil k/2\rceil$ positive and $\lfloor k/2\rfloor$ negative
eigenvalues in every leading $k$-block, again the required signs.
Every inequality excludes zero; a failed enclosure is inconclusive.

For $x\ge3$, use the normalized Newton basis and compare the actual
matrix with $J_M=\operatorname{diag}((-1)^iM_i)$.
Since $M_i\ge1$, $\|J_M^{-1}\|_\infty\le1$.
Every nonleading reference term decays with a strictly positive
exponent. The branch error is bounded by
\[
 \frac{2\kappa W_iW_j
       e^{-[\Omega-(a_i+a_j)/2]x}}{(1-e^{-\Omega x})^2},\qquad
 W_i^2=\sum_l|c_l|\{Q_i(r_l)^2+Q_i(-r_l)^2\}.
\]
The certificate verifies all decay exponents and bounds every error
row sum at $x=3$ by less than $0.834$. These bounds decrease on the
whole ray. Each leading-block homotopy from $J_M$ therefore stays
nonsingular, proving all required determinant signs for $x\ge3$.

\paragraph{Coverage.}
The certificate uses the following
exact cover, where $x_* =491697171693/1638400000000$:
\begin{center}\small
\begin{tabular}{lrr}
\toprule
Range & Selected closed pieces & Method\\\midrule
$[0,1/1000]$ & 98 & Central heat expansion\\
$[1/1000,3/1000]$ & 81 & Near Fourier remainder\\
$[3/1000,3/100]$ & 71 & Jacobi basis\\
$[3/100,1/10]$ & 101 & Jacobi basis\\
$[1/10,1/5]$ & 160 & Jacobi basis\\
$[1/5,x_*]$ & 133 & Jacobi basis\\
$[x_*,1/2]$ & 17 & Newton basis\\
$[1/2,1]$ & 13 & Newton basis\\
$[1,3]$ & 9 & Newton basis\\\midrule
$[0,3]$ & 683 & Complete finite cover\\
$[3,\infty)$ & --- & Decreasing analytic bound\\
\bottomrule
\end{tabular}
\end{center}
An assembly step checks exact rational coverage and endpoints, a common
order and selected factor list, admissible parameters, and one strict
tail record. The cover uses $683$ of $684$ accepted intervals.

Setup and the main checks use $2000$ decimal digits, with $1000$
digits for the far-region matrix checks. The Jacobi components use
$N=8$, $K=49$; the far components use $K=17$, retaining either the
linear common correction or independent moment enclosures, with
every omitted term included in the error bound.
Central coefficients are stored in
outward-rounded $256$-bit balls only after their signed congruences.
The largest accepted central Neumann bound is less than $0.990$;
the tail row-sum bound is approximately $0.833815816692$.
A separate completed audit regenerated all $270$ selected zeros at
$180$ digits and checked their identities and cutoffs.

\end{document}